\documentclass[11pt,a4paper]{amsart} 

\usepackage[french,english]{babel}	
\usepackage[T1]{fontenc}		
\usepackage[utf8]{inputenc}		
\usepackage{enumerate}    
\usepackage{lmodern}			
\usepackage{dsfont}			
\usepackage{amsmath, amsthm, amssymb, amsthm,bbm,bm,mathrsfs,mathtools}	
\usepackage[pdfborder={0 0 0}]{hyperref}
\usepackage[toc,page]{appendix} 
\usepackage{epsfig,color}
\numberwithin{equation}{section}

\newtheorem{theorem}{Theorem}[section] 
\newtheorem{lemma}[theorem]{Lemma}
\newtheorem{proposition}[theorem]{Proposition}
\newtheorem{propdef}[theorem]{Proposition-Definition}

\newtheorem{corollary}[theorem]{Corollary}

\newtheorem{uniprob}[theorem]{Universal property}
\newtheorem{theo}{Theorem}
\newtheorem*{theo*}{Theorem}

\newcommand{\diff}{\,\mathrm{d}}
\newcommand{\Tr}{\operatorname{Tr}}
\newcommand{\tr}{\operatorname{tr}}
\newcommand{\Log}{\operatorname{Log}}

\newcommand{\End}{\operatorname{End}}
\newcommand{\D}{\operatorname{D}}
\newcommand{\Ad}{\operatorname{Ad}}

\newcommand{\Id}{\mathrm{Id}}
\newcommand{\id}{\mathrm{id}}
\newcommand{\hol}{\mathrm{hol}}

\begin{document}

\title{Free convolution operators and free Hall transform}
\author{Guillaume Cébron}
\thanks{\textit{E-mail adress:} guillaume.cebron@etu.upmc.fr}
\date{2013}

\maketitle
\vspace{-0.7cm}
\begin{abstract}We define an extension of the polynomial calculus on a $W^*$-probability space by introducing an algebra $\mathbb{C}\{X_i:i\in I\}$ which contains polynomials. This extension allows us to define transition operators for additive and multiplicative free convolution. It also permits us to characterize the free Segal-Bargmann transform and the free Hall transform introduced by Biane, in a manner which is closer to classical definitions. Finally, we use this extension of polynomial calculus to prove two asymptotic results on random matrices: the convergence for each fixed time, as $N$ tends to $\infty$, of the $*$-distribution of the Brownian motion on the linear group $GL_N(\mathbb{C})$ to the $*$-distribution of a free multiplicative circular Brownian motion, and the convergence of the classical Hall transform on $U(N)$ to the free Hall transform.

\end{abstract}
\vspace{-0.1cm}
\tableofcontents
\newpage
\section*{Introduction}
\addtocontents{toc}{\protect\setcounter{tocdepth}{-1}}

\subsection*{Conditional expectations}

Throughout this paper, $\left(\mathcal{A},\tau\right)$ will denote a $W^*$-probability space, two random variables $A,B\in \mathcal{A}$ will be said to be free if the von Neumann   algebras generated by $A$ and $B$ are free, and $\tau\left(\cdot|B\right)$ will denote the conditional expectation from $\mathcal{A}$ to the von Neumann algebra generated by $B$ (see Section~\ref{probaspace}).

In~\cite{Biane1998}, Biane shows how to compute some conditional expectations in free product of von Neumann algebras: let us recall Theorem 3.1 of~\cite{Biane1998}, which is in practice accompanied by a characterization of the Feller Markov kernel involved.
\begin{theo*}[Biane~\cite{Biane1998}]
Let $\left(\mathcal{A},\tau\right)$ be a $W^*$-probability space. Let $A,B\in \mathcal{A}$ be two self-adjoint random variables which are free. Then there exists a Feller Markov kernel $K = k(x,du)$ on $\mathbb{R}\times\mathbb{R}$ such that for any Borel bounded function $f$ on $\mathbb{R}$,
$$\tau \left(\left.f (A +B  )\right|B \right) = Kf (B )$$
(where $Kf (x) = \int_{\mathbb{R}} f (u)k(x,\diff u)$).
\end{theo*}
In this theorem, the Feller Markov kernel $K$ depends on both $A$ and $B\in \mathcal{A}$. Informally, we might expect that for any Borel bounded function $f$ on $\mathbb{R}$, there exists an object $R_Af$ which depends only on $A$ such that, for all $B\in\mathcal{A}$ free from $A$,
$$\tau \left(\left.f (A +B  )\right|B \right) = R_Af (B ).$$
Let us summarize how we partly answer this question in Section~\ref{condexp}.

We extend the notion of polynomial calculus to a more general calculus. Let $A\in \mathcal{A}$. We are more interested by random variables of the type
$\tau\left(P_1\left(A\right)\right)\cdots\tau\left(P_n\left(A\right)\right) \cdot P_0\left(A\right),$
with $n\in \mathbb{N}$, and $P_0,\ldots,P_n\in \mathbb{C}\left[X\right]$, than by random variables of the type $P\left(A\right)$ with $P\in \mathbb{C}\left[X\right]$. One could object that the random variables are exactly the same. However, the interest resides in the fact that, for all $n\in \mathbb{N},P_0,\ldots,P_n\in \mathbb{C}\left[X\right]$, the map
$$A\mapsto  \tau\left(P_1\left(A\right)\right)\cdots\tau\left(P_n\left(A\right)\right) \cdot P_0\left(A\right)$$
cannot be represented by a polynomial calculus nor by a functional calculus, but these sorts of maps are essential in the following development.

Let us denote by $\mathbb{C}\{X\}$ the vector space generated by the formal vectors
$$\left\{P_0 \tr\left(P_1\right)\cdots\tr\left(P_n\right):n\in \mathbb{N},P_0,\ldots,P_n\in \mathbb{C}\left[X\right]\right\}.$$
For $P=P_0 \tr\left(P_1\right)\cdots\tr\left(P_n\right)\in \mathbb{C}\{X\}$ and $A\in \mathcal{A}$, let us denote by $P(A)$ the random variable $P(A)=\tau\left(P_1\left(A\right)\right)\cdots\tau\left(P_n\left(A\right)\right) \cdot P_0\left(A\right)$ and we extend this notation to all $\mathbb{C}\{X\}$ by linearity. We construct the vector space $\mathbb{C}\{X\}$ and the $\mathbb{C}\{X\}$-calculus more precisely in Section~\ref{polcalcex}. In the same section, we construct similarly the space $\mathbb{C}\{X,X^{-1}\}$ (which naturally contains the space of Laurent polynomial $\mathbb{C}\left[X,X^{-1}\right]$) and the $\mathbb{C}\{X,X^{-1}\}$-calculus. More generally, for any index set $I$, we construct the space $\mathbb{C}\{X_i: i\in I\}$ and the $\mathbb{C}\{X_i:i\in I\}$-calculus.
We are now able to formulate the following theorem, which is a version of Theorem~\ref{freekernel} for one variable.
\begin{theo}
Let $A\in \mathcal{A}$.\label{freekernelone} There exists a linear operator $\Delta_A:\mathbb{C}\left\{X\right\}\rightarrow \mathbb{C}\left\{X\right\}$ such that, for all polynomials $P\in \mathbb{C}[X]$, and all $B\in \mathcal{A}$ free from $A$, we have
$$\tau\left(P\left(A+B\right)|B\right)=\left(e^{\Delta_{A}}P\right)\left(B\right).$$
\end{theo}
Theorem~\ref{freekernelone} (or Theorem~\ref{freekernel}) has two advantages. Firstly, it deals with arbitrary non-commutative random variables, and not only with self-adjoint variables. Secondly, it introduces a transition kernel $e^{\Delta_A}$ which depends only on the variable $A$ with which we convolved.

In~\cite{Biane1998}, Biane established multiplicative versions of his theorem. There is also a multiplicative version of Theorem~\ref{freekernelone} in Theorem~\ref{freekernelmult}. Let us formulate a version for one variable.
\begin{theo}
Let $A\in \mathcal{A}$. Then there exists a linear operator $\D_A:\mathbb{C}\left\{X\right\}\rightarrow \mathbb{C}\left\{X\right\}$ such that, for all polynomials $P\in \mathbb{C}[X]$, and all $B\in \mathcal{A}$ free from $A$, we have\label{freekernelmultone}
$$\tau\left(P\left(AB\right)|B\right)=\left(e^{\D_{A}}P\right)\left(B\right).$$
\end{theo}
Section~\ref{condexp} is devoted to the proofs of Theorems~\ref{freekernel} and~\ref{freekernelmult}, which are multivariate versions of Theorems~\ref{freekernelone} and~\ref{freekernelmultone}.

\subsection*{Free Hall transform}

In Section~\ref{fsb}, we use and extend Theorems~\ref{freekernel} and~\ref{freekernelmult} to give another description of the free Segal-Bargmann transform in Theorem~\ref{fsbt} and of the free Hall transform in Theorem~\ref{fht}. Let us explain the result for the Hall transform.

\subsubsection*{Classical Hall transform}


Let $N\in \mathbb{N}^*$. We endow the Lie algebra $\frak{u}(N)$ of the unitary group $U(N)$ with the inner product $\left\langle X,Y\right\rangle_{\frak{u}(N)}=N \Tr (X^*Y)$, and the Lie algebra $\frak{gl}_N(\mathbb{C})$ of the linear group $GL_N(\mathbb{C})$ with the real-valued inner product $\left\langle X,Y\right\rangle_{\frak{gl}_N(\mathbb{C})}=N \Re \Tr (X^*Y)$ (see Section~\ref{granma}). These scalar products determine right-invariant Laplace operators $\Delta_{U(N)}$ and $\Delta_{GL_N(\mathbb{C})}$ respectively.
The Brownian motion on $U(N)$ is a Markov process $(U^{(N)}_t)_{t\geq 0}$ on $U(N)$, starting at the identity, and with generator $\frac{1}{2}\Delta_{U(N)}$.
For all $t\geq 0$, we denote by $L^2(U^{(N)}_t)$ the Hilbert space
$$\left\{f\left(U^{(N)}_t\right): f \text{ is a complex Borel function on $U(N)$ such that }\mathbb{E}\left[\left|f\left(U^{(N)}_t\right)\right|^2\right]<+\infty\right\}.$$
The Brownian motion on $GL_N(\mathbb{C})$ is a Markov process $(G^{(N)}_t)_{t\geq 0}$ on $GL_N(\mathbb{C})$, with generator $\frac{1}{4}\Delta_{GL_N(\mathbb{C})}$, and starting at the identity. Observe that, for convenience, we have taken a definition of the Brownian motion on $GL_N(\mathbb{C})$ which includes an unusual factor of $2$ in its generator. In such a way, the Brownian motion proceeds at "half speed" on $GL_N(\mathbb{C})$, just as a standard complex Brownian motion on $\mathbb{C}$ is sometimes defined to be a Markov process with generator $\frac{1}{4}(\partial_x^2+\partial_y^2)$.
For all $t\geq 0$, we denote by $L^2_{\hol}(G^{(N)}_t)$ the Hilbert space
$$\left\{F\left(G^{(N)}_t\right): F \text{ is a holomorphic function on $GL_N(\mathbb{C})$ such that }\mathbb{E}\left[\left|F\left(G^{(N)}_t\right)\right|^2\right]<+\infty\right\}.$$
The fact that $L^2_{\hol}(G^{(N)}_t)$ is a Hilbert space is not trivial. It is a part of Hall's theorem which may be stated as follows (see~\cite{Driver1995},~\cite{Hall1994} and Section~\ref{hallclas}).
\begin{theo*}[Hall~\cite{Hall1994}]Let $t> 0$. Let $f$ be a Borel function on $U(N)$ such that $f(U^{(N)}_t)\in L^2(U^{(N)}_t)$. The function $e^{\frac{t}{2}\Delta_{U(N)}}f$ has an analytic continuation to a holomorphic function on $GL_N(\mathbb{C})$, also denoted by $e^{\frac{t}{2}\Delta_{U(N)}}f$. Moreover, $(e^{\frac{t}{2}\Delta_{U(N)}}f)(G^{(N)}_t)\in L^2_{\hol}(G^{(N)}_t)$ and the linear map
$$B_t:f\left(U^{(N)}_t\right)\mapsto \left(e^{\frac{t}{2}\Delta_{U(N)}}f\right)\left(G^{(N)}_t\right)$$
is an isomorphism of Hilbert spaces between $L^2(U^{(N)}_t)$ and $L^2_{\hol}(G^{(N)}_t)$. In particular, for all bounded Borel function $f$, we have $B_t(f(U^{(N)}_t))=F(G^{(N)}_t),$
where $F$ is the analytic continuation of $U\mapsto \mathbb{E}\left[f\left(U^{(N)}_t U\right)\right]$ to all $GL_N(\mathbb{C})$.
\end{theo*}
The transform $B_t$ is referred to as the Segal-Bargmann transform, the Segal-Bargmann-Hall transform, or the Hall transform. We choose to use the third name. It should be remarked that this formulation of Hall's theorem is quite different (but equivalent) from the original one: the point of view from which it is considered is a probabilistic one. Indeed, we identify the space of square integrable Borel functions with respect to the law of $U_t^{(N)}$ with $L^2(U^{(N)}_t)$, and the space of square integrable holomorphic functions with respect to the law of $G_t^{(N)}$ with $L^2_{\hol}(G^{(N)}_t)$. This identification between random variables and their functional representations is described in Section~\ref{rvfr}.

\subsubsection*{Free Hall transform}Let $\left(\mathcal{A},\tau\right)$ be a $W^*$-probability space. 
Let $\left(U_t\right)_{t\geq 0}$ be a free unitary Brownian motion, and let $\left(G_t\right)_{t\geq 0}$ be a free circular multiplicative Brownian motion (see~\cite{Biane1997}, or Section~\ref{fsc}). Let $t\geq 0$. We denote by $L^2(U_t,\tau)$ the Hilbert completion of the *-algebra generated by $U_t$ and $U_t^{-1}$ for the norm $\|\cdot\|_2:A\mapsto \tau\left(A^* A\right)^{1/2}$, and by $L^2_{\hol}(G_t,\tau)$ the Hilbert completion of the algebra generated by $G_t$ and $G_t^{-1}$ for the same norm $\|\cdot\|_2$.

In~\cite{Biane1997}, Biane defined the free Hall transform. This definition is based on a restriction of the free Segal-Bargmann transform in infinite dimensions, as Gross and Malliavin did in the classical case in~\cite{Gross1996}. In Theorem~\ref{fhtp}, a simplified version of Theorem~\ref{fht}, we give a description of the free Hall transform that is more direct and very close to the classical description thanks to the space $\mathbb{C}\{X,X^{-1}\}$ (see Section~\ref{hallfree} for more details about those two theorems).
\begin{theo*}[Biane~\cite{Biane1997}]
Let $t> 0$. There exists a linear transform $\mathcal{G}_t$ of the space of Laurent polynomials $\mathbb{C}\left[X,X^{-1}\right]$ such that $\mathcal{F}_t:P \Big(U_t \Big) \mapsto \mathcal{G}_t(P) \Big(G_t \Big)  $ is an isometric map which extends to a Hilbert space isomorphism $\mathcal{F}_t$ between $L^2(U_t,\tau)$ and $L^2_{\hol}(G_t,\tau)$.
\end{theo*}
\begin{theo}
Let $t> 0$.\label{fhtp}
For all $P\in \mathbb{C}\left[X\right]$, $\mathcal{F}_t(P (U_t ) )= (e^{\D_{U_t}}P) (G_t )  $, where $D_{U_t}$ is given in Theorem \ref{freekernelmultone}. Moreover, if $\left(U_t\right)_{t\geq 0}$ and $\left(G_t\right)_{t\geq 0}$ are free, for all $P\in \mathbb{C}\{X,X^{-1}\}$, $$\mathcal{F}_t\Big(P(U_t)\Big) =\tau\Big(P(U_t G_t)\Big|G_t\Big).$$
\end{theo}

\subsection*{Random matrices}
Section~\ref{granma} contains an application of our formalism and of our previous results to random matrices. The two main results are summarized in Theorem~\ref{main}: the large-$N$ limit for each fixed time of the noncommutative distribution of the Brownian motion $(G^{(N)}_t)_{t\geq 0}$ on $GL_N(\mathbb{C})$, and the large-$N$ limit of the Hall transform for $U(N)$. It should be mentioned that, at the same time as the author, Kemp has studied the first question in~\cite{Kemp2013} and~\cite{Kemp2013a}, and Driver, Hall and Kemp have studied the second question in~\cite{Driver2013}. While the approaches are quite similar, they do not entirely overlap.

In~\cite{Biane1997}, Biane suggests boosting $B_t$ to $B_t\otimes \Id_{M_N(\mathbb{C})}$ with the aim of studying the action of $B_t\otimes \Id_{M_N(\mathbb{C})}$ on random variables given by the functional calculus for matrices. Let us explain how the space $\mathbb{C}\{X,X^{-1}\}$ allows us to better understand the action of $B_t\otimes \Id_{M_N(\mathbb{C})}$ on variables given by the polynomial calculus. For all $N\in \mathbb{N}^*$, we endow $M_N(\mathbb{C})$ with the inner product $\langle X,Y\rangle_{M_N(\mathbb{C})}=\frac{1}{N} \Tr (X^*Y)$. Let us identify the space
$$\left\{f\left(U^{(N)}_t\right): f \text{ is a Borel function from $U(N)$ to $M_N(\mathbb{C})$ such that }\mathbb{E}\left[\left\|f\left(U^{(N)}_t\right)\right\|_{M_N(\mathbb{C})}^2\right]<+\infty\right\}$$
with $L^2(U^{(N)}_t)\otimes M_N(\mathbb{C})$ and the space
\begin{multline*}
\left\{\vphantom{\mathbb{E}\left[\left\|F\left(G^{(N)}_t\right)\right\|_{M_N(\mathbb{C})}^2\right]}F\left(G^{(N)}_t\right): F \text{ is a holomorphic function from $GL_N(\mathbb{C})$ to $M_N(\mathbb{C})$}\right.\\[-0.36cm]
\hspace{8.8cm}\left.\text{such that }\mathbb{E}\left[\left\|F\left(G^{(N)}_t\right)\right\|_{M_N(\mathbb{C})}^2\right]<+\infty\right\}
\end{multline*}
with $L^2_{\hol}(G^{(N)}_t)\otimes M_N(\mathbb{C})$. For all $t> 0$ and $N\in \mathbb{N}^*$, we denote by $B_t^{(N)}$ the Hilbert space isomorphism $B_t\otimes \Id_{M_N(\mathbb{C})}$ from $L^2(U^{(N)}_t)\otimes M_N(\mathbb{C})$ into $L^2_{\hol}(G^{(N)}_t)\otimes M_N(\mathbb{C})$.

Unfortunately, the space of random variables $\{P(U^{(N)}_t)\}_{ P\in \mathbb{C}[X,X^{-1}]}$ is not transformed by $B_t^{(N)}$ into the space $\{P(G^{(N)}_t)\}_{ P\in \mathbb{C}[X,X^{-1}]}$. The $\mathbb{C}\{X,X^{-1}\}$-calculus offers us larger spaces which are stable under $B_t^{(N)}$. Indeed,  the space of random variables $\{P(U^{(N)}_t)\}_{ P\in \mathbb{C}\{X,X^{-1}\}}$ is transformed by $B_t^{(N)}$ into the space of random variables $\{P(G^{(N)}_t)\}_{ P\in \mathbb{C}\{X,X^{-1}\}}$ (Proposition~\ref{btpun}).

The use of the $\mathbb{C}\{X,X^{-1}\}$-calculus also allows us to study the limit in large dimension. It is already known that the free unitary Brownian motion is the limit in distribution of the Brownian motion on $U(N)$ (see~\cite{Biane1997a}, \cite{LEVY2008}, \cite{Rains1997} and \cite{Sengupta2008}), which is the first item of the following theorem. For the two other items, see Theorem~\ref{Glim} and Theorem~\ref{Flim}. As mentioned above, the concurrent papers~\cite{Driver2013} and~\cite{Kemp2013,Kemp2013a} address respectively the third item and the second item with complementary techniques and points of view.
\begin{theo}Let $t\geq 0$.\label{main}
\begin{enumerate}
\item For all $n\in \mathbb{N}$, and all Laurent polynomials $P_0,\ldots,P_n\in \mathbb{C}[X,X^{-1}]$, we have
$$\displaystyle\lim_{N\rightarrow \infty}\mathbb{E}\left[ \tr\left(P_0\left(U_t^{(N)}\right)\right)\cdots \tr\left(P_n\left(U_t^{(N)}\right)\right) \right]=\tau\left(P_0\left(U_t\right)\right)\cdots\tau\left(P_n\left(U_t\right)\right) .$$
\item For all $n\in \mathbb{N}$, and all polynomials $P_0,\ldots,P_n\in \mathbb{C}\langle X,X^*,X^{-1},{X^*}^{-1}\rangle$, we have
$$\displaystyle\lim_{N\rightarrow \infty}\mathbb{E}\left[ \tr\left(P_0\left(G_t^{(N)}\right)\right)\cdots \tr\left(P_n\left(G_t^{(N)}\right)\right) \right]=\tau\left(P_0\left(G_t\right)\right)\cdots\tau\left(P_n\left(G_t\right)\right) .$$
\item For all Laurent polynomial $P\in \mathbb{C}[X,X^{-1}]$, and $t>0$, as $N \rightarrow \infty$, we have
$$ \left\|B_t^{(N)}\Big(P\left(U_t^{(N)}\right)\Big)-\mathcal{G}_t(P)\Big(G_t^{(N)} \Big)\right\|^2_{L^2_{\hol}(G^{(N)}_t)\otimes M_N(\mathbb{C})}=O(1/N^2).$$
\end{enumerate}
\end{theo}
In fact, the second item can be strengthened to the convergence of the full process $(G^{(N)}_t)_{t\geq 0}$ to $\left(G_t\right)_{t\geq 0}$ (i.e. the convergence of all finite dimensional distributions).
Indeed, since the increments of $\left(G_t\right)_{t\geq 0}$ are free, it remains to prove that the increments $(G^{(N)}_t)_{t\geq 0}$ are asymptotically free, which is a consequence of the $\Ad(U(N))$-invariance of the law of each increment of $(G^{(N)}_t)_{t\geq 0}$. A complete proof is present in~\cite{Kemp2013a}.

\addtocontents{toc}{\protect\setcounter{tocdepth}{2}}

\section{Functional calculus extension}
\label{polcalcex}
In this section, we define the algebra $\mathbb{C}\{X_i:i\in I\}$. We denote by $\mathbb{C}\left\langle X_i:i\in I\right\rangle$ the space of polynomials in the non commuting indeterminates $(X_i)_{i\in I}$. The algebra $\mathbb{C}\{X_i:i\in I\}$ is an extension of the algebra $\mathbb{C}\left\langle X_i:i\in I\right\rangle$. Intuitively, $\mathbb{C}\{X_i:i\in I\}$ is the free algebra generated by $I$ indeterminates $(X_i)_{i\in I}$ and an indeterminate center-valued expectation $\tr$.

We also define a $\mathbb{C}\{X_i:i\in I\}$-calculus which extends the polynomial calculus, but which depends also on the data of another center-valued trace $\tau$. This algebra $\mathbb{C}\{X_i:i\in I\}$ and its functional calculus is the basis of all others sections.

\subsection{The algebra $\mathbb{C}\{X_i:i\in I\}$}The algebra $\mathbb{C}\{X_i:i\in I\}$ will be defined by a universal property, which allows us to forget about its construction and to focus on its properties.

We present first the universal property and its immediate consequences before introducing the algebra $\mathbb{C}\{X_i:i\in I\}$ as its unique solution.

\subsubsection{Universal property}
Let $\mathcal{A}$ be a unital complex algebra. The center of $\mathcal{A}$ is the unital complex algebra $Z_\mathcal{A}$ formed by elements of $\mathcal{A}$ which commute with all the elements in $\mathcal{A}$.
 The algebra $\mathcal{A}$ is then a $Z_\mathcal{A}$-module. 
A center-valued expectation $\tau$ is 
a linear function from $\mathcal{A}$ to $Z_\mathcal{A}$ such that
\begin{enumerate}
\item for all $A,B\in \mathcal{A}$, we have $\tau(\tau(A)B)=\tau( A)\tau(B)$;
\item $\tau(1_\mathcal{A})=1_\mathcal{A}$.
\end{enumerate}
Let us remark that the restriction for $\tau$ to be a morphism of $Z_\mathcal{A}$-modules is not required, since it is not needed in this paper.
\begin{uniprob}\label{up}
Let $I$ be an arbitrary index set. Let $\mathcal{X}$ be an algebra endowed with a center-valued expectation $\tr$, and with $I$ specified elements $\left(X_i\right)_{i\in I}$. The triplet $\left(\mathcal{X},\tr,\left(X_i\right)_{i\in I}\right)$ possesses the universal property~\ref{up} for index set $I$ if for all algebras $\mathcal{A}$ endowed with a center-valued expectation $\tau$, and with $I$ elements $\left(A_i\right)_{i\in I}$, there exists a unique algebra homomorphism $f$ from $\mathcal{X}$ to $\mathcal{A}$ such that
\begin{enumerate}
\item for all $i\in I$, we have $f(X_i)=A_i$;
\item for all $X\in \mathcal{X}$, we have $\tau(f(X))=f(\tr(X))$.
\end{enumerate}
\end{uniprob}

Such a homomorphism will be called an $I$-adapted homomorphism from $\left(\mathcal{X},\tr,\left(X_i\right)_{i\in I}\right)$ to $\left(\mathcal{A},\tau,\left(A_i\right)_{i\in I}\right)$, or more simply from $\mathcal{X}$ to $\mathcal{A}$. If an $I$-adapted homomorphism is bijective, we will call it an $I$-adapted isomorphism.

The nature of the universal property~\ref{up} induces some properties on its solutions summed up in the following proposition.
\begin{proposition}Let $I$ be an arbitrary index set, and let $\left(\mathcal{X},\tr,\left(X_i\right)_{i\in I}\right)$ possess the universal property~\ref{up} for $I$.
\begin{enumerate}
\item The $I$-adapted isomorphim $\id_{\mathcal{X}}$ is the unique $I$-adapted algebra automorphism on $\mathcal{X}$.
\item $\left(\mathcal{X},\tr,\left(X_i\right)_{i\in I}\right)$ is unique in the following sense: if $\left(\mathcal{Y},\widetilde{\tr},\left(Y_i\right)_{i\in I}\right)$ also possesses the universal property~\ref{up} for $I$, there exists a unique $I$-adapted isomorphism from $\mathcal{X}$ to $\mathcal{Y}$.
\item There exists a unique algebra homomorphism $f$, which is injective, from $\mathbb{C}\left\langle Y_i:i\in I\right\rangle$ to $\mathcal{X}$ such that $f(Y_i)=X_i$ for all $i\in I$.
\end{enumerate}\label{UniProp}
\end{proposition}

Thus, if $I$ is an arbitrary index set, and $\left(\mathcal{X},\tr,\left(X_i\right)_{i\in I}\right)$ possesses the universal property~\ref{up} for $I$, we will always see $\mathbb{C}\left\langle X_i:i\in I\right\rangle$ as a subalgebra of $\mathcal{X}$.
\begin{proof} The first assertion is immediate.

Let $\left(\mathcal{Y},\widetilde{\tr},\left(Y_i\right)_{i\in I}\right)$ also possess the universal property~\ref{up} for $I$. There exists a unique $I$-adapted homomorphism $f$ from $\mathcal{X}$ to $\mathcal{Y}$ and an $I$-adapted homomorphism $g$ from $\mathcal{Y}$ to $\mathcal{X}$. Then $g\circ f$ is an $I$-adapted homomorphism from $\mathcal{X}$ to $\mathcal{X}$, and consequently, $g\circ f=\id_{\mathcal{X}}$. Similarly, $f\circ g=\id_{\mathcal{Y}}$, and therefore, $f$ is an $I$-adapted isomorphism from $\mathcal{X}$ to $\mathcal{Y}$.

For the third assertion, we endow $\mathbb{C}\left\langle Y_i:i\in I\right\rangle$ with the center-valued expectation $\tau$ such that $\tau(M)=1$ for all monomials $M\in \mathbb{C}\left\langle Y_i:i\in I\right\rangle$. There exists a unique algebra homomorphism $f$ from $\mathbb{C}\left\langle Y_i:i\in I\right\rangle$ to $\mathcal{X}$ such that $f(Y_i)=X_i$ for all $i\in I$. There exists an $I$-adapted homomorphism $g$ from $\mathcal{X}$ to $\mathbb{C}\left\langle Y_i:i\in I\right\rangle$. Finally, $g\circ f$ is an algebra automorphism on $\mathbb{C}\left\langle Y_i:i\in I\right\rangle$ such that $g\circ f(Y_i)=Y_i$ for all $i\in I$, and thus it is equal to $\id_{\mathbb{C}\left\langle Y_i:i\in I\right\rangle}$. We deduce that $f$ is injective.\end{proof}

\subsubsection{Definition}
We now state an existence result in the following proposition.
\begin{propdef}\label{ni}
Let $I$ be an arbitrary index set. There exists an object satisfying the universal property~\ref{up} for $I$. This unique (up to an $I$-adapted isomorphism) object will be denoted by $\left(\mathbb{C}\{X_i:i\in I\},\tr,\left(X_i\right)_{i\in I}\right)$. Furthermore,
$$\{M_0 \tr M_1 \cdots \tr  M_n: n\in \mathbb{N}, M_0,\ldots, M_n \text{ are monomials of }\mathbb{C}\left\langle X_i:i\in I\right\rangle\}$$
is a basis of $\mathbb{C}\{X_i:i\in I\}$, called the canonical basis.
\end{propdef}

Let us recall that $\mathbb{C}\left\langle X_i:i\in I\right\rangle$ is viewed as a subalgebra of $\mathbb{C}\{X_i:i\in I\}$. Thus, the set
$$\{M_0 \tr M_1 \cdots \tr  M_n: n\in \mathbb{N}, M_0,\ldots M_n, \text{ are monomials of }\mathbb{C}\left\langle X_i:i\in I\right\rangle\}$$
is unambiguously defined in $\mathbb{C}\{X_i:i\in I\}$. Furthermore, an $I$-adapted isomorphism does not modify the definition of this set. Thus, Proposition-Definition~\ref{ni} tells us that this set forms a basis in one particular realization and consequently, forms a basis in any realization of the universal property~\ref{up}.
\begin{proof}
We will present a construction of $\left(\mathbb{C}\{X_i:i\in I\},\tr,\left(X_i\right)_{i\in I}\right)$ in the appendix, which satisfies the characterization of the canonical basis, and a proof of its universal property.
\end{proof}
\begin{proposition}Let $J \subset I$ be two arbitrary index sets. There exists a unique $J$-adapted morphism, which is injective, from $\left(\mathbb{C}\{X_i:i\in J\},\tr,\left(X_i\right)_{i\in J}\right)$ to $\left(\mathbb{C}\{X_i:i\in I\},\tr,\left(X_i\right)_{i\in J}\right)$.
\end{proposition}

If $J \subset I$ are two arbitrary index sets, we will always see $\mathbb{C}\{X_i:i\in J\}$ as a subalgebra of $\mathbb{C}\{X_i:i\in I\}$.
\begin{proof}
Let $\left(Y_i\right)_{i\in I}$ be elements of $\mathbb{C}\{X_i:i\in J\}$ such that $Y_i=X_i$ for all $i\in J$ (one can set $Y_i=1_{\mathbb{C}\{X_i:i\in J\}}$ for $i\notin J$).

There exists a unique $J$-adapted morphism $f$ from $\mathbb{C}\{X_i:i\in J\}$ to $\mathbb{C}\{X_i:i\in I\}$. There exists an $I$-adapted homomorphism $g$ from $\left(\mathbb{C}\{X_i:i\in I\},\tr,\left(X_i\right)_{i\in I}\right)$ to $\left(\mathbb{C}\{X_i:i\in J\},\tr,\left(Y_i\right)_{i\in I}\right)$. Finally, $g\circ f$ is a $J$-adapted automorphism on $\mathbb{C}\{X_i:i\in J\}$, and thus $g\circ f$ is equal to $\id_{\mathbb{C}\{X_i:i\in J\}}$. We deduce that $f$ is injective.

\end{proof}

\subsubsection{Degrees}

Let us define a notion of degree on $\mathbb{C}\{X_i:i\in I\}$.

Let $n\in \mathbb{N}$, and $M_0,\ldots,M_n \in  \mathbb{C}\left\langle X_i:i\in I\right\rangle$ be monomials whose degrees are respectively $k_0,\ldots,k_n \in \mathbb{N}$. The degree of $M_0\tr \left( M_1\right) \cdots \tr \left( M_n\right)$ is defined to be $k_0+\ldots+k_n $. For all $P\in \mathbb{C}\{X_i:i\in I\}$, the degree of $P$ is defined to be the maximal degree of the elements of its decomposition in the canonical basis.

For $d \in \mathbb{N}$, we denote by $\mathbb{C}_d\{X_i:i\in I\}$ the subspace of $\mathbb{C}\{X_i:i\in I\}$ whose elements have degrees less than $d$. We have $\mathbb{C}\{X_i:i\in I\}=\displaystyle\cup_{d=0}^{\infty}\mathbb{C}_d\{X_i:i\in I\}$. If $I$ is finite, each space $\mathbb{C}_d\{X_i:i\in I\}$ is a finite-dimensional space. In particular, the space $\mathbb{C}\{X_i:i\in I\}$ is the union of finite-dimensional spaces: $$\mathbb{C}\{X_i:i\in I\}=\bigcup_{\substack{d\in \mathbb{N}\\ J\subset I,\ J\ \mathrm{finite}}}\mathbb{C}_d\{X_i:i\in J\}.$$

\subsubsection{Another product}
Let us define on $\mathbb{C}\{X_i:i\in I\}$ a second product which is bilinear and associative.

For all $P, Q\in \mathbb{C}\{X_i:i\in I\}$, the product $P\cdot_{\tr} Q$ is defined by
$P \cdot_{\tr} Q=P \tr Q$. The bilinearity is due to the linearity of $\tr$ and the associativity is simply due to the fact that, for all $P,Q$ and $R\in \mathbb{C}\left\langle X_i:i\in I\right\rangle$, we have
$$P \cdot_{\tr} (Q \cdot_{\tr} R)=P\tr(Q\tr(R))=P\tr Q \tr R=(P\tr Q) \tr R=(P \cdot_{\tr} Q) \cdot_{\tr} R.$$
Thus, $\left(\mathbb{C}\{X_i:i\in I\}, \cdot_{\tr}  \right)$ is a unital complex algebra. Moreover, this algebra is generated by the monomials of $ \mathbb{C}\left\langle X_i:i\in I\right\rangle$. Indeed, for all $n\in \mathbb{N}$, $P_0,\ldots,P_n\in \mathbb{C}\{X_i:i\in I\}$, we have\label{Anotprod}
$$P_0\tr P_1 \cdots \tr  P_n=P_0\cdot_{\tr} P_1\cdot_{\tr} \cdots\cdot_{\tr} P_n.$$

\subsection{The $\mathbb{C}\{X_i:i\in I\}$-calculus}

Let $\mathcal{A}$ be a unital complex algebra, and $\tau$ be a linear functional on $\mathcal{A}$ such that $\tau\left(1_{\mathcal{A}}\right)=1$.

We define a $\mathbb{C}\{X_i:i\in I\}$-calculus on $\left(\mathcal{A},\tau \right)$ in this way. Let $\mathbf{A}=\left(A_i\right)_{i\in I}$ be a family of elements in $\mathcal{A}$.  Viewing $\mathbb{C}$ as contained in $\mathcal{A}$ as $\mathbb{C}\cdot 1_{\mathcal{A}}$, the algebra $\mathcal{A}$ is endowed with a center-valued expectation $\tau$, and with $I$ elements $\left(A_i\right)_{i\in I}$. Thus, there exists a unique algebra homomorphism $f$ from $\mathbb{C}\{X_i:i\in I\}$ to $\mathcal{A}$ such that
\begin{enumerate}
\item for all $i\in I$, we have $f(X_i)=A_i$;
\item for all $X\in \mathbb{C}\{X_i:i\in I\}$, we have $\tau (f(X))=f(\tr(X))$.
\end{enumerate}
For $P\in \mathbb{C}\{X_i:i\in I\}$, we say that $f(P)$ is the element of $\mathcal{A}$ obtained by substitution of $\left(A_i\right)_{i\in I}$ for the indeterminates $\left(X_i\right)_{i\in I}$ and the substitution of $\tau$ for $\tr$ in $P$, and we denote this element by $P(A_i:i\in I)=f(P)$. Thus, the map $P\mapsto P(A_i:i\in I)$ is the unique algebra homomorphism from $\left(\mathbb{C}\{X_i:i\in I\},\cdot \right)$ to $\mathcal{A}$ such that
\begin{enumerate}
\item for all $j\in I$, we have $X_j(A_i:i\in I)=A_j$;
\item for all $P\in \mathbb{C}\{X_i:i\in I\}$, we have $\tau (P(A_i:i\in I))=(\tr P)(A_i:i\in I)$.
\end{enumerate}

\subsubsection{The algebra $\mathbb{C}\{X_i,X_i^*:i\in I\}$}From an arbitrary index set $I$, we construct the index set $\tilde{I}=I\cup (I\times \{*\})$. For all $i\in I$, the element $X_{(i,*)}\in\mathbb{C}\{X_i:i\in \tilde{I}\}$ will be denoted by $X^*_i$, and the algebra $\mathbb{C}\{X_i:i\in \tilde{I}\}$ will be denoted $\mathbb{C}\{X_i,X_i^*:i\in I\}$.
We define a *-algebra structure on $\mathbb{C}\{X_i,X_i^*:i\in I\}$ with the involution * given naturally by taking for all $i\in I$, $(X_i)^*=X_i^*$, and for all $P \in \mathbb{C}\{X_i,X_i^*:i\in I\}$, $(\tr P)^*=\tr (P^*)$.\label{involution}

Let $\mathcal{A}$ be a unital complex *-algebra, and $\tau$ be a linear functional on $\mathcal{A}$ such that $\tau\left(1_{\mathcal{A}}\right)=1$. Let us define a $\mathbb{C}\{X_i,X_i^*:i\in I\}$-calculus on $\left(\mathcal{A},\tau \right)$ in this way. Let $\left(A_i\right)_{i\in I}\in \mathcal{A}^{I}$. For all $i\in I$, set $A_{(i,*)}=A_i^*$. For $P\in \mathbb{C}\{X_i,X_i^*:i\in I\}$, let us denote $P(A_i:i\in \tilde{I})$ by $P(A_i:i\in I)$. Thus, the map $P\mapsto P(A_i:i\in I)$ is a *-algebra homomorphism from $\left(\mathbb{C}\{X_i,X_i^*:i\in I\},\cdot \right)$ to $\mathcal{A}$.

\subsubsection{The algebra $\mathbb{C}\{X_i,X_i^*,X_i^{-1},{X^*_i}^{-1}:i\in I\}$}
Similarly, from an arbitrary index set $I$, we construct the index set $\tilde{I}=I\cup (I\times \{*\})\cup (I\times \{-1\})\cup(I\times \{-*\})$. For all $i\in I$, the element $X_{(i,*)}\in\mathbb{C}\{X_i:i\in \tilde{I}\}$ will be denoted by $X^*_i$, the element $X_{(i,-1)}\in\mathbb{C}\{X_i:i\in \tilde{I}\}$ will be denoted by $X^{-1}_i$, and the element $X_{(i,-*)}\in\mathbb{C}\{X_i:i\in \tilde{I}\}$ will be denoted by ${X^*_i}^{-1}$. Finally, the algebra $\mathbb{C}\{X_i:i\in \tilde{I}\}$ will be denoted by $\mathbb{C}\{X_i,X_i^*,X_i^{-1},{X^*_i}^{-1}:i\in I\}$.
We define a *-algebra structure on $\mathbb{C}\{X_i,X_i^*,X_i^{-1},{X^*_i}^{-1}:i\in I\}$ with the involution * given naturally by taking for all $i\in I$, $(X_i)^*=X_i^*$ and $(X_i^{-1})^*={X_i^*}^{-1}$, and for all $P \in\mathbb{C}\{X_i,X_i^*,X_i^{-1},{X^*_i}^{-1}:i\in I\}$, $\left(\tr P\right)^*=\tr (P^*)$.

Let $\mathcal{A}$ be a unital complex *-algebra, and $\tau$ be a linear functional on $\mathcal{A}$ such that $\tau\left(1_{\mathcal{A}}\right)=1$. Let $P\in\mathbb{C}\{X_i,X_i^*,X_i^{-1},{X^*_i}^{-1}:i\in I\}$, and $\left(A_i\right)_{i\in I}\in \mathcal{A}^{I}$ be a family of invertible elements. For all $i\in I$,  set $A_{(i,*)}=A_i^*$, $A_{(i,-1)}=A_i^{-1}$ and $A_{(i,-*)}={A^*_i}^{-1}$. For $P\in \mathbb{C}\{X_i,X_i^*:i\in I\}$, let us abuse notation slightly and denote $P(A_i:i\in \tilde{I})$ also by $P(A_i:i\in I)$, as this should cause no confusion. Thus, the map $P\mapsto P(A_i:i\in I)$ is a *-algebra homomorphism from $\mathbb{C}\{X_i,X_i^*,X_i^{-1},{X^*_i}^{-1}:i\in I\}$ to $\mathcal{A}$.

\subsubsection{Factorization by the distribution}There exists a useful factorization of the $\mathbb{C}\{X_i:i\in I\}$-calculus by the $\mathbb{C}\left\langle X_i:i\in I\right\rangle$-calculus. Let us explain how it works.\label{fact}

Let $\mathcal{A}$ be a unital complex algebra, and $\tau$ be a linear functional on $\mathcal{A}$ such that $\tau\left(1_{\mathcal{A}}\right)=1$ and $\tau\geq 0$ (i.e. $\tau(a a^*)\geq 0$ for all $a\in \mathcal{A}$).
Such a space is called a non-commutative probability space. Elements of $\mathcal{A}$ are called (non-commutative) random variables. Let $\mathbf{A}=\left(A_i\right)_{i\in I}$ be a family of non-commutative variables of $\mathcal{A}$. The map$$\begin{array}{crcl}\mu_{\mathbf{A}}:&\mathbb{C}\left\langle X_i:i\in I\right\rangle & \rightarrow &\mathbb{C}  \\&P & \mapsto & \tau(P(\mathbf{A}))\end{array}$$ will be called the distribution of $\mathbf{A}$. The algebra $\mathbb{C}\left\langle X_i:i\in I\right\rangle$ is then endowed with a center-valued expectation $\mu_{\mathbf{A}}$, and with $I$ specified elements $\left(X_i\right)_{i\in I}$. Thus, there exists a unique algebra homomorphism $f$ from $\mathbb{C}\{X_i:i\in I\}$ to $\mathbb{C}\left\langle X_i:i\in I\right\rangle$ such that
\begin{enumerate}
\item for all $i\in I$, we have $f(X_i)=X_i$;
\item for all $X\in \mathbb{C}\{X_i:i\in I\}$, we have $\mu_{\mathbf{A}} (f(X))=f(\tr(X))$.
\end{enumerate}
For $P\in \mathbb{C}\{X_i:i\in I\}$, we say that $f(P)$ is the element of $\mathbb{C}\left\langle X_i:i\in I\right\rangle$ obtained by substitution of $\mu_\mathbf{A}$ for $\tr$ in $P$, and we denote this element by $\left.P\right|_\mathbf{A}$.

However, we can now use polynomial calculus, by substitution of $\left(A_i\right)_{i\in I}$ for $\left(X_i\right)_{i\in I}$ in $\left.P\right|_\mathbf{A}$. Because the homomorphism from $\mathbb{C}\{X_i:i\in I\}$ to $\mathcal{A}$ given by $P\mapsto \left.P\right|_\mathbf{A}(\mathbf{A})$ is an $I$-adapted homomorphism, we have that $\left.P\right|_\mathbf{A}(\mathbf{A})=P(\mathbf{A})$ using the universal property of $\mathbb{C}\{X_i:i\in I\}$.

\section{Computation of some conditional expectations}\label{condexp}

In this section, we show the existence of operators on $\mathbb{C}\{X_i:i\in I\}$ which play the role of transition kernels in the context of free convolution. The first result, Theorem~\ref{freekernel}, deals with additive free convolution whereas the second one, Theorem~\ref{freekernelmult}, deals with multiplicative free convolution. Despite the analogy of the two theorems, the proofs are completely different. This is to be expected, since the non-commutativity means that the direct connection $e^{x+y}=e^x\cdot e^y$ between addition and multiplication is lost.

\subsection{Generalities}\label{probaspace}

\subsubsection{$W^*$-probability spaces}Let $\left(\mathcal{A},\tau\right)$ be a non-commutative probability space such that $\mathcal{A}$ is a von Neumann algebra, and $\tau$ is a faithful normal tracial state. That is to say $\tau$ is a linear functional such that $\tau(1_\mathcal{A})=1$, and
\begin{enumerate}
 \item for all $A\in \mathcal{A}$, if $A\geq 0$, then $ \tau(A)\geq 0$ (positivity),
   \item for all $A,B\in \mathcal{A}$, $\tau(AB)=\tau(BA)$ (traciality),
   \item $\tau$ is continuous for the ultraweak topology (normality),
   \item for all $A\in \mathcal{A}$, if $\tau(A^* A)=0$, then $A=0$ (faithfulness).
\end{enumerate}
We call $\left(\mathcal{A}, \tau \right)$ a $W^*$-probability space. For all $\mathbf{A}=(A_i)_{i\in I}\in \mathcal{A}^{I}$,  we denote by $W^*(\mathbf{A})$ the von Neumann subalgebra of $\mathcal{A}$ generated by $(A_i)_{i\in I}$.

\subsubsection{Freeness}Let $I$ be a set of indices. For all $i \in I$, let $\mathcal{B}_i$ be a von Neumann subalgebra of $\mathcal{A}$. These algebras are called free if, for all $n\in \mathbb{N}$, and all indices $i_1 \neq i_2 \neq \ldots \neq i_n$, whenever $A_j \in \mathcal{B}_{i_j}$ and $\tau(A_j) = 0$ for all $1\leq j \leq n$, we have $\tau(A_1\cdots A_n) = 0$.

For all $(A_i)_{i\in I}\in \mathcal{A}^I$, we say that the elements $A_i$ are free for $i\in I$ if the algebras $W^*(A_i)$ are free for $i\in I$. For example, the operators $1_\mathcal{A}$ and $0_\mathcal{A}$ are free from all $A\in \mathcal{A}$.

\subsubsection{Conditional expectation}If $\mathcal{B} \subset \mathcal{A}$ is a von Neumann subalgebra, there exists a unique conditional expectation from $\mathcal{A}$ to $\mathcal{B}$ with respect to $\tau$, which we denote by $\tau\left(.|\mathcal{B}\right)$. This map is a weakly continuous, completely positive, identity preserving, contraction, and it is characterized by the property that, for any $A \in \mathcal{A}$ and $B \in \mathcal{B}$, $\tau(AB) = \tau(\tau(A|\mathcal{B})B)$
(see e.g.~\cite{Accardi1982} and~\cite{Takesaki1972}). For any index set $I$, and $\mathbf{A}=(A_i)_{i\in I}\in \mathcal{A}^I$, we denote by $\tau\left(\cdot|\mathbf{A}\right)$ the conditional expectation $\tau(\cdot|W^*(\mathbf{A}))$.

\subsection{Free cumulants}This section contains a succinct presentation of the theory of the free cumulants due to Speicher (see e.g.~\cite{Nica2006},~\cite{Speicher1994},~\cite{Speicher1997} and~\cite{Speicher1998}).

\subsubsection{Non-crossing partitions}

Let $S$ be a totally ordered set. A partition of the set $S$ is said to have a crossing if there exist $i,j,k,l \in S$, with $i < j < k < l$, such that $i$ and $k$ belong to some block of the partition and $j$ and $l$ belong to another block. If a partition has no crossings, it is called non-crossing. The set of all non-crossing partitions of $S$ is denoted by $NC(S)$. It is a lattice with respect to the fineness relation defined as follows: for all $\pi_1$ and $\pi_2\in NC(S)$, $\pi_1 \preceq \pi_2$ if every block of $\pi_1$ is contained in a block of $\pi_2$.

Let $n\in \mathbb{N}$. When $S = \left\{1, \ldots , n\right\}$, with its natural order, we will use the notation $NC (n)$. This set has been first considered by Kreweras in~\cite{Kreweras1972}. We denote by $0_n$ and $1_n$ respectively the minimal element  $\left\{\{1\}, \ldots , \{n\}\right\}$ of $NC (n)$, and the maximal element  $\left\{\{1, \ldots , n\}\right\}$ of $NC (n)$.

\subsubsection{Free cumulants}\label{Freecumulants}

For all $n\in \mathbb{N}$, $S\subset\left\{1, \ldots , n\right\}$, $\pi \in NC(S)$, and $A_1,\ldots,A_n\in \mathcal{A}$, set
$$\tau\left[\pi\right]\left(A_1,\ldots,A_n\right)=\displaystyle\prod_{V\in \pi} \tau\left(A_V\right)$$
where $A_V = A_{j_1}\cdots A_{j_k} $ if $V = \{j_1,... ,j_k\}$ is a block of the partition $\pi$, with $j_1 < j_2 < \ldots < j_k$. This way, $\tau\left[\pi\right]$ is a $n$-linear form on $\mathcal{A}$.

We use now the theory of Möbius inversion on lattices (see e.g. \cite{Stanley2011}). We denote by $\mu$ the Möbius function of the poset $\left(NC (n),\preceq\right)$, which is defined on $$\left\{\left(\sigma,\pi\right):\sigma\preceq \pi\right\}\subset NC (n)\times NC (n).$$
For all $n\in \mathbb{N}$, and $A_1,\ldots,A_n\in \mathcal{A}$, the free cumulant $\kappa\left(A_1,\ldots,A_n\right)$ is defined by
$$\kappa\left(A_1,\ldots,A_n\right)=\displaystyle\sum_{\pi\in NC(n)}\mu(\pi,1_n)\tau\left[\pi\right]\left(A_1,\ldots,A_n\right). $$
If $A_1=\cdots=A_n=A$, we call $\kappa\left(A_1,\ldots,A_n\right)$ the free cumulant of order $n$ of $A$, and we denote it by $\kappa_n\left(A\right)$.
For all $n\in \mathbb{N}$, $S\subset\left\{1, \ldots , n\right\}$, $\pi \in NC(S)$, and $A_1,\ldots,A_n\in \mathcal{A}$, set $$\kappa\left[\pi\right]\left(A_1,\ldots,A_n\right)=\displaystyle\prod_{V\in \pi} \kappa\left(A_V\right)$$
where $A_V = \left(A_{j_1},\ldots, A_{j_k}\right) $ if $V = \{j_1,... ,j_k\}$ is a block of the partition $\pi$, with $j_1 < j_2 < \ldots < j_k$. Similarly to $\tau$, $\kappa\left[\pi\right]$ is a $n$-linear form on $\mathcal{A}$. If $A_1=\cdots=A_n=A$, we denote $\kappa\left[\pi\right]\left(A_1,\ldots,A_n\right)$ by $\kappa\left[\pi\right](A)$.

For all $n\in \mathbb{N}$, $S\subset\left\{1, \ldots , n\right\}$, $\pi \in NC(S)$, and $A_1,\ldots,A_n\in \mathcal{A}$, we have (by the definition of the Möbius functions) the following relations
\begin{eqnarray}
\tau\left[\pi\right]\left(A_1,\ldots,A_n\right)&=&\displaystyle\sum_{\substack{\sigma\in NC(S)\\ \sigma \preceq\pi}}\kappa\left[\sigma\right]\left(A_1,\ldots,A_n\right)\label{taucum},\\
\kappa\left[\pi\right]\left(A_1,\ldots,A_n\right)&=&\displaystyle\sum_{\substack{\sigma\in NC(S)\\ \sigma \preceq\pi}}\mu(\sigma,\pi)\tau\left[\sigma\right]\left(A_1,\ldots,A_n\right)\label{cumtau}.
\end{eqnarray}
The importance of the free cumulants is in large part due to the following characterization of freeness.
\begin{proposition}
Let $\left(\mathcal{B}_i \right)_{i\in I}$ be subalgebras of $\mathcal{A}$. \label{freeness}They are free if and only if their mixed cumulants vanish. That is to say: for all $n\in \mathbb{N}^*$, all $i_1,\ldots, i_n \in I$ and all $A_1, \ldots , A_n \in \mathcal{A}$ such that $A_j$ belongs to some $\mathcal{B}_{i_j}$ for all $1\leq j \leq n$, whenever there exists some $j$ and $j' $ with $i_j\neq i_{j'}$, we have $\kappa(A_1,\ldots ,A_n) = 0$.
\end{proposition}
By linearity, this property has the following consequence for the computation of cumulants.
\begin{corollary}
For all $A_1,\ldots,A_n\in \mathcal{A}$ free from $B_1,\ldots,B_n \in \mathcal{A}$, we have\label{addfree} $$\kappa(A_1+B_1,\ldots ,A_n+B_n)=\kappa(A_1,\ldots ,A_n)+\kappa(B_1,\ldots ,B_n) .$$
\end{corollary}

\subsubsection{Semi-circular variables} Let $t\geq 0$. A non-commutative random variable $S_t$ is called semi-circular of variance $t$ if $S_t$ is self-adjoint and the free cumulants of $S_t$ are $\kappa_1(S_t)=0$, $\kappa_2(S_t)=t$ and $\kappa_n(S_t)=0$ for all $n>2$.\label{sc}

The distribution of $S_t$ is given by \eqref{taucum}: for all $n\in \mathbb{N}^*$, we have 
$\tau(S_t^n)=\sum_{\pi\in NC(n)} \kappa[\pi](S_t)$. The terms in the sum which are nonzero correspond to the non-crossing partitions composed of pairs of elements. For all $n\in \mathbb{N}^*$, we denote by $NC_2(n)$ the subset of $NC(n)$ of non-crossing partitions composed of pairs of elements. The cardinality of $NC_2(2n)$ is called the $n$-th Catalan number, and is denoted by $C_n$. Thus, for all $n\in \mathbb{N}^*$, we have 
$\tau(S_t^{2n})=\sum_{\pi\in NC_2(2n)} t^n=t^n C_n,$
and
$\tau(S_t^{2n-1})=0.$

\subsection{Additive transition operators}
Let $I$ be an arbitrary index set.
Let $\mathbf{A}=(A_i)_{i\in I}\in \mathcal{A}^I$. Let us define a derivation $\Delta_{\mathbf{A}}$ associated to $\mathbf{A}$ on $\left(\mathbb{C}\{X_i:i\in I\},\cdot_{\tr}\right)$ in the following way. 
For all $n\in \mathbb{N}$ and $ i(1),\ldots,i(n)\in I$, we set
$$\begin{array}{rccl}
\Delta_{\mathbf{A}} \left( X_{i(1)}\cdots X_{i(n)}\right)&=& \vspace{-0,5cm} \displaystyle\sum_{\substack{1\leq m\leq n\\ 1\leq k_1< \ldots < k_m\leq n}} & \kappa\left(A_{i(k_1)},\ldots,A_{i(k_m)}\right)X_{i(1)}\cdots X_{i(k_1-1)}\\
 & & &\hspace{1cm} \cdot \tr\left( X_{i(k_1+1)}\cdots X_{i(k_2-1)} \right)\\
  & & &\hspace{1cm}\cdot \tr\left( X_{i(k_2+1)}\cdots X_{i(k_3-1)} \right)\\
  & & &\hspace{1cm}\cdots\\
& & &\hspace{1cm}\cdot \tr\left( X_{i(k_{m-1}+1)}\cdots X_{i(k_m-1)} \right)\\
& & &\hspace{1cm} \cdot  X_{i(k_m+1)}\cdots X_{i(n)}
\end{array}$$
and we extend $\Delta_{\mathbf{A}}$ to all $\mathbb{C}\{X_i:i\in I\}$ by linearity and by the relation of derivation
$$\forall P,Q \in \mathbb{C}\{X_i:i\in I\},\ \Delta_{\mathbf{A}} \left(P\tr Q\right)=\left(\Delta_{\mathbf{A}} P\right)\tr Q+P\tr \left(\Delta_{\mathbf{A}} Q\right).$$

\begin{proposition}
\label{com}Let $\mathbf{A}=\left(A_i\right)_{i\in I}$ and $\mathbf{B}=\left(B_i\right)_{i\in I}\in \mathcal{A}^I$ be two families of elements of $\mathcal{A}$.
Then $\Delta_{\mathbf{A}}$ and $ \Delta_{\mathbf{B}}$ commute. Moreover, if $\mathbf{A}$ and $\mathbf{B}$ are free, we have $\Delta_{\mathbf{A}+\mathbf{B}}=\Delta_{\mathbf{A}}+\Delta_{\mathbf{B}}$.
\end{proposition}
\begin{proof}We recall that, from Section~\ref{Anotprod}, the monomials of $ \mathbb{C}\left\langle X_i:i\in I\right\rangle$ generate the algebra $\left(\mathbb{C}\{X_i:i\in I\},\cdot_{\tr} \right)$. Let $M$ be a monomial of $\mathbb{C}\left\langle X_i:i\in I\right\rangle$. Let  $n\in \mathbb{N}$ and $ i(1),\ldots,i(n)\in I$ such that $M=X_{i(1)}\cdots X_{i(n)}$. For all $1\leq k<l\leq n$, we will denote by $\tr_{k,l}$ the element $\tr( X_{i(k+1)}\cdots X_{i(l-1)})$. Since $\Delta_{\mathbf{B}}$ is a derivation for $\cdot_{\tr}$, we have
$$\begin{array}{rccl}
\Delta_{\mathbf{B}}\left( \Delta_{\mathbf{A}}\left( X_{i(1)}\cdots X_{i(n)}\right)\right)&=& \vspace{-0,5cm} \displaystyle\sum_{\substack{1\leq m\leq n\\ 1\leq k_1< \ldots < k_m\leq n}} & \kappa\left(A_{i(k_1)},\ldots,A_{i(k_m)}\right)\cdot \tr_{k_1,k_2}\cdots \tr_{k_{m-1},k_m}\\
& & &\hspace{1cm} \cdot  \Delta_{\mathbf{B}}\left(X_{i(1)}\cdots X_{i(k_1-1)} X_{i(k_m+1)}\cdots X_{i(n)}\right)\\
& & + \vspace{-0,7cm}\displaystyle\sum_{\substack{1\leq m\leq n\\ 1\leq k_1< \ldots < k_m\leq n\\ 1\leq p\leq m-1}}&\kappa\left(A_{i(k_1)},\ldots,A_{i(k_m)}\right)X_{i(1)}\cdots X_{i(k_1-1)}\\
 & & &\hspace{1cm} \cdot \tr_{k_1,k_2}\cdots \Delta_{\mathbf{B}} ( \tr_{k_p,k_{p+1}}) \cdots  \tr_{k_{m-1},k_m}\\
& & &\hspace{1cm} \cdot  X_{i(k_m+1)}\cdots X_{i(n)}.
\end{array}$$
Let us fix $1\leq m\leq n$ and $1\leq k_1< \ldots < k_m\leq n$ in the first sum above and compute the corresponding term of the sum:

$$\begin{array}{llcll}
 \kappa\left(A_{i(k_1)},\ldots,A_{i(k_m)}\right)\cdot\Big(&& \hspace{-1cm} \vspace{-0,5cm} \displaystyle\sum_{\substack{1\leq m'\leq n\\ 1\leq l_1< \ldots < l_{m'} <k_1}}  \hspace{-0,5cm}& \kappa\left(B_{i(l_1)},\ldots,B_{i(l_{m'})}\right)X_{i(1)}\cdots X_{i(l_1-1)}&\\
  &&&\hspace{1cm} \cdot \tr_{l_1,l_2}\cdots \tr_{l_{m'-1},l_{m'}}\cdot  X_{i(l_{m'}+1)}\cdots X_{i(k_1-1)}&\\
  &&&\hspace{1cm} \cdot \tr_{k_1,k_2}\cdots \tr_{k_{m-1},k_m} \cdot X_{i(k_m+1)}\cdots X_{i(n)}&\\
&+& \hspace{-1cm}\vspace{-0,5cm} \displaystyle\sum_{\substack{1\leq m'\leq n\\ k_m< l_1< \ldots < l_{m'}\leq n}}  \hspace{-0,5cm}& \kappa\left(B_{i(l_1)},\ldots,B_{i(l_{m'})}\right)X_{i(1)}\cdots X_{i(k_1-1)}\\
& &&\hspace{1cm} \cdot \tr_{k_1,k_2}\cdots \tr_{k_{m-1},k_m} \cdot  X_{i(k_{m}+1)}\cdots X_{i(l_1-1)}&\\
 &&&\hspace{1cm} \cdot \tr_{l_1,l_2}\cdots \tr_{l_{m'-1},l_{m'}} \cdot X_{i(l_{m'}+1)}\cdots X_{i(n)}&\\
&+& \hspace{-1cm} \vspace{-1cm}\displaystyle\sum_{\substack{1\leq m'\leq n\\1\leq q\leq m'-1\\ 1\leq l_1< \ldots l_q< k_1\\k_m<l_{q+1}<\ldots<l_{m'}\leq n}}\hspace{-0,5cm}&\kappa\left(B_{i(l_1)},\ldots,B_{i(l_{m'})}\right)X_{i(1)}\cdots X_{i(l_1-1)}&\\
  &&&\hspace{1cm} \cdot \tr_{l_1,l_2}\cdots \tr_{l_{q-1},l_q}\cdot \tr_{k_1,k_2}\cdots \tr_{k_{m-1},k_m}&\\
& &&\hspace{1cm}  \cdot \tr\left(X_{i(l_q+1)}\cdots X_{i(k_1-1)} X_{i(k_m+1)}\cdots X_{i(l_{q+1}-1)}\right)&\\
 &&&\hspace{1cm} \cdot \tr_{l_{q+1},l_{q+2}}\cdots \tr_{l_{m'-1},l_{m'}} \cdot  X_{i(l_{m'}+1)}\cdots X_{i(n)}&\Big).
\end{array}$$
When summing over the indices $1\leq k_1< \ldots < k_m\leq n$, the elements $\mathbf{A}$ and $\mathbf{B}$ play symmetric roles in the two first sums, while the last sum gives us the term
$$\begin{array}{rccl} & &  \vspace{-0,7cm}\displaystyle\sum_{\substack{1\leq m'\leq n\\ 1\leq l_1< \ldots < l_{m'}\leq n\\ 1\leq q\leq m'-1}}&\kappa\left(B_{i(k_1)},\ldots,B_{i(k_m)}\right)X_{i(1)}\cdots X_{i(l_1-1)}\\
 & & &\hspace{1cm} \cdot \tr_{l_1,l_2}\cdots \Delta_{\mathbf{A}} ( \tr_{l_q,l_{q+1}}) \cdots  \tr_{l_{m'-1},l_{m'}}\\
& & &\hspace{1cm} \cdot  X_{i(l_{m'}+1)}\cdots X_{i(n)}
\end{array}$$
which is exactly the term needed to achieve the symmetry of the roles of $\mathbf{A}$ and $\mathbf{B}$ in $\Delta_{\mathbf{B}} \Delta_{\mathbf{A}}( M)$. Thus, $\Delta_{\mathbf{A}}\Delta_{\mathbf{B}}(M)=\Delta_{\mathbf{B}}\Delta_{\mathbf{A}}(M)$ for all monomials of $\mathbb{C}\langle X_i: i\in I\rangle$, and the commutativity is extended by recurence on all $\mathbb{C}\{X_i: i\in I\}$ because they are both derivations. Indeed, for all $P,Q \in \mathbb{C}\{X_i: i\in I\}$,
\begin{eqnarray*}\Delta_{\mathbf{A}} \Delta_{\mathbf{B}}\left(P\tr Q\right)=\left(\Delta_{\mathbf{A}}\Delta_{\mathbf{B}} P\right)\tr \left(Q\right)+ P\tr \left(\Delta_{\mathbf{A}}\Delta_{\mathbf{B}}Q\right)+\left(\Delta_{\mathbf{A}}P\right)\tr \left(\Delta_{\mathbf{B}} Q\right)+\left(\Delta_{\mathbf{B}} P\right)\tr \left(\Delta_{\mathbf{A}}Q\right).\end{eqnarray*}

If $\mathbf{A}$ and $\mathbf{B}$ are free, we have from Corollary~\ref{addfree} that $\Delta_{\mathbf{A}+\mathbf{B}}M=(\Delta_{\mathbf{A}}+\Delta_{\mathbf{B}})M$ for all monomials $M$ of $\mathbb{C}\langle X_i: i\in I\rangle$, and since $\Delta_{\mathbf{A}+\mathbf{B}}$ and $\Delta_{\mathbf{A}}+\Delta_{\mathbf{B}}$ are both derivations, they are equal on $\mathbb{C}\{X_i: i\in I\}$.
\end{proof}
Let $\mathbf{A}=(A_i)_{i\in I}\in \mathcal{A}^I$. The operator $\Delta_{\mathbf{A}}$ makes the degree strictly decrease. Thus, the operator $e^{\Delta_{\mathbf{A}}}$ on $\mathbb{C}\{X_i:i\in I\}$ is defined by the formal series $\sum_{k=0}^{\infty}\frac{1}{k!} \Delta_{\mathbf{A}}^k.$
For all $P\in \mathbb{C}_d\{X_i:i\in I\}$, $e^{\Delta_{\mathbf{A}}} P$ is the finite sum $\sum_{k=0}^{d} \frac{1}{k!}  \Delta_{\mathbf{A}}^k P$.
The operator $\Delta_{\mathbf{A}}$ is a derivation, and we therefore have the Leibniz formula
\begin{eqnarray*}\forall k\in \mathbb{N}, \forall P,Q \in \mathbb{C}\{X_i: i\in I\},\left(  \Delta_{\mathbf{A}}\right)^k \left(P\tr Q\right) =\displaystyle\sum_{l=0}^k \binom{k}{l} \left(  \Delta_{\mathbf{A}}^l P \right)\tr \left(  \Delta_{\mathbf{A}}^{k-l} Q\right),\end{eqnarray*}
from which we deduce, using the standard power series argument, that the operator $e^{ \Delta_{\mathbf{A}}}$ is multiplicative in the following sense:
\begin{eqnarray*}\forall P,Q \in \mathbb{C}\{X_i: i\in I\},\ e^{ \Delta_{\mathbf{A}}} \left(P\tr Q\right)=\left(e^{ \Delta_{\mathbf{A}}} P\right)\tr \left(e^{\Delta_{\mathbf{A}}}Q\right).\end{eqnarray*}
\begin{theorem}
Let $I$ be an arbitrary index set. Let $\mathbf{A}=(A_i)_{i\in I}\in \mathcal{A}^I$. For all $P\in \mathbb{C}\{ X_i:i\in I\}$, and all $\mathbf{B}=(B_i)_{i\in I}\in \mathcal{A}^I$ free from $(A_i)_{i\in I}$, we have
$$\tau\left(\left.P\left(\mathbf{A}+\mathbf{B}\right)\right|\mathbf{B}\right)=(e^{\Delta_{\mathbf{A}}}P)\left(\mathbf{B}\right).$$
In particular, for all $P\in \mathbb{C}\{ X_i:i\in I\}$, we have
$$ \tau\left(P\left(\mathbf{A}\right)\right)=(e^{\Delta_{\mathbf{A}}}P)\left(0\right).$$
\label{freekernel}
\end{theorem}

\subsubsection*{Example} Let $t\geq 0$. Let $S_t$ be a semi-circular random variable of variance $t$. For all $B\in \mathcal{A}$ free from $S_t$, let us compute $\tau((S_t+B)^3|B)$.
We have 
$\Delta_{S_t} X^3=2t X +t \tr(X)$ and $(\Delta_{S_t})^2 X^3=\Delta_{S_t}(2t X +t \tr(X))=0$.
Thus, $
e^{\Delta_{S_t}}(X^3)=X^3 +  \Delta_{S_t}X^3+0
=X^3+2t X +t \tr(X).
$ Using Theorem~\ref{freekernelmult}, we have, for all $B\in \mathcal{A}$ free from $S_t$, $$\tau\left(\left(S_t+B\right)^3|B\right)= B^3+2t B +t \tau(B).$$

\begin{proof}[Proof of Theorem~\ref{freekernel}]
One could prove Theorem~\ref{freekernel} directly but it would be very combinatorial and we prefer to present here a more dynamical proof where the combinatorics only appear infinitesimally.
We shall prove first a dynamical lemma and a weaker version of Theorem~\ref{freekernel} before proving the theorem in all generality.

Let us first prove Lemma~\ref{eqdiff}, which will be useful in the rest of the paper.
\begin{lemma}Let $J$ be a finite index set and $d\in \mathbb{N}$. Let $L$ be a linear operator on $\mathbb{C}_d\{X_i:i\in J\}$ and $\left(\phi_t\right)_{t\geq 0}$ be linear functionals on $\mathbb{C}_d\{X_i:i\in J\}$.
Let us assume that, for all $P\in \mathbb{C}_d\{X_i:i\in J\}$, $t\mapsto \phi_t(P)$ is differentiable on $[0,\infty)$, and for all $t\geq 0$,
$$\frac{\diff}{\diff t} \phi_t(P)=\phi_t(LP).$$
If there exists $\mathbf{A}=(A_i)_{i\in J}\in \mathcal{A}^I$ such that $\phi_0(P)=P\left(\mathbf{A}\right)$, then for all $t\geq 0$ and all $P\in \mathbb{C}_d\{X_i:i\in J\}$, we have\label{eqdiff}
$$\phi_t(P)=e^{tL}P\left(\mathbf{A}\right). $$
\end{lemma}
\begin{proof}Let us fix a finite basis $\{P_b:b\in B\}$ of $\mathbb{C}_d\{X_i:i\in J\}$. For all $P\in \mathbb{C}_d\{X_i:i\in J\}$, let us denote by $(\alpha_b(P))_{b\in B}$ the coefficients of $P$ in the basis $\{P_b:b\in B\}$.
We claim that the functions $(t\mapsto\phi_t(P_b))_{b\in B}$ and $(t\mapsto e^{tL}P_b(\mathbf{A}))_{b\in B}$ are both the unique solution to the multidimensional differential equation
\begin{equation*}
\left\{
\begin{array}{r c l}
    y_b(0) &=& P_b\left(\mathbf{A}\right), \\
    y_b' &=& \displaystyle\sum_{c\in B}\alpha_c\left(L P_b\right)\cdot y_c, \forall b\in B. 
\end{array}
\right.\end{equation*}
Indeed, let $b\in B$. We have
$$\frac{\diff}{\diff t}e^{tL}P_b\left(\mathbf{A}\right)= e^{tL}LP_b\left(\mathbf{A}\right)=\displaystyle\sum_{c\in B}\alpha_c\left(L P_b\right)\cdot e^{tL}P_c\left(\mathbf{A}\right)$$
and
$$\frac{\diff}{\diff t} \phi_t(P_b)=\phi_t(L P_b)=\displaystyle\sum_{c\in B}\alpha_c\left(LP_b\right)\cdot \phi_t(P_c).$$
Furthermore, $e^{0 L}P_b\left(\mathbf{A}\right)=\phi_0(P_b)=P_b\left(\mathbf{A}\right)$. Hence, we have $\phi_t(P_b)=e^{tL}P_b\left(\mathbf{A}\right)$ for all $b\in B$. We extend the relation $\phi_t(P)=e^{tL}P(\mathbf{A})$ from $P\in \{P_b:b\in B\}$ to all $\mathbb{C}_d\{X_i:i\in J\}$ by linearity.
\end{proof}
The following lemma is a weak version of Theorem~\ref{freekernel}.\label{freekernelproof}
\begin{lemma}\label{freekernelweak}
Let $I$ be an arbitrary index set. Let $\mathbf{A}=(A_i)_{i\in I}\in \mathcal{A}^I$. For all $P\in \mathbb{C}\{ X_i:i\in I\}$, we have
$$ \tau\left(P\left(\mathbf{A}\right)\right)=(e^{\Delta_{\mathbf{A}}}P)\left(0\right).$$
\end{lemma}
\begin{proof}For all $n\in \mathbb{N}$ and $\pi\in NC(n)$, let us denote by $| \pi|$ the number of blocks in $\pi$.

Let $t\geq 0$. 
For all $n\in \mathbb{N}$ and $i(1),\ldots, i(n)\in I$, we set (with the convention $0^0=1$)
\begin{eqnarray*}
\phi_t\left(X_{i(1)}\cdots X_{i(n)}\right)&=&\displaystyle\sum_{ \pi\in NC(n)} t^{| \pi|}\kappa[\pi](A_{i(1)} ,\ldots, A_{i(n)}  )
\end{eqnarray*}
and we extend $\phi_t$ to all $\left(\mathbb{C}\{ X_i:i\in I\},\cdot_{\tr} \right)$ by linearity and by the multiplicative relation
$$\forall P,Q \in \mathbb{C}\{ X_i:i\in I\}, \phi_t( P\tr Q)= \phi_t( P)\tau ( \phi_t( Q)).$$
We remark that, using \eqref{taucum}, we find $\tau\left(M\left(\mathbf{A}\right)\right)=\phi_1(M)$ for all monomials $M$ of $\mathbb{C}\langle X_i:i\in I\rangle $. Moreover, the map $P\mapsto \tau(P(\mathbf{A}))$ satisfies the same multiplicative relation as $\phi_1$. Indeed, we have
$$\forall P,Q \in \mathbb{C}\{ X_i:i\in I\},\ \tau\Big(P\tr Q\left(\mathbf{A}\right)\Big)= \tau\Big(P\left(\mathbf{A}\right)\Big)\tau \Big( \tau\Big(Q\left(\mathbf{A}\right)\Big)\Big).$$
It follows that $\phi_1(P)= \tau\left(P\left(\mathbf{A}\right)\right)$ for all $P \in \mathbb{C}\{ X_i:i\in I\}$. Thus, it remains to prove that, for all $\mathbb{C}\{ X_i:i\in I\}$, we have $\phi_1(P)=(e^{\Delta_{\mathbf{A}}}P)\left(0\right)$. In order to use Lemma~\ref{eqdiff} in the third step, we will prove that, for all $P\in \mathbb{C}\{ X_i:i\in I\}$, $\phi_0(P)=P(0)$ in the first step, and that, for all $t\geq 0$ and all $P\in  \mathbb{C}\{ X_i:i\in I\}$, we have $\frac{\diff}{\diff t} \phi_t(P)=\phi_t(\Delta_{\mathbf{A}}P)$ in the second step.

\subsubsection*{Step 1}
For all $P\in \mathbb{C}\{ X_i:i\in I\}$, $\phi_0(P)=P(0)$. Indeed, $\phi_0(1)=1=M(0)$ for $M=1\in \mathbb{C}\langle X_i:i\in I\rangle $ and $\phi_0(M)=0=M(0)$ for all monomials $M$ of $\mathbb{C}\langle X_i:i\in I\rangle $ with non-zero degree. 
We infer the equality $\phi_0(P)=P(0)$ to all $ \mathbb{C}\{ X_i:i\in I\}$, since the evaluation satisfies the same multiplicative relation as $\phi_0$. Indeed, we have
$$\forall P,Q \in \mathbb{C}\{ X_i:i\in I\},\ (P\tr Q)\left(0\right)= P\left(0\right)\tau \left(Q\left(0\right)\right).$$

\subsubsection*{Step 2}
We prove now that, for all $t\geq 0$ and all $P\in \mathbb{C}\{ X_i:i\in I\}$, we have $\frac{\diff}{\diff t} \phi_t(P)=\phi_t(\Delta_{\mathbf{A}}P)$.

Let $t\geq 0$ and $M$ be a monomial of $  \mathbb{C}\langle X_i:i\in I\rangle$. Let us fix $n\in \mathbb{N}$ and $i(1),\ldots, i(n)\in I$ such that $M=X_{i(1)}\cdots X_{i(n)}$.
We have
\begin{eqnarray*}
\frac{\diff}{\diff t}\phi_t\left(X_{i(1)}\cdots X_{i(n)}\right)
&=&\displaystyle\sum_{ \pi\in NC(n)} | \pi| t^{| \pi|-1}\kappa[\pi]\left(A_{i(1)} ,\ldots, A_{i(n)}  \right)\\
&=&\displaystyle\sum_{ \pi\in NC(n)} \hspace{0.3cm}\displaystyle\sum_{V\in \pi}  t^{| \pi|-1}\kappa[\pi]\left(A_{i(1)} ,\ldots, A_{i(n)}  \right)\\
&=&\displaystyle\sum_{V \subset \{1,\ldots,n\}} \hspace{0.3cm}\displaystyle\sum_{\substack{ \pi\in NC(n)\\ V\in \pi}}  t^{| \pi|-1}\kappa[\pi]\left(A_{i(1)} ,\ldots, A_{i(n)}  \right).
\end{eqnarray*}
Let us fix a subset $V=\{k_1,\ldots,k_m\}$ of $\{1,\ldots,n\}$ such that $1\leq k_1< \ldots < k_m\leq n$.
Let us denote $W_1=\{k_1+1,\ldots,k_2-1\},W_2=\{k_2+1,\ldots,k_3-1\},\ldots,W_{m-1}=\{k_{m-1}+1,\ldots,k_m-1\}$ and $W_m=\{1,\ldots,k_1-1,k_m+1,\ldots,k_n\}$.
To each partition $\pi\in NC(n)$ such that $V\in \pi$, we associate for all $1\leq i\leq m$ the partition $\pi_i \in NC(W_i)$ induced by $\pi$ on $W_i$. Conversely, to each $m$-tuple $\pi_1 \in NC(W_1), \ldots, \pi_m \in NC(W_m)$, we associate the non-crossing partition $\pi=\{V\}\cup \pi_1\cup\cdots\cup \pi_m \in NC(n)$. Because of the non-crossing condition, this leads to a bijective correspondance
$$\{ \pi\in NC(n):\ V\in \pi\}  \leftrightarrow \{ (\pi_1,\ldots, \pi_m)\in NC(W_1)\times \cdots \times NC(W_m)\}$$
which allows us to sum separately each non-crossing partition on each subset $W_1,\ldots,W_{m}$.

We have to examine now how the terms of the sum are transformed. Let us fix $\pi\in NC(n)$ such that $V\in \pi$. Let us define $(\pi_1,\ldots, \pi_m)\in NC(W_1)\times \cdots \times NC(W_m)$ as before. We have
$$t^{| \pi|-1}=t^{\left| \{V\}\cup \pi_1\cup\cdots\cup \pi_m\right|-1}=t^{|\pi_1|}\cdots t^{|\pi_m|}$$
and
$$
\kappa[\pi](A_{i(1)} ,\ldots, A_{i(n)}  )=
\kappa(A_{i(k_1)} ,\ldots, A_{i(k_m)} )\cdot \kappa[ \pi_1](A_{i(1)} ,\ldots, A_{i(n)} ) \cdots \kappa[ \pi_m](A_{i(1)} ,\ldots, A_{i(n)} ).
$$
We infer
\begin{eqnarray*}
& &\hspace{-1cm}\displaystyle\sum_{\substack{ \pi\in NC(n)\\ V\in \pi}}  t^{| \pi|-1}\kappa[\pi](A_{i(1)} ,\ldots, A_{i(n)}  )\\
&=&
\kappa(A_{i(k_1)} ,\ldots, A_{i(k_m)} )\\
& &
\hspace{1cm}\cdot\left( \displaystyle\sum_{ \pi_1\in NC(W_1)} t^{| \pi|}\kappa[\pi](A_{i(1)} ,\ldots, A_{i(n)}  )\right)\cdots\left( \displaystyle\sum_{ \pi_m\in NC(W_m)} t^{| \pi|}\kappa[\pi](A_{i(1)} ,\ldots, A_{i(n)})  \right)\\
&=&\kappa(A_{i(k_1)} ,\ldots, A_{i(k_m)} )\\
& &
\hspace{1cm}\cdot\phi_t\left( X_{i(k_1+1)}\cdots X_{i(k_2-1)}\right)\\
  &&\hspace{1cm}\cdots\\
 &&\hspace{1cm}\cdot \phi_t\left( X_{i(k_{m-1}+1)}\cdots X_{i(k_m-1)} \right)\\
  &&\hspace{1cm} \cdot  \phi_t\left( X_{i(1)}\cdots X_{i(k_1-1)} X_{i(k_m+1)}\cdots X_{i(n)}\right).
  \end{eqnarray*}
Finally, we use the definition of $\Delta_{\mathbf{A}}$ and the multiplicative relation of $\phi_t$ to infer
$$\begin{array}{rcl}
&&\hspace{-2cm}\displaystyle\frac{\diff}{\diff t}\phi_t\left(X_{i(1)}\cdots X_{i(n)}\right)\\
&=& \displaystyle\sum_{\substack{1\leq m\leq n\\1\leq k_1< \ldots < k_m\leq n}}\hspace{0.3cm} \displaystyle\sum_{\substack{ \pi\in NC(n)\\ \{k_1,\ldots, k_m\}\in \pi}}  t^{| \pi|-1}\kappa[\pi](A_{i(1)} ,\ldots, A_{i(n)}  )\\
&=&\vspace{-0,5cm} \displaystyle\sum_{\substack{1\leq m\leq n\\1\leq k_1< \ldots < k_m\leq n}}  \kappa\left(A_{i(k_1)},\ldots,A_{i(k_m)}\right)\\
 && \hspace{3cm} \cdot \tau\left(\phi_t\left( X_{i(k_1+1)}\cdots X_{i(k_2-1)} \right)\right)\\
  &&\hspace{3cm}\cdots\\
 &&\hspace{3cm}\cdot \tau\left(\phi_t\left( X_{i(k_{m-1}+1)}\cdots X_{i(k_m-1)} \right)\right)\\
  &&\hspace{3cm} \cdot  \phi_t\left( X_{i(1)}\cdots X_{i(k_1-1)} X_{i(k_m+1)}\cdots X_{j(n)}\right)\\
  &=&\phi_t\left(\Delta_{\mathbf{A}} \left( X_{i(1)}\cdots X_{i(n)}\right)\right).
\end{array}$$
We extend the equality $\frac{\diff}{\diff t}\phi_t(P)=\phi_t\left(\Delta_{\mathbf{A}} P\right)$ from monomials of $\mathbb{C}\langle X_i:i\in I\rangle$ to all $P\in \mathbb{C}\{X_i:i\in I\}$ by linearity and by the following induction. If $P$ and $Q\in\mathbb{C}\{X_i:i\in I\}$ verify $\frac{\diff}{\diff t}\phi_t(P)=\phi_t\left(\Delta_{\mathbf{A}} P\right)$ and $\frac{\diff}{\diff t}\phi_t(Q)=\phi_t\left(\Delta_{\mathbf{A}} Q\right)$, we have
$$\begin{array}{rcl}
\displaystyle\frac{\diff}{\diff t} \phi_t\left(P\tr Q\right)&=&\displaystyle\frac{\diff}{\diff t}\Big(\phi_t\left(P\right)\tau\left(\phi_t(Q)\right)\Big) \\
&=&\left(\displaystyle\frac{\diff}{\diff t} \phi_t\left(P\right)\right)\tau\left(\phi_t(Q)\right)+\phi_t\left(P\right)\tau\left(\displaystyle\frac{\diff}{\diff t}\phi_t(Q)\right)\\
&=& \phi_t\left(\Delta_{\mathbf{A}} P\right)\tau\left(\phi_t(Q)\right)+\phi_t\left(P\right)\tau\left(\phi_t\left(\Delta_{\mathbf{A}} Q\right)\right)\\
&=& \phi_t\Big( \Delta_{\mathbf{A}} P\cdot\tau\left(Q\right)+P\cdot\tau\left(\Delta_{\mathbf{A}} Q\right)\Big)\\
&=&\phi_t\Big(\Delta_{\mathbf{A}} \left(P\tr Q\right)\Big).
\end{array}$$

\subsubsection*{Step $3$}Let $P\in \mathbb{C}\{X_i:i\in I\}$. There exists a finite index set $J\subset I$ and $d\in \mathbb{N}$ such that $P\in\mathbb{C}_d\{X_i:i\in J\}$. We remark that $\Delta_{\mathbf{A}}$ is a linear operator on $\mathbb{C}_d\{X_i:i\in J\}$ and that  $(\phi_t)_{t\geq 0}$ are linear functionals. We deduce from Lemma~\ref{eqdiff} and Steps $1$ and $2$ that, for all $t\geq 0$, we have $\phi_t(P)=(e^{t\Delta_{\mathbf{A}}}P)\left(0\right)$. In particular, $(e^{\Delta_{\mathbf{A}}}P)\left(0\right)=\phi_1(P)=\tau\left(P\left(\mathbf{A}\right)\right)$.
\end{proof}

Let us finish the proof of Theorem~\ref{freekernel}. Let $\mathbf{A}=(A_i)_{i\in I}\in \mathcal{A}^I$ and $\mathbf{B}=(B_i)_{i\in I}\in \mathcal{A}^I$ be free, and $P\in \mathbb{C}\{ X_i:i\in I\}$. We remark that, thanks to Proposition~\ref{com}, we know that $\Delta_{\mathbf{A}+\mathbf{B}}=\Delta_{\mathbf{A}}+\Delta_{\mathbf{B}}$, and that $\Delta_{\mathbf{A}}$ and $\Delta_{\mathbf{B}}$ commute. We deduce directly from Lemma~\ref{freekernelweak} that
 \begin{equation}
 \tau\left(P(\mathbf{A}+\mathbf{B})\right)=(e^{\Delta_{\mathbf{A}+\mathbf{B}}}P)\left(0\right)
 =(e^{\Delta_{\mathbf{A}}+\Delta_{\mathbf{B}}}P)\left(0\right)
 =(e^{\Delta_{\mathbf{B}}}e^{\Delta_{\mathbf{A}}}P)\left(0\right)
 =\tau\left((e^{\Delta_{\mathbf{A}}}P)\left(\mathbf{B}\right)\right).\label{freekernelaux}
\end{equation}
We will now use the following characterization of conditional expectation. The element $\tau(P(\mathbf{A}+\mathbf{B})|\mathbf{B})$ is the unique element of $W^*(\mathbf{B})$ such that, for all $B_{{i_0}} \in W^*(\mathcal{B})$, $$\tau\left(P(\mathbf{A}+\mathbf{B})B_{{i_0}}\right) = \tau\left(\tau\left(\left.P(\mathbf{A}+\mathbf{B})\right|\mathbf{B}\right)B_{{i_0}}\right).$$
Since $(e^{\Delta_{\mathbf{A}}}P)\left(\mathbf{B}\right)\in W^*(\mathbf{B})$, it remains to prove that, for all $B_{{i_0}} \in W^*(\mathcal{B})$,
$$\tau\left(P(\mathbf{A}+\mathbf{B})B_{{i_0}}\right)= \tau\left(e^{\Delta_{\mathbf{A}}}P\left(\mathbf{B}\right)B_{{i_0}}\right).$$
In order to use \eqref{freekernelaux}, we prefer to work on $\mathbb{C}\{ X_i:i\in I\cup\{{i_0}\}\}$. Let $R_{{i_0}}:P\mapsto PX_{i_0}$ be the operator of right multiplication by $X_{i_0}$ on $\mathbb{C}\{ X_i:i\in I\cup\{{i_0}\}\}$.
Let $A_{i_0}=0$ and $B_{{i_0}} \in W^*(\mathcal{B})$.
On one hand, we have $P(\mathbf{A}+\mathbf{B})B_{{i_0}}=(R_{{i_0}}P)(\mathbf{A}+\mathbf{B},A_{i_0}+B_{{i_0}})$, and using \eqref{freekernelaux}, we have $\tau\left(P(\mathbf{A}+\mathbf{B})B_{{i_0}}\right)=\tau\left(\left(e^{\Delta_{\mathbf{A},A_{i_0}}}R_{{i_0}}(P)\right)\left(\mathbf{B},B_{{i_0}}\right)\right).$ On the other hand, $\tau\left(e^{\Delta_{\mathbf{A}}}P\left(\mathbf{B}\right)B_{{i_0}}\right)=\tau\left(\left(R_{{i_0}}e^{\Delta_{\mathbf{A},A_{i_0}}}(P)\right)\left(\mathbf{B},B_{{i_0}}\right)\right)$. Thus, it remains to prove that the operators $\Delta_{\mathbf{A},A_{i_0}}$ and $R_{{i_0}}$ commute. Let us check this on monomials. For all $n\in \mathbb{N}$ and $ i(1),\ldots,i(n)\in I\cup\{{i_0}\}$, let us fix $i(n+1)={i_0}$. For all $1\leq k<l\leq n+1$, we will denote by $\tr_{k,l}$ the element $\tr( X_{i(k+1)}\cdots X_{i(l-1)})$. Because all free cumulants involving $A_{{i_0}}=0$ are equal to zero, we have
$$\begin{array}{rccl}
&&\hspace{-3cm}\Delta_{\mathbf{A}}R_{{i_0}} \left( X_{i(1)}\cdots X_{i(n)}\right)& \\
&=& \vspace{-0,5cm} \displaystyle\sum_{\substack{1\leq m\leq n\\1\leq k_1< \ldots < k_m\leq n+1}} & \kappa\left(A_{i(k_1)},\ldots,A_{i(k_m)}\right)X_{i(1)}\cdots X_{i(k_1-1)}\\
 & & &\vspace{0,5cm} \hspace{1cm}\cdot \tr_{k_1,k_2}\cdots \tr_{k_{m-1},k_m} \cdot  X_{i(k_m+1)}\cdots X_{i(n+1)}\\
&=&\vspace{-0,5cm} \displaystyle\sum_{\substack{1\leq m\leq n\\1\leq k_1< \ldots < k_m\leq n}} & \kappa\left(A_{i(k_1)},\ldots,A_{i(k_m)}\right)X_{i(1)}\cdots X_{i(k_1-1)}\\
 & & & \vspace{0,2cm} \hspace{1cm}\cdot \tr_{k_1,k_2}\cdots \tr_{k_{m-1},k_m} \cdot  X_{i(k_m+1)}\cdots X_{i(n)}X_{i_0}\\
 &=&R_{{i_0}} \Delta_{\mathbf{A}}\left( X_{i(1)}\cdots X_{i(n)}\right) .&
\end{array}$$
The commutativity is extended to all $\mathbb{C}\{ X_i:i\in I\cup\{{i_0}\}\}$ by an immediate induction, because for all $P,Q\in \mathbb{C}\{ X_i:i\in I\cup\{{i_0}\}\}$, we have
$$\Delta_{\mathbf{A}}R_{{i_0}}(P\tr Q)=\Big(\Delta_{\mathbf{A}}R_{{i_0}}(P)\Big)\tr Q+R_{{i_0}}\Big(P\tr(\Delta_{\mathbf{A}}(Q))\Big) $$
and
\begin{equation*}
R_{{i_0}}\Delta_{\mathbf{A}}(P\tr Q)=\Big(R_{{i_0}}\Delta_{\mathbf{A}}(P)\Big)\tr Q+R_{{i_0}}\Big(P\tr(\Delta_{\mathbf{A}}(Q))\Big)  .\qedhere
\end{equation*}
\end{proof}

\subsection{Free log-cumulants}Consider two free random variables $A,B\in \mathcal{A}$. The distribution of $AB$ is determined by the distributions of $A$ and $B$ separately. The machinery for doing this computation, in the same spirit as Corollary~\ref{addfree}, is the $S$-transform (see Section 3.6 in~\cite{Voiculescu1992}).

Since we have to deal with more than two variables, it is necessary for our computations to introduce the free log-cumulants. There are some quantities introduced by Mastnak and Nica in~\cite{Mastnak2010}, which have nice behavior with the product of free elements of $\mathcal{A}$.
 The name "log-cumulants", not present in~\cite{Mastnak2010}, comes from analogous objects with that name in the context of classical multiplicative convolution~\cite{NICOLAS}. In this section, we give a definition of such quantities, in a slightly different presentation from~\cite{Mastnak2010}. The link is made in~Section~\ref{mastnak}, and we present some of the properties of the free log-cumulants in Section~\ref{frprop}.

\subsubsection{The Kreweras complementation map}There exist several equivalent definitions of the Kreweras complementation map. See the book \cite{Nica2006} for a detailed presentation.

Let $n\in \mathbb{N}^*$. We endow the set $\{1,1',\ldots,n,n'\}$ with the cylic order $(1,1',\ldots,n,n')$.\label{noncrosspart} For all $\pi_1,\pi_2 \in NC(n)$, we form a not necessarily non-crossing partition $\pi_1\cup \pi_2$ of $\{1,1',\ldots,n,n'\}$ by identifying $\pi_2$ as a partition of $\{1',\ldots,n'\}$, and merging $\pi_1$ and $\pi_2$. Let $\pi=\{V_1,\ldots, V_l\}\in NC(n)$. For all $1\leq i\leq l$, we denote by $V_i'\subset \{1',\ldots,n'\}$ the image of $V_i$ via the isomorphism $(1,\ldots,n)\simeq(1',\ldots,n')$. We denote by $\tilde{\pi}$ the non-crossing partition $\tilde{\pi}=\{V_1\cup V_1',\ldots, V_n\cup V_n'\}$ of $\{1,1',\ldots,n,n'\}$. It is in fact the "completion" of $\pi\cup \pi$ (i.e. it is the smallest non-crossing partition that coarsens both $\pi$ as a partition of $\{1,\ldots,n\}$ and $\pi$ as a partition of $\{1',\ldots,n'\}$).

\begin{figure}[h]
   \centering
   \includegraphics{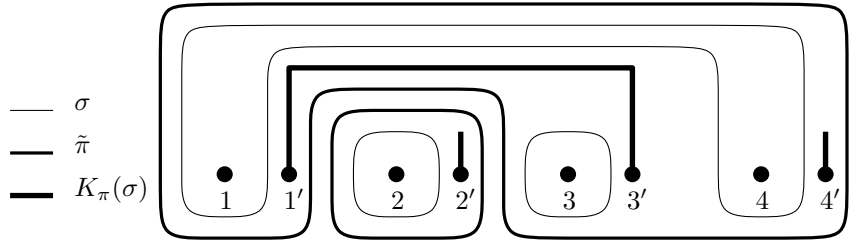} 
   \caption{With $\pi=\{\{1,3,4\},\{2\}\}$ and $\sigma=\{\{1,4\},\{2\},\{3\}\}$, we have $K_{\pi}(\sigma)=\{\{1,3\},\{2\},\{4\}\}$.}
   \label{fig1}
\end{figure}

Now, consider two partitions $\sigma\preceq \pi \in NC(n)$. The partition $K_{\pi}(\sigma)$ is by definition the largest element of $NC(n)$ such that $\sigma\cup K_{\pi}(\sigma)$ is a non-crossing partition of $\{1,1',\ldots,n,n'\}$ and $\sigma\cup K_{\pi}(\sigma)\preceq \tilde{\pi}$. The map $K_{\pi}$ is a bijection from $\{\sigma \in NC(n): \sigma\preceq \pi \}$ onto itself called the Kreweras complementation map with respect to $\pi$ (see Figure~\ref{fig1} for an example). If $\pi=1_n$, we set $K(\sigma)=K_{1_n}(\sigma)$.\label{Kreweras}
\subsubsection{Chains of non-crossing partitions}Let $n\in \mathbb{N}$. A multi-chain in the lattice $NC(n)$ is a tuple of the form $\Gamma = (\pi_0,\ldots,\pi_l)$ with $\pi_0,\ldots,\pi_l\in NC(n) $ such that $ \pi_0 \preceq \pi_1 \preceq \cdots \preceq \pi_l$ (notice that we do not impose $\pi_0=0_n$ or $\pi_l=1_n$, unlike in~\cite{Mastnak2010}). The positive integer $l$ appearing is called the length of the multi-chain, and is denoted by $|\Gamma|$. If $ \pi_0 \neq \pi_1 \neq \cdots \neq \pi_l$, we say that $\Gamma$ is a chain in $NC(n)$. If, for all $1\leq i\leq l$, $K_{\pi_i}(\pi_{i-1})$ has exactly one block which has more than two elements, we say that $\Gamma$ is a simple chain in $NC(n)$.

\begin{figure}[h] 
   \centering
   \includegraphics{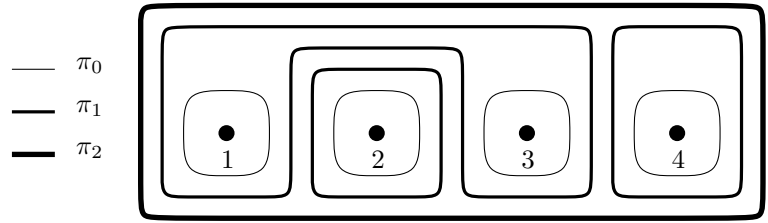} 
   \caption{The chain $(\pi_0,\pi_1,\pi_2)$ with $\pi_0=0_4$, $\pi_1=\{\{1,3\},\{2\},\{4\}\}$ and $\pi_2=1_4$.}
   \label{fig2}
\end{figure}

Let us describe examples. In Figure~\ref{fig2}, it is showed how a multi-chain can be represented graphically as nesting partitions.
It could be helpful to imagine that, in a multi-chain $\Gamma = (\pi_0,\ldots,\pi_l)$, each step from $\pi_{i-1}$ to $\pi_{i}$ is achieved by gluing some blocks of $\pi_{i-1}$.

\begin{figure}[h] 
   \centering
   \includegraphics{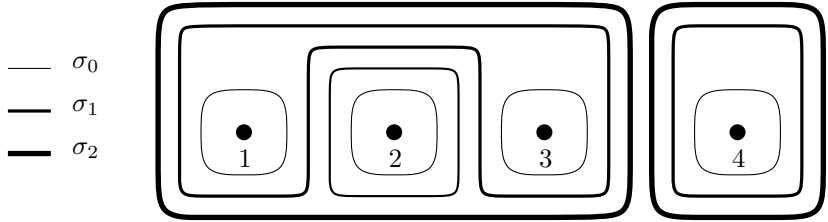} 
   \caption{The simple chain $(\sigma_0,\sigma_1,\sigma_2)$ with $\sigma_0=0_4$, $\sigma_1=\{\{1,3\},\{2\},\{4\}\}$ and $\sigma_2=\{\{1,2,3\},\{4\}\}$.}
   \label{fig3}
\end{figure}

Notice that, for all $1\leq i\leq l$, one can visualize $K_{\pi_i}(\pi_{i-1})$ by putting $n$ dots labelled by $1',\ldots,n'$ at the adequate places (which depend on $i$), just as in Figure~\ref{fig1}. For example, in Figure~\ref{fig2}, $K_{\pi_1}(\pi_{0})=\{\{1,3\},\{2\},\{4\}\}$ and $K_{\pi_2}(\pi_{1})=\{\{1,2\},\{3,4\}\}$. Thus, the multi-chain $(\pi_0,\pi_1,\pi_2)$ is a chain but not a simple chain. A simple chain is obtained by gluing slowly the blocks of $\pi_1$. In Figure~\ref{fig3}, $K_{\sigma_1}(\sigma_{0})=\{\{1,3\},\{2\},\{4\}\}$ and $K_{\sigma_2}(\sigma_{1})=\{\{1,2\},\{3\},\{4\}\}$. Thus, the multi-chain $(\sigma_0,\sigma_1,\sigma_2)$ is a simple chain.
\subsubsection{Free log-cumulants}For all $n\geq 2$, and $A_1,\ldots,A_n \in \mathcal{A}$ such that $\tau(A_1)=\ldots=\tau(A_n)=1$, set
$$L\kappa(A_1,\ldots,A_n)=\displaystyle\sum_{\substack{\Gamma \text{ chain in }NC(n)\\ \Gamma=(\pi_0,\ldots,\pi_{|\Gamma|})\\\pi_0=0_n, \pi_{|\Gamma|}=1_n}} \frac{(-1)^{1+|\Gamma|}}{|\Gamma|} \displaystyle\prod_{i=1}^{|\Gamma|} \kappa\left[K_{\pi_i}(\pi_{i-1})\right]\left(A_1,\ldots,A_n\right).$$
For all $n\geq 2$, and $A_1,\ldots,A_n \in \mathcal{A}$ such that $\tau(A_1),\ldots,\tau(A_n)$ are non-zero, set
$$L\kappa(A_1,\ldots,A_n)=L\kappa\left(\frac{A_1}{\tau(A_1)},\ldots,\frac{A_n}{\tau(A_n)}\right).$$
If $n\geq 2$ and $A_1=\cdots=A_n=A$, we call $L\kappa\left(A_1,\ldots,A_n\right)$ the free log-cumulant of order $n$ of $A$, and we denote it by $L\kappa_n\left(A\right)$. Let us define also $L\kappa_1(A)$, or $L\kappa(A)$, the free log-cumulant of order $1$ of $A$, by $\Log(\tau(A))$, using the principal value $\Log$, the complex logarithm whose imaginary part lies in the interval $(-\pi, \pi]$.

\subsection{The graded Hopf algebra $\mathcal{Y}^{(k)}$} We present here the formalism necessary for a detailed understand of the free log-cumulants and their relation with the free cumulants. We use the notation and results of Mastnak and Nica, and we refer to~\cite{Mastnak2010} for further detailed explanations.

Let $\{1,\ldots,k\}\subset \mathbb{N}$ be a fix index set, and $[k]^*=\cup_{n=0}^\infty\{1,\ldots,k\}^n$ be the set of all words of finite length over the alphabet $\{1,\ldots,k\}$. Let us denote by $\mathcal{Y}^{(k)}$ the commutative algebra of polynomials $\mathbb{C}\Big[Y_w:w\in[k]^*,2\leq |w| \Big]$.
By convention, for a word $w \in [k]^*$ such that $|w| = 1$, set $Y_w=1$. Moreover, for each word $w \in [k]^*$ such that $1\leq |w| = n$, and all $\pi=\{S_1,\ldots,S_k\}\in NC(n)$, set $Y_{w;\pi}=Y_{w_1}\cdots Y_{w_k} \in \mathcal{Y}^{(k)}$ where, for all $1\leq i\leq k$, $w_i$ is the word $w$ restricted to $S_i$.\label{mastnak}

\subsubsection{Definitions}

The comultiplication $\Delta: \mathcal{Y}^{(k)}\rightarrow \mathcal{Y}^{(k)}\otimes \mathcal{Y}^{(k)}$ is the unital algebra homomorphism
uniquely determined by the requirement that, for every $w \in [k]^*$ with $2\leq |w| $, we have $$\Delta(Y_w)=\displaystyle\sum_{\pi\in NC(n)}Y_{w;\pi}\otimes Y_{w;K(\pi)}.$$
The counit $\varepsilon: \mathcal{Y}^{(k)}\mapsto \mathbb{C}$ is the unital algebra homomorphism uniquely determined by the requirement that, for every $w \in [k]^*$ with $2\leq |w| $, we have $\varepsilon(Y_w)=0$.

Proposition 3.6 of~\cite{Mastnak2010} reveals the bialgebra structure of $\mathcal{Y}^{(k)}$. The algebra $\mathcal{Y}^{(k)}$ is a bialgebra when we endow it with the comultiplication $\Delta$ and the counit $\varepsilon$, which means that
\begin{enumerate}
\item $(\Delta\otimes \id)\circ \Delta=(\id\otimes \Delta)\circ \Delta:\mathcal{Y}^{(k)}\rightarrow\mathcal{Y}^{(k)}\otimes \mathcal{Y}^{(k)}\otimes \mathcal{Y}^{(k)}$ (coassociativity),
\item $(\varepsilon\otimes \id)\circ \Delta=(\id\otimes \varepsilon)\circ \Delta=\id$ (counity),
\item $\Delta$ and $\varepsilon$ are unital algebra homomorphisms.
\end{enumerate}
Moreover, $\mathcal{Y}^{(k)}$ is a graded bialgebra. More precisely, for all $n\in \mathbb{N}$, the homogeneous elements of degree $n$ are the elements of $$\mathcal{Y}^{(k)}_n=\text{span} \left\{Y_{w_1}\cdots Y_{w_q}:\left.\begin{array}{l} 1\leq q, w_1,\ldots,w_q\in [k]^*\text{ with }\\|w_1|,\ldots, |w_q|\geq 1\text{ and }|w_1|+\ldots+ |w_q|=n+q\end{array}\right. \right\}.$$
This grading makes the bialgebra $\mathcal{Y}^{(k)}$ a graded connected bialgebra, which means that
\begin{enumerate}
\item $\mathcal{Y}^{(k)}_0=\mathbb{C}$,
\item $\mathcal{Y}^{(k)}=\displaystyle\oplus_{n=0}^{\infty}\mathcal{Y}^{(k)}_n$
\item for all $m,n\in \mathbb{N}$, $\mathcal{Y}^{(k)}_n\cdot \mathcal{Y}^{(k)}_m\subset \mathcal{Y}^{(k)}_{m+n}$ and $\Delta(\mathcal{Y}^{(k)}_n)\subset \displaystyle\oplus_{i=0}^{n}\mathcal{Y}^{(k)}_i\otimes \mathcal{Y}^{(k)}_{n-i}$.
\end{enumerate}

\subsubsection{The convolution}
Let $m:\mathcal{Y}^{(k)}\otimes \mathcal{Y}^{(k)}\rightarrow \mathcal{Y}^{(k)}$ be the linear map given by multiplication in $\mathcal{Y}^{(k)}$. Let $\xi$ and $\eta$ be two endomorphisms of $\mathcal{Y}^{(k)}$. The convolution product $\xi\ast\eta$ is defined by the formula
$$\xi\ast\eta=m\circ(\xi\otimes\eta)\circ\Delta.$$
It is an associative product on endomorphisms. Let us emphasize that a linear functional $\xi:\mathcal{Y}^{(k)}\rightarrow \mathbb{C}$ is considered as an endomorphism of $\mathcal{Y}^{(k)}$ since $\mathbb{C}$ is identified with the subalgebra $\mathbb{C}\cdot 1_{\mathcal{Y}^{(k)}}$ of $\mathcal{Y}^{(k)}$. Let $\xi$ be an endomorphism of $\mathcal{Y}^{(k)}$ such that $\xi(1)=0$. Because $\Delta$ respects the grading of $\mathcal{Y}^{(k)}$, the convolution powers $(\xi^{\ast l})_{l\in \mathbb{N}}$ of $\xi$ are locally nilpotent. More precisely, for all $n\in\mathbb{N}$,  $\xi^{\ast l}$ vanishes on $\mathcal{Y}^{(k)}_n$ whenever $l>n$.
Thus, for any sequence $(\alpha_l)_{l\in \mathbb{N}}$ in $\mathbb{C}$, we can unambiguously define $\sum_{l=0}^\infty\alpha_l\xi^{\ast l}.$ For all $n\in\mathbb{N}$ and all $Y\in \mathcal{Y}^{(k)}$ of degree $n$, we have $$\left(\sum_{l=0}^\infty\alpha_l\xi^{\ast l}\right)Y=\left(\sum_{l=0}^n\alpha_l\xi^{\ast l}\right)Y.$$
In particular, for any endomorphism $\xi$ of $\mathcal{Y}^{(k)}$ such that $\xi(1)=0$, let us define $\exp_{\ast}(\xi)$ by $\sum_{l=0}^\infty\frac{1}{l!}\xi^{\ast l}$, and for any endomorphism $\eta$ of $\mathcal{Y}^{(k)}$ such that $\eta(1)=1$, let us define $\log_{\ast}(\eta)$ by $-\sum_{l=0}^\infty\frac{1}{l}(\varepsilon-\eta)^{\ast l}$. The exponentiation $\exp_{\ast} $ maps the set of endomorphims $\xi$ of $\mathcal{Y}^{(k)}$ such that $\xi(1)=0$ bijectively onto the set of endomorphisms $\eta$ of $\mathcal{Y}^{(k)}$ such that $\eta(1)=1$, and the inverse of this bijection is given by $\log_{\ast}$.

Let us indicate that $\mathcal{Y}^{(k)}$ is a Hopf algebra. Indeed, let us define the antipode $S$ by the series $\varepsilon+\sum_{l=1}^\infty(\varepsilon-\id)^{\ast l}$. Then the antipode $S$ is such that $S\ast \id = \varepsilon = \id\ast S$.

\subsubsection{Free cumulants}\label{nica}
Let $\mathbf{A}=(A_1,\ldots,A_k)\in \mathcal{A}^k$. The character $\chi_\mathbf{A}:\mathcal{Y}^{(k)}\rightarrow \mathbb{C}$ associated to $\mathbf{A}$ is defined as follows (note that in~\cite{Mastnak2010}, $\chi_\mathbf{A}$ would be denoted by $\chi_{\mu_\mathbf{A}}$ where $\mu_\mathbf{A}$ is the distribution of $\mathbf{A}$).
The linear functional $\chi_\mathbf{A}$ is multiplicative, and for all $w=(i(1),\ldots,i(n))\in [k]^*$ such that $2\leq n $, we have $\chi_\mathbf{A}(Y_w)=\kappa(A_{i(1)},\ldots, A_{i(n)})$.

Let us suppose that $\kappa(A_1)=\ldots=\kappa(A_k)=1$. For all word $w=(i(1),\ldots,i(n)) \in [k]^*$ such that $1\leq  n$, and all $\pi\in NC(n)$, we have\begin{equation}
\chi_\mathbf{A}(Y_{w;\pi}) =\kappa[\pi]\left(A_{i(1)},\ldots, A_{i(n)}\right)
\end{equation}
by the multiplicative properties of $\chi_\mathbf{A}$ and of the free cumulants.
Proposition~4.5 of~\cite{Mastnak2010} links the free cumulants and the free log-cumulants as follows. For all $w=(i(1),\ldots,i(n))\in [k]^*$ such that $2\leq n $, we have \begin{equation}
\log_{\ast}\chi_\mathbf{A}(Y_w)=L\kappa\left(A_{i(1)},\ldots, A_{i(n)}\right).\label{chilog}
\end{equation}
Moreover, $\log_{\ast}\chi_\mathbf{A}$ is an infinitesimal character, which means that, for all $Y_1,Y_2\in \mathcal{Y}^{(k)}$, we have
$$\log_{\ast}\chi_\mathbf{A}(Y_1Y_2)=\log_{\ast}\chi_\mathbf{A}(Y_1)\cdot \varepsilon(Y_2)+\varepsilon(Y_1)\cdot \log_{\ast}\chi_\mathbf{A}(Y_2) .$$

Let $w=(i(1),\ldots,i(n))\in [k]^*$ be such that $2\leq n $. Let $\pi \in NC(n)$ be such that $\pi$ has exactly one block which has at least two elements. Let $ \{j_1,\ldots ,j_N\}$ be this block of $\pi$, with $j_1 < j_2 < \ldots < j_N$. Let us denote by $L\kappa\left[\pi\right](A_{i(1)},\ldots,A_{i(n)})$ the free log-cumulant $L\kappa(A_{i(j_1)},\ldots,A_{i(j_N)})$. 
\begin{lemma}
Let $(A_1,\ldots,A_k)\in \mathcal{A}^k$ be such that $\kappa(A_1)=\ldots=\kappa(A_k)=1$. Let $w=(i(1),\ldots,i(n))\in [k]^*$ be such that $2\leq n $ and $\pi\in NC(n)$. We have
\begin{equation*}\kappa[\pi]\left(A_{i(1)},\ldots,A_{i(n)}\right)=\displaystyle\sum_{\substack{\Gamma \text{ simple chain in }NC(n)\\ \Gamma=(\pi_0,\ldots,\pi_{|\Gamma|})\\\pi_0=0_n, \pi_{|\Gamma|}=\pi}} \frac{1}{|\Gamma| !} \displaystyle\prod_{i=1}^{|\Gamma|} L\kappa\Big[K_{\pi_i}(\pi_{i-1})\Big]\left(A_{i(1)},\ldots,A_{i(n)}\right).\end{equation*}\label{cumlogcumptwo}
\end{lemma}
\begin{proof}
We have $\kappa[\pi](A_{i(1)},\ldots, A_{i(n)})=\chi_\mathbf{A}(Y_{w;\pi})=\exp_{\ast}(\log_{\ast}\chi_\mathbf{A})(Y_{w;\pi})$.
By definition of $\exp_{\ast}$, we have $$\kappa[\pi]\left(A_{i(1)},\ldots, A_{i(n)}\right)=\displaystyle\sum_{l=0}^\infty\frac{1}{l!}(\log_{\ast}\chi_\mathbf{A})^{\ast l}(Y_{w;\pi}).$$
For all $1\leq l$, let $\Delta^l:\mathcal{Y}^{(k)}\rightarrow(\mathcal{Y}^{(k)})^{\otimes l}$ denote the iterate of $\Delta$. Following step-by-step the proof of Proposition~4.2 in~\cite{Mastnak2010} but with an arbitrary endpoint $\pi_{l}$, we have
$$\Delta^l(Y_{w;\pi})=\displaystyle\sum_{\substack{\Gamma \text{ multi-chain in }NC(n)\\ \Gamma=(\pi_0,\ldots,\pi_{l})\\ \pi_{l}=\pi}}  Y_{w;K_{\pi_1}(\pi_{0})}\otimes Y_{w;K_{\pi_2}(\pi_{1})}\otimes \cdots\otimes Y_{w;K_{\pi_l}(\pi_{l-1})} .$$
Thus, we have
\begin{eqnarray*}\kappa[\pi]\left(A_{i(1)},\ldots, A_{i(n)}\right)&=&\displaystyle\sum_{l=0}^\infty\frac{1}{l!}\displaystyle\sum_{\substack{\Gamma \text{ multi-chain in }NC(n)\\ \Gamma=(\pi_0,\ldots,\pi_{l})\\ \pi_0=0_n,\pi_{l}=\pi}}\displaystyle\prod_{i=1}^{|\Gamma|} \log_{\ast}\chi_\mathbf{A}(Y_{w;K_{\pi_i}(\pi_{i-1})})\\
&=&\displaystyle\sum_{\substack{\Gamma \text{ multi-chain in }NC(n)\\ \Gamma=(\pi_0,\ldots,\pi_{|\Gamma|})\\\pi_0=0_n, \pi_{|\Gamma|}=\pi}} \frac{1}{|\Gamma| !} \displaystyle\prod_{i=1}^{|\Gamma|} \log_{\ast}\chi_\mathbf{A}(Y_{w;K_{\pi_i}(\pi_{i-1})}).
\end{eqnarray*}
It remains to prove that, for all $\pi\in NC(n)$, $\log_{\ast}\chi_\mathbf{A}(Y_{w;\pi})=L\kappa[\pi](A_{i(1)},\ldots,A_{i(n)})$ if $\pi$ has exactly one block which has more than two elements, and $0$ otherwise.
Let $\pi\in NC(n)$ be such that $\pi$ has exactly one block which has at least two elements. Let $\{j_1,\ldots, j_N\}$ be this block of $\pi$, with $j_1 < j_2 < \ldots <j_N$. Using \eqref{chilog}, we have $$\log_{\ast}\chi_\mathbf{A}(Y_{w;\pi})=\log_{\ast}\chi_\mathbf{A}\left(Y_{(i(j_1),\ldots,i(j_N))}\right)=L\kappa\left(A_{i(j_1)},\ldots, A_{i(j_N)}\right)=L\kappa[\pi]\left(A_{i(1)},\ldots,A_{i(n)}\right).$$
If $\pi=0_n \in NC(n)$, we have $\log_{\ast}\chi_\mathbf{A}(Y_{w;\pi})=\log_{\ast}\chi_\mathbf{A}(1)=0$.
Finally, if $\pi\in NC(n)$ has two blocks which have at least two elements, there exist $Y_1,Y_2\in \mathcal{Y}^{(k)}\setminus \mathcal{Y}^{(k)}_0$ such that $Y_{w;\pi}=Y_1Y_2$. We have $\varepsilon(Y_1)=\varepsilon(Y_2)=0$. Thus, 
\begin{equation*}
\log_{\ast}\chi_\mathbf{A}(Y_{w;\pi})=\log_{\ast}\chi_\mathbf{A}(Y_1)\varepsilon(Y_2)+\varepsilon(Y_1)\log_{\ast}\chi_\mathbf{A}(Y_2)=0.\qedhere
\end{equation*}
\end{proof}

\subsection{Properties of the free log-cumulants}This section presents some important properties of the free log-cumulants which arise from the Hopf algebra structure presented in the previous section. However, in their formulations, the different results do not involve the understanding of the previous section.\label{frprop}

\subsubsection{Cumulants from log-cumulants}

Let $n\geq 2$, and $A_1,\ldots,A_n \in \mathcal{A}$ be such that $\tau(A_1),\ldots,\tau(A_n)$ are non-zero. Let $\pi \in NC(n)$ be such that $\pi$ has exactly one block which has at least two elements. Let $ \{i_1,... ,i_k\}$ be this block of $\pi$, with $i_1 < i_2 < \ldots < i_k$. Let us denote by $L\kappa\left[\pi\right]\left(A_1,\ldots,A_n\right)$ the free log-cumulant $L\kappa(A_{i_1},\ldots,A_{i_k})$. 
\begin{proposition}
Let $n\in \mathbb{N}^*$, and $A_1,\ldots,A_n \in \mathcal{A}$ be such that $\tau(A_1),\ldots,\tau(A_n)$ are non-zero. For all $\pi \in NC(n)$, we have
\begin{equation}\kappa[\pi](A_1,\ldots,A_n)=e^{L\kappa(A_1)}\cdots e^{L\kappa(A_n)}\hspace{-1,3cm}\displaystyle\sum_{\substack{\Gamma \text{ simple chain in }NC(n)\\ \Gamma=(\pi_0,\ldots,\pi_{|\Gamma|})\\\pi_0=0_n, \pi_{|\Gamma|}=\pi}} \frac{1}{|\Gamma| !} \displaystyle\prod_{i=1}^{|\Gamma|} L\kappa\Big[K_{\pi_i}(\pi_{i-1})\Big]\left(A_1,\ldots,A_n\right).\label{cumlogcum}\end{equation}\label{cumlogcump}
\end{proposition}
\begin{proof}The degenerate case where $n=1$ is verified because in this case, $\pi=\{\{1\}\}$,$$\kappa[\pi](A_1)=\tau(A_1)=e^{L\kappa(A_1)}\text{ and }\displaystyle\sum_{\substack{\Gamma \text{ simple chain in }NC(1)\\ \Gamma=(\pi_0,\ldots,\pi_{|\Gamma|})\\\pi_0=0_n, \pi_{|\Gamma|}=\pi}} \frac{1}{|\Gamma| !} \displaystyle\prod_{i=1}^{|\Gamma|} L\kappa\Big[K_{\pi_i}(\pi_{i-1})\Big]\left(A_1\right)=1.$$
Let us suppose that $2\leq n$. Since $\kappa[\pi](A_1,\ldots,A_n)=e^{L\kappa(A_1)}\cdots e^{L\kappa(A_n)}\kappa[\pi](A_1/\tau(A_1),\ldots,A_n/\tau(A_n))$, and the free log-cumulants of $\left(A_1,\ldots,A_n\right)$ involved are equal to those of $(A_1/\tau(A_1),\ldots,A_n/\tau(A_n))$, we can assume that $\tau(A_1)=\ldots=\tau(A_n)=1$, or equivalently that $\kappa(A_1)=\ldots=\kappa(A_n)=1$. We end the proof using Lemma~\ref{cumlogcumptwo}.\end{proof}

From \eqref{taucum} and \eqref{cumlogcum}, we deduce the following corollary.
\begin{corollary}
Let $n\in \mathbb{N}^*$, and $A_1,\ldots,A_n \in \mathcal{A}$ be such that $\tau(A_1),\ldots,\tau(A_n)$ are non zero. We have
\begin{equation}
\tau(A_1\cdots A_n)=e^{L\kappa(A_1)}\cdots e^{L\kappa(A_n)}\hspace{-1,3cm}\displaystyle\sum_{\substack{\Gamma \text{ simple chain in }NC(k)\\ \Gamma=(\pi_0,\ldots,\pi_{|\Gamma|}),\pi_0=0_n}} \frac{1}{|\Gamma| !} \displaystyle\prod_{i=1}^{|\Gamma|} L\kappa\Big[K_{\pi_i}(\pi_{i-1})\Big]\left(A_1,\ldots,A_n\right).\label{taulogcum}
\end{equation}
\end{corollary}

\subsubsection{Log-cumulants and freeness}The $LS$-transform of $\mathbf{A}$ is the series in the non-commuting variables $z_1,\ldots,z_k$ defined by
\begin{eqnarray*}
LS_{\mathbf{A}}(z_1,\ldots, z_k) &=& \displaystyle\sum_{\substack{2\leq n\\ w=(i(1),\ldots,i(n))\in [k]^*}}\log_{\ast} \chi_\mathbf{A}(Y_w) \cdot z_{i(1)}\cdots z_{i(n)}\\
&=&\displaystyle\sum_{\substack{2\leq n\\ (i(1),\ldots,i(n))\in [k]^*}} L\kappa\left(A_{i(1)},\ldots,A_{i(n)}\right)z_{i(1)}\cdots z_{i(n)}.
\end{eqnarray*}

The following two propositions are results analogous to Proposition~\ref{freeness} and Corollary~\ref{addfree}. They are reformulations of Proposition 5.4 and Corollary 1.5 of~\cite{Mastnak2010}.
\begin{proposition}
Let $\left(\mathcal{B}_i \right)_{i\in I}$ be subalgebras of $\mathcal{A}$. They are free if and only if their mixed free log-cumulants vanish. That is to say: for all $n\in \mathbb{N}^*$, all $i_1,\ldots, i_n \in I$ and all $A_1, \ldots , A_n \in \mathcal{A}$ such that $\tau(A_1)=\ldots=\tau(A_n)=1$ and such that $A_j$ belongs to some $\mathcal{B}_{i_j}$ for all $1\leq j \leq n$, whenever there exists some $j$ and $j' $ with $i_j\neq i_{j'}$, we have $L\kappa(A_1,\ldots ,A_n) = 0$.\label{logfreeness}
\end{proposition}
\begin{proof}
Proposition 5.4 of~\cite{Mastnak2010} says that $A_1, \ldots ,A_k$ are free if and only if one has $LS_{\mathbf{A}}(z_1,\ldots z_k)=LS_{A_1}(z_1)+\ldots+LS_{A_k}(z_k)$, or equivalently, that $A_1, \ldots ,A_k$ are free if and only if $L\kappa(A_{i(1)},\ldots,A_{i(n)})=0$ each time there exists some $j$ and $j'$ such that $i_j\neq i_{j'}$.
Proposition \ref{logfreeness} follows immediately.
\end{proof}
\begin{proposition}
Let $A$ and $B\in \mathcal{A}$ be free.\label{addlogfree}

We have $L\kappa_1(AB)\equiv L\kappa_1(A)+L\kappa_1(B) \pmod{2i\pi}$, and for all $ n\geq 2$: $$L\kappa_n(AB)=L\kappa_n(A)+L\kappa_n(B) .$$
\end{proposition}
\begin{proof}Let $A$ and $B\in \mathcal{A}$ be free and such that $\tau(A)$ and $\tau(B)$ are non-zero. We have first $\tau(AB)=\tau(A)\tau(B)$, thus $L\kappa_1(AB)\equiv L\kappa_1(A)+L\kappa_1(B) \pmod{2i\pi}$. Set $\tilde{A}=A/\tau(A)$ and $\tilde{B}=B/\tau(B)$. Corollary 1.5 of~\cite{Mastnak2010} says that $LS_{\tilde{A}\tilde{B}}=LS_{\tilde{A}}+LS_{\tilde{B}}$. Thus, for all $2\leq n$, $
L\kappa_n\left(AB\right)=L\kappa_n(\tilde{A}\tilde{B})=L\kappa_n(\tilde{A})+L\kappa_n(\tilde{B})=L\kappa_n\left(A\right)+L\kappa_n\left(B\right).$
\end{proof}

\subsubsection{$S$-transform}In the $1$-dimensional case, the free log-cumulants can be recovered from the $S$-transform. Indeed, let $A\in \mathcal{A}$ be such that $\tau(A)=1$. Let us consider the $R$-transform of $A$, i.e. the formal series $R_A(z)=\sum_{n=1}^\infty \kappa_n(A) z^n$. Let $S_A$ be the $S$-transform of $A$: it is the formal series $S_A$ such that $zS_A(z)$ is the inverse under composition of $R_A(z)$.
We remark that $\kappa_1(A)=\tau(A)=1$, and by consequence, we have $S_A(0)=1/\kappa_1(A)=1$. Thus we can define the formal logarithm of $S_A$ as the formal series
$\log S_A (z) =-\sum_{n=1}^\infty \frac{1}{n} (1-S_A(z))^n$. Corollary 6.12 of~\cite{Mastnak2010} establishes  then $-z \log S_A (z) = \sum_{n=2}^\infty L\kappa_n(A) z^n=LS_{A}(z)$.

\subsubsection{Free unitary Brownian motion} \label{logcumut}Let $t\geq 0$. A non-commutative random variable $U_t$ is a free unitary Brownian motion in distribution at time $t$ if $U_t$ is unitary and the free log-cumulants of $U_t$ are $L\kappa_1(U_t)=-t/2$, $L\kappa_2(U_t)=-t$ and $L\kappa_n(U_t)=0$ for all $n>2$.\label{fu} The distribution of $U_t$ is given by \eqref{taulogcum}: for all $n\in \mathbb{N}^*$, we have
$$\tau\left(U_t^n\right)=e^{-\frac{nt}{2}}\displaystyle\sum_{\substack{\Gamma \text{ simple chain in }NC(n)\\ \Gamma=(\pi_0,\ldots,\pi_{|\Gamma|}),\pi_0=0_n}} \frac{1}{|\Gamma| !} \displaystyle\prod_{i=1}^{|\Gamma|} L\kappa\Big[K_{\pi_i}(\pi_{i-1})\Big]\left(U_t\right).$$
A simple chain $\Gamma = (0_n,\pi_1,\ldots,\pi_{|\Gamma|})$ in $NC(n)$ is called an increasing path if, for all $1\leq i\leq l$, the block of $K_{\pi_i}(\pi_{i-1})$ which has more than two elements has exactly two elements. Proposition 6.6 of~\cite{LEVY2008} tells us that, for all $k\geq0$, the number of increasing paths of length $k$ in $NC(n)$ is exactly $\binom{n}{k+1}n^{k-1}$ if $k\leq n-1$ and $0$ if $k\geq n$. Thus,
$$\tau\left(U_t^n\right)=e^{-\frac{nt}{2}}\displaystyle\sum_{k=0}^{n-1}\left( \displaystyle\sum_{\substack{\Gamma \text{ increasing path in }NC(n)\\ |\Gamma|=k}} \frac{(-t)^k}{k !}\right)=e^{-\frac{nt}{2}}\displaystyle\sum_{k=0}^{n-1}\frac{(-t)^k}{k !}\binom{n}{k+1}n^{k-1}.$$
In~\cite{Biane1997a}, Biane proved that it is indeed the distribution of a free unitary Brownian motion $(U_t)_{t\geq 0}$ at time $t$ as defined in Section~\ref{fu}.

\subsection{Multiplicative transition operators}\label{dut}
Let $\mathbf{A}=(A_i)_{i\in I}\in \mathcal{A}^I$. Let us define a derivation $\D_{\mathbf{A}}$ associated to $\mathbf{A}$ on $\left(\mathbb{C}\{X_i:i\in I\},\cdot_{\tr}\right)$ in the following way. 
For all $n\in \mathbb{N}$ and $ i(1),\ldots,i(n)\in I$, we set
$$\begin{array}{rccl}
\D_{\mathbf{A}} \left( X_{i(1)}\cdots X_{i(n)}\right)&=& \vspace{-0,5cm}  \displaystyle\sum_{\substack{1\leq m\leq n\\ 1\leq k_1< \ldots < k_m\leq n}} & L\kappa\left(A_{i(k_1)},\ldots,A_{i(k_m)}\right)X_{i(1)}\cdots X_{i(k_1-1)}\\
 & & & \hspace{1cm}\cdot \tr\left( X_{i(k_1)}\cdots X_{i(k_2-1)} \right)\\
  & & &\hspace{1cm}\cdot \tr\left( X_{i(k_2)}\cdots X_{i(k_3-1)} \right)\\
  & & &\hspace{1cm}\cdots\\
& & &\hspace{1cm}\cdot \tr\left( X_{i(k_{m-1})}\cdots X_{i(k_m-1)} \right)\\
& & &\hspace{1cm} \cdot  X_{i(k_m)}\cdots X_{i(n)}
\end{array}$$
and we extend $\D_{\mathbf{A}}$ to all $\mathbb{C}\{X_i:i\in I\}$ by linearity and by the relation of derivation
$$\forall P,Q \in \mathbb{C}\{X_i:i\in I\},\ \D_{\mathbf{A}} \left(P\tr Q\right)=\left(\D_{\mathbf{A}} P\right)\tr Q+P\tr \left(\D_{\mathbf{A}} Q\right).$$
For any finite index set $J\subset I$ and $d\in \mathbb{N}$, the finite-dimensional space $\mathbb{C}_d\{X_i:i\in J\}$ is invariant for the operator $\D_{\mathbf{A}}$. Thus, we can define $e^{\D_{\mathbf{A}}}$ on each of those spaces. The operator $e^{\D_{\mathbf{A}}}$ on $\mathbb{C}\{X_i:i\in I\}$ is defined by the series $\sum_{k=0}^{\infty}\frac{1}{k!} \D_{\mathbf{A}}^k.$
For all $P\in \mathbb{C}_d\{X_i:i\in J\}$, $e^{\D_{\mathbf{A}}} P$ is the convergent sum $\sum_{k=0}^{\infty} \frac{1}{k!}  \D_{\mathbf{A}}^k P$.
The operator $\D_{\mathbf{A}}$ is a derivation, and we have the Leibniz formula
$$\forall k\in \mathbb{N}, \forall P,Q \in \mathbb{C}\{X_i: i\in I\},\left(  \D_{\mathbf{A}}\right)^k \left(P\tr Q\right) =\displaystyle\sum_{l=0}^k \binom{k}{l} \left(  \D_{\mathbf{A}}^l P \right)\tr \left(  \D_{\mathbf{A}}^{k-l} Q\right),$$
from which we deduce, using the standard power series argument, that the operator $e^{ \D_{\mathbf{A}}}$ is multiplicative in the following sense:
$$\forall P,Q \in \mathbb{C}\{X_i: i\in I\},\ e^{ \D_{\mathbf{A}}} \left(P\tr Q\right)=\left(e^{ \D_{\mathbf{A}}} P\right)\tr \left(e^{\D_{\mathbf{A}}}Q\right).$$
Let us denote by $\tau(\mathbf{A})$ the family $(\tau(A_i))_{i\in I}$ and by $\mathbf{A}/\tau(\mathbf{A})$ the family $(A_i/\tau(A_i))_{i\in I}$.
\begin{proposition}Let $\mathbf{A}=\left(A_i\right)_{i\in I}$ and $\mathbf{B}=\left(B_i\right)_{i\in I}\in \mathcal{A}^I$ be such that $\tau(A_i)\neq 0$ and $\tau(B_i)\neq 0$ for all $i\in I$. The operators\label{propcom}
$\D_{\tau(\mathbf{A})}$, $\D_{\tau(\mathbf{B})}$ and $\D_{\mathbf{A}/\tau(\mathbf{A})}$ commute, $\D_{\mathbf{A}}=\D_{\tau(\mathbf{A})}+\D_{\mathbf{A}/\tau(\mathbf{A})}$,  and $e^{\D_{\mathbf{A}}}=e^{\D_{\mathbf{A}/\tau(\mathbf{A})}}e^{\D_{\tau(\mathbf{A})}}$.

\end{proposition}
\begin{proof}The free log-cumulants of $\tau(\mathbf{A})$ are zero except $L\kappa(\tau(A_i))=L\kappa(A_i)$ for all $i\in I$ and the free log-cumulants of $\mathbf{A}/\tau(\mathbf{A})$ are those of $\mathbf{A}$, except $L\kappa(A_i/\tau(A_i))=0$ for all $i\in I$.

We recall that, from Section~\ref{Anotprod}, the monomials of $ \mathbb{C}\langle X_i:i\in I\rangle$ generate the algebra $(\mathbb{C}\{X_i:i\in I\},\cdot_{\tr} )$. Let  $n\in \mathbb{N}$ and $ i(1),\ldots,i(n)\in I$ such that $M=X_{i(1)}\cdots X_{i(n)}$. For all $1\leq k<l\leq n$, we denote by $\tr_{k,l}$ the element $\tr( X_{i(k)}\cdots X_{i(l-1)})$. We compute
$$\begin{array}{rcl}
& &\hspace{-1,5cm}\D_{\tau(\mathbf{B})}\D_{\mathbf{A}/\tau(\mathbf{A})} \left( X_{i(1)}\cdots X_{i(n)}\right)\\
&=&\vspace{-0,3cm}\displaystyle\sum_{1\leq l \leq n} \displaystyle\sum_{\substack{2\leq m\leq n\\ 1\leq k_1< \ldots < k_m\leq n}}  L\kappa\left(B_{i(l)}\right)L\kappa\left(A_{i(k_1)},\ldots,A_{i(k_m)}\right)\\
& &\hspace{4cm} \cdot X_{i(1)}\cdots X_{i(k_1-1)} \cdot \tr_{k_1,k_2}\cdots   \tr_{k_{m-1},k_m} \cdot  X_{i(k_m)}\cdots X_{i(n)}\\
&=&\D_{\mathbf{A}/\tau(\mathbf{A})} \D_{\tau(\mathbf{B})}\left( X_{i(1)}\cdots X_{i(n)}\right).
\end{array}$$
Thus, $\D_{\tau(\mathbf{B})}\D_{\mathbf{A}/\tau(\mathbf{A})}M=\D_{\mathbf{A}/\tau(\mathbf{A})}\D_{\tau(\mathbf{B})}M$ for all monomials $M\in \mathbb{C}\langle X_i: i\in I\rangle$, and we extend the commutativity on all $\mathbb{C}\{X_i: i\in I\}$ by induction because $\D_{\tau(\mathbf{B})}$ and $\D_{\mathbf{A}/\tau(\mathbf{A})}$ are derivations. Similarly, we verify the commutativity of $\D_{\tau(\mathbf{A})}$, $\D_{\tau(\mathbf{B})}$ and $\D_{\mathbf{A}/\tau(\mathbf{A})}$ on monomials of $P\in \mathbb{C}\langle X_i: i\in I\rangle$, and we extend it on all $\mathbb{C}\{X_i: i\in I\}$ by induction.

Finally, the operators $\D_{\mathbf{A}}$ and $\D_{\tau(\mathbf{A})}+\D_{\mathbf{A}/\tau(\mathbf{A})}$ are two derivations which coincide on monomials, so $\D_{\mathbf{A}}=\D_{\tau(\mathbf{A})}+\D_{\mathbf{A}/\tau(\mathbf{A})} $, and $e^{\D_{\mathbf{A}}}=e^{\D_{\mathbf{A}/\tau(\mathbf{A})}}e^{\D_{\tau(\mathbf{A})}}$ is a direct consequence of the two first assertions.
\end{proof}
For all $\mathbf{A}=(A_i)_{i\in I}\in \mathcal{A}^I$ and $\mathbf{B}=(B_i)_{i\in I}\in \mathcal{A}^I$, let us denote $(A_iB_i)_{i\in I}\in \mathcal{A}^I$ by $\mathbf{A}\mathbf{B}$.
\begin{theorem}
Let $I$ be an arbitrary index set. Let $\mathbf{A}=(A_i)_{i\in I}\in \mathcal{A}^I$ be such that $\tau(A_i)\neq 0$ for all $i\in I$. For all $P\in \mathbb{C}\{ X_i:i\in I\}$, and all $\mathbf{B}=(B_i)_{i\in I}\in \mathcal{A}^I$ free from $(A_i)_{i\in I}$ and such that $\tau(B_i)\neq 0$ for all $i\in I$, we have\label{freekernelmult}
$$\tau\left(\left.P\left(\mathbf{A}\mathbf{B}\right)\right|\mathbf{B}\right)=(e^{\D_{\mathbf{A}}}P)\left(\mathbf{B}\right).$$
In particular, for all $P\in \mathbb{C}\{ X_i:i\in I\}$, we have
\begin{equation*} \tau\left(P\left(\mathbf{A}\right)\right)=(e^{\D_{\mathbf{A}}}P)\left(1\right).\end{equation*}
\end{theorem}

Note that, from Proposition~\ref{propcom} and Theorem~\ref{freekernelmult}, we have for all $P\in \mathbb{C}\{ X_i:i\in I\}$
\begin{eqnarray*}\tau\left(P\left(\mathbf{A}\right)\right)=(e^{\D_{\mathbf{A}}}P)\left(1\right)=(e^{\D_{\mathbf{A}/\tau(\mathbf{A})}}(e^{\D_{\tau(\mathbf{A})}}P))(1).\end{eqnarray*}
In fact, the last expression $e^{\D_{\mathbf{A}/\tau(\mathbf{A})}}e^{\D_{\tau(\mathbf{A})}}P$ is easier to calculate because in practice $\D_{\mathbf{A}/\tau(\mathbf{A})}$ is locally nilpotent, i.e. it is nilpotent on any finite dimensional space, and $ \D_{\tau(\mathbf{A})}$ simply multiplies each element of the canonical basis of $P\in \mathbb{C}\{ X_i:i\in I\}$ by a factor which depends on its degree (see the example in the next section).

\begin{proof}[Proof of Theorem~\ref{freekernelmult}]\label{freekernelmultproof}We start by proving a weak version of Theorem~\ref{freekernelmult} in the following proposition.
\label{sssmone}
\begin{proposition}\label{freekernelmultweak}
Let $I$ be an arbitrary index set. Let $\mathbf{A}=\left(A_i\right)_{i\in I}\in \mathcal{A}^I$. Let $P\in \mathbb{C}\{ X_i:i\in I\}$, and $\mathbf{B}=\left(B_i\right)_{i\in I}\in \mathcal{A}^I$ be free from $\left(A_i\right)_{i\in I}$. Let us assume that, for all $i\in I$, $\tau(A_i)=\tau(B_i)=1$. We have$$ \tau\left(P\left(\mathbf{A}\mathbf{B}\right)\right)=\tau\left((e^{\D_{\mathbf{A}}}P)\left(\mathbf{B}\right)\right).$$
\end{proposition}
\begin{proof}For any $P\in \mathbb{C}\{ X_i:i\in I\}$, $P$ depends on only finitely-many indices. Thus, we can suppose that $I$ is finite. Let us say that $I=\{1,\ldots,k\}\subset \mathbb{N}$.

Let us assume first that $P$ is a monomial. Let $w=(i(1),\ldots,i(n))\in [k]^*$ be such that $P=X_{i(1)}\cdots X_{i(n)}$.
Thanks to \eqref{taucum}, we have
$$
\tau\left(A_{i(1)}B_{i(1)} \cdots A_{i(n)}B_{i(n)}\right)=\displaystyle\sum_{\pi\in NC(2n)}\kappa\left[\pi\right]\left(A_{i(1)},B_{i(1)}, \cdots, A_{i(n)},B_{i(n)}\right)
.$$
The elements $\mathbf{A}$ and $\mathbf{B}$ are free. By Proposition~\ref{freeness}, the only elements of $NC(2n)$ which contribute to the sum are the partitions of the form $\pi_1 \cup \pi_2$ where $\pi_1 \in NC(\{1,3,\ldots,2n-1\})$ and $\pi_2 \in NC(\{2,4,\ldots,2n\})$. We use now the definition of $K(\pi)$. We have
\begin{eqnarray*}
& &\hspace{-2cm}\tau\left(A_{i(1)}B_{i(1)} \cdots A_{i(n)}B_{i(n)}\right)\\
&=&\displaystyle\sum_{\pi_1\in NC(\{1,3,\ldots,2n-1\}) }\displaystyle\sum_{\substack{\pi_2 \in NC(\{2,4,\ldots,2n\})\\ \pi_1 \cup \pi_2 \in NC(2n)}} \kappa\left[\pi_1 \cup \pi_2\right]\left(A_{i(1)},B_{i(1)} ,\ldots,A_{i(n)},B_{i(n)}\right) \\
&=&\displaystyle\sum_{\pi_1 \in NC(n) }\displaystyle\sum_{\substack{\pi_2 \in NC(n)\\ \pi_2 \preceq K(\pi_1) \in NC(n)}} \kappa\left[\pi_1\right] \left(A_{i(1)} ,\ldots, A_{i(n)}\right)\cdot \kappa\left[\pi_2\right] \left(B_{i(1)}  ,\ldots, B_{i(n)}\right).
\end{eqnarray*}
Applying now \eqref{taucum}, we see that
$$\tau\left(A_{i(1)}B_{i(1)} \cdots A_{i(n)}B_{i(n)}\right)=\displaystyle\sum_{\pi \in NC(n) }\kappa\left[\pi \right]\left(A_{i(1)} ,\ldots, A_{i(n)}\right)\cdot \tau\Big[K(\pi)\Big] \left(B_{i(1)}  ,\ldots, B_{i(n)}\right).$$
We recognize here the comultiplication $\Delta$. More precisely, let us define $\eta_\mathbf{B}$ to be the character associated to $\mathbf{B}$ as below. The linear functional $\eta_\mathbf{B}$ is multiplicative, and for all $w=(i(1),\ldots,i(n))\in [k]^*$, we have $\eta_\mathbf{B}(Y_w)=\tau(B_{i(1)},\ldots, B_{i(n)})$.
For all $w=(i(1),\ldots,i(n))\in [k]^*$ we have
\begin{eqnarray*}\tau\left(A_{i(1)}B_{i(1)} \cdots A_{i(n)}B_{i(n)}\right)&=&\displaystyle\sum_{\pi \in NC(n) }\kappa\left[\pi \right]\left(A_{i(1)} ,\ldots, A_{i(n)}\right)\cdot \tau\Big[K(\pi)\Big] \left(B_{i(1)}  ,\ldots, B_{i(n)}\right)\\
&=&\displaystyle\sum_{\pi \in NC(n) }\chi_\mathbf{A}(Y_{w;\pi})\eta_\mathbf{B}(Y_{w;K(\pi)})\\
&=&(\chi_\mathbf{A}\ast \eta_\mathbf{B})(Y_{w}).
\end{eqnarray*}
We arrive in the Hopf algebra $\mathcal{Y}^{(k)}$, and consequently we introduce the algebra homomorphism $\rho: \left(\mathbb{C}\{X_i:i\in I\},\cdot_{\tr}\right),\mapsto (\mathcal{Y}^{(k)},\cdot)$. For all monomials $X_{i(1)}\cdots X_{i(n)}\in \mathbb{C}\{X_i:i\in I\}$, we set $\rho (X_{i(1)}\cdots X_{i(n)})=Y_{(i(1),\ldots,i(n))}$, and we extend $\rho$ by linearity and by products.

Since $P\mapsto \tau\left(P\left(\mathbf{A}\mathbf{B}\right)\right)$ is multiplicative and $\chi_\mathbf{A}\ast \eta_\mathbf{B}$ is also multiplicative as a convolution of two characters, the relation $\tau\left(P\left(\mathbf{A}\mathbf{B}\right)\right)=(\chi_\mathbf{A}\ast \eta_\mathbf{B})(\rho(P))$ extends from monomials to all $\mathbb{C}\{X_i:i\in I\}$. Similarly, the relation $\tau\left(P\left(\mathbf{B}\right)\right)= \eta_\mathbf{B}(\rho(P))$ extends from monomials to all $\mathbb{C}\{X_i:i\in I\}$.

Here we prove a general result which relates the composition exponentiation with the convolution exponentiation $\exp_\ast$.
\begin{lemma}For all linear maps  $\xi$ from $\mathcal{Y}^{(k)}$ to $\mathbb{C}$, and for all endomorphisms $\eta$ of $\mathcal{Y}^{(k)}$, we have
$\xi\ast\eta=\eta\circ (\xi\ast \id)$. Moreover, if $\xi(1)=0$, we have $\exp_\ast(\xi)\ast\id=e^{\xi\ast \id} $.\label{expexp}
\end{lemma}
\begin{proof}
For all $Y_1,Y_2\in \mathcal{Y}^{(k)}$, we have$$m(\xi(Y_1)\otimes \eta(Y))=\xi(Y_1)\eta(Y_2)=\eta\Big(\xi(Y_1)Y_2\Big)=\eta\circ m(\xi(Y_1)\otimes Y_2).$$
Thus, $m\circ(\xi\otimes\eta)=\eta\circ m\circ (\xi\otimes \id)$, and eventually, $$\xi\ast\eta=m\circ(\xi\otimes\eta)\circ \Delta=\eta\circ m\circ (\xi\otimes \id)\circ\Delta=\eta\circ (\xi\ast \id).$$
For the second relation, let us suppose that $\xi(1)=0$. It suffices to apply the previous result with $\eta=\xi\ast \id$. By an immediate induction, for all $l\in \mathbb{N}$, we have
$$\xi^{\ast l}\ast\id=\xi^{\ast (l-1)}\ast(\xi\ast \id)=(\xi\ast \id)\circ (\xi^{\ast (l-1)}\ast\id)=\ldots= (\xi\ast \id)^{(l-1)}\circ (\xi\ast\id)=(\xi\ast \id)^l,$$
and consequently $\exp_\ast(\xi)\ast\id=\displaystyle\sum_{l=0}^\infty\frac{1}{l!}\xi^{\ast l}\ast \id=\displaystyle\sum_{l=0}^\infty\frac{1}{l!}(\xi\ast \id)^l=e^{\xi\ast \id}$.
\end{proof}
Because $\log_{\ast}\chi_\mathbf{A}$ maps $\mathcal{Y}^{(k)}$ to $\mathbb{C}$, using Lemma~\ref{expexp} twice, for all $P\in \mathbb{C}\{X_i:i\in I\}$ we have $$\tau\left(P\left(\mathbf{A}\mathbf{B}\right)\right)=\chi_\mathbf{A}\ast \eta_\mathbf{B}(\rho(P))=\eta_\mathbf{B}\circ (\chi_\mathbf{A}\ast \id)\circ\rho(P)=\eta_\mathbf{B}\circ e^{\log_{\ast}\chi_\mathbf{A}\ast \id}\circ \rho(P).$$
For all $P\in \mathbb{C}\{X_i:i\in I\}$, we also have that $\tau\left((e^{\D_{\mathbf{A}}}P)\left(\mathbf{B}\right)\right)=\eta_\mathbf{B}\circ \rho\circ e^{\D_{\mathbf{A}}}(P)$. Therefore, it remains to prove that $(\log_{\ast}\chi_\mathbf{A}\ast \id)\circ \rho=\rho\circ \D_{\mathbf{A}}$. Let us check this first on monomials. For all $w=(i(1),\ldots,i(n))\in [k]^*$, we have
\begin{eqnarray*}
(\log_{\ast}\chi_\mathbf{A}\ast \id)\circ \rho\left(X_{i(1)}\cdots X_{i(n)}\right)&=&\log_{\ast}\chi_\mathbf{A}\ast \id(Y_w)\\
&=&\displaystyle\sum_{\pi \in NC(n) } \log_{\ast}\chi_\mathbf{A} (Y_{w;\pi})\cdot Y_{w;K(\pi)}.\\
\end{eqnarray*}
In the proof of Proposition~\ref{cumlogcumptwo}, we showed that, for all $\pi\in NC(n)$, $\log_{\ast}\chi_\mathbf{A}(Y_{w;\pi})=L\kappa[\pi](A_{i(1)},\ldots,A_{i(n)})$ if $\pi$ has exactly one block which has more than two elements, and $0$ otherwise. In the case where $\pi$ has exactly one block which has more than two elements, let us denote by $ \{k_1,\cdots,k_m\}$ this block of $\pi$, with $k_1 < k_2 < \ldots < k_m$. We have $K(\pi)=\{\{1,\ldots,k_1-1,k_m,\ldots,n\},\{k_1,\ldots,k_2-1\},\ldots,\{k_{m-1},\ldots,k_m-1\}\}$, which we denote by $\pi(k_1,\ldots, k_m)$. Thus, we have
\begin{eqnarray*}
(\log_{\ast}\chi_\mathbf{A}\ast \id)\circ \rho(X_{i(1)}\cdots X_{i(k)})&=&\displaystyle\sum_{\pi \in NC(n) } \log_{\ast}\chi_\mathbf{A} (Y_{w;\pi})\cdot Y_{w;K(\pi)}\\
&=& \displaystyle\sum_{\substack{2\leq m\\ 1\leq k_1< \ldots < k_m\leq n}}  L\kappa\left(A_{i(k_1)},\ldots,A_{i(k_m)}\right)Y_{w;\pi(k_1,\ldots, k_m)}\\
&=&\rho(\D_{\mathbf{A}}X_{i(1)}\cdots X_{i(k)}).
\end{eqnarray*}
Since $\D_{\mathbf{A}}$ is a derivation, $(\log_{\ast}\chi_\mathbf{A}\ast \id)\circ \rho=\rho\circ \D_{\mathbf{A}}$ will be a direct consequence of the fact that $\log_{\ast}\chi_\mathbf{A}\ast \id$ is a derivation. One can verify this directly, but we remark that $e^{\log_{\ast}\chi_\mathbf{A}\ast \id}=\chi_\mathbf{A}\ast \id$, and because $\chi_\mathbf{A}\ast \id-\id$ is locally nilpotent because $\chi_\mathbf{A}\ast \id-\id$ makes the degree strictly decrease, we also have $\log_{\ast}\chi_\mathbf{A}\ast \id=\log(\chi_\mathbf{A}\ast \id)=-\sum_{l=0}^\infty\frac{1}{l}(\id-\chi_\mathbf{A}\ast \id)^{l}$. Finally, $\chi_\mathbf{A}\ast \id$ is multiplicative, so $\log_{\ast}\chi_\mathbf{A}\ast \id=\log(\chi_\mathbf{A}\ast \id)$ is a derivation. This concludes the proof of Proposition~\ref{freekernelmultweak}.
\end{proof}
Let us finish the proof of Theorem~\ref{freekernelmult}\label{sssmtwo}.
Let $\mathbf{A}=(A_i)_{i\in I}\in \mathcal{A}^I$ be such that $\tau(A_i)\neq 0$ for all $i\in I$ and let $\mathbf{B}=(B_i)_{i\in I}\in \mathcal{A}^I$ be free from $(A_i)_{i\in I}$ such that $\tau(B_i)\neq 0$ for all $i\in I$.

\subsubsection*{Step 1}Proposition~\ref{propcom} says us that $\D_{\tau(\mathbf{A})}$, $\D_{\tau(\mathbf{B})}$ and $\D_{\mathbf{A}/\tau(\mathbf{A})}$ commute, and that $e^{\D_{\mathbf{A}}}=e^{\D_{\mathbf{A}/\tau(\mathbf{A})}}e^{\D_{\tau(\mathbf{A})}}$. Let $(i(1),\ldots,i(n))\in I$, and $P=X_{i(1)}\cdots X_{i(n)}$. Using Proposition~\ref{freekernelmultweak}, we have
\begin{eqnarray*}
 \tau\left(P\left(\mathbf{A}\mathbf{B}\right)\right)&=&
 \tau\left(\tau\left(A_{i(1)}\right)\tau\left(B_{i(1)}\right)\cdots \tau\left(A_{i(n)}\right)\tau\left(B_{i(n)}\right) P\left(\mathbf{A}/\tau(\mathbf{A})\cdot \mathbf{B}/\tau(\mathbf{B})\right)\right)\\
 &=&\tau\left(e^{\D_{\tau(\mathbf{B})}}e^{\D_{\tau(\mathbf{A})}}P\left(\mathbf{A}/\tau(\mathbf{A})\cdot \mathbf{B}/\tau(\mathbf{B})\right)\right)\\
 &=& \tau\left(e^{\D_{\mathbf{A}/\tau(\mathbf{A})}}e^{\D_{\tau(\mathbf{B})}}e^{\D_{\tau(\mathbf{A})}}P\left(\mathbf{B}/\tau(\mathbf{B})\right)\right)\\
 &=&\tau\left(e^{\D_{\tau(\mathbf{B})}}e^{\D_{\mathbf{A}}}P\left(\mathbf{B}/\tau(\mathbf{B})\right)\right).
\end{eqnarray*}
This relation extends to all $\mathbb{C}\{ X_i:i\in I\}$ by induction.

In particular, and because $e^{\D_{\mathbf{A}}}P\in \mathbb{C}\{ X_i:i\in I\}$, we have
\begin{equation}
\tau\left(e^{\D_{\mathbf{A}}}P\left(\mathbf{B}\right)\right)= \tau\left(e^{\D_{\mathbf{A}}}P\left(\tau(\mathbf{B})\mathbf{B}/\tau(\mathbf{B})\right)\right) = \tau\left(e^{\D_{\tau(\mathbf{B})}}e^{\D_{\mathbf{A}}}P\left(\mathbf{B}/\tau(\mathbf{B})\right)\right)= \tau\left(P\left(\mathbf{A}\mathbf{B}\right)\right).\label{eqfreekernelmultweak}
\end{equation}

\subsubsection*{Step 2}

We will use the following characterization of conditional expectation. The element $\tau\left(\left.P(\mathbf{A}\mathbf{B})\right|\mathbf{B}\right)$ is the unique element of $W^*(\mathbf{B})$ such that, for all $B_{{i_0}} \in W^*(\mathcal{B})$, $$\tau\left(P(\mathbf{A}\mathbf{B})B_{{i_0}}\right) = \tau\left(\tau\left(\left.P(\mathbf{A}\mathbf{B})\right|\mathbf{B}\right)B_{{i_0}}\right),$$
and since $e^{\D_{\mathbf{A}}}P\left(\mathbf{B}\right)\in W^*(\mathbf{B})$, it remains to prove that, for all $B_{{i_0}} \in W^*(\mathcal{B})$,
$$\tau\left(P(\mathbf{A}\mathbf{B})B_{{i_0}}\right)= \tau\left((e^{\D_{\mathbf{A}}}P)\left(\mathbf{B}\right)B_{{i_0}}\right).$$

In order to use Lemma~\ref{freekernelmultweak}, we work on $\mathbb{C}\{ X_i:i\in I\cup\{{i_0}\}\}$. Let $R_{{i_0}}:P\mapsto PX_{i_0}$ be the operator of right multiplication by $X_{i_0}$ on $\mathbb{C}\{ X_i:i\in I\cup\{{i_0}\}\}$.
Let $A_{i_0}=1$ and $B_{{i_0}} \in W^*(\mathcal{B})$.
On one hand, we have $P(\mathbf{A}\mathbf{B})B_{{i_0}}=(R_{{i_0}}P)(\mathbf{A}\mathbf{B},A_{i_0} B_{{i_0}})$, and using \eqref{eqfreekernelmultweak}, we have $\tau\left(P(\mathbf{A} \mathbf{B})B_{{i_0}}\right)=\tau\left(\left(e^{\D_{\mathbf{A},A_{i_0}}}R_{{i_0}}P\right)\left(\mathbf{B},B_{{i_0}}\right)\right).$ On the other hand, $\tau\left((e^{\D_{\mathbf{A}}}P)\left(\mathbf{B}\right)B_{{i_0}}\right)=\tau\left(\left(R_{{i_0}}e^{\D_{\mathbf{A},A_{i_0}}}(P)\right)\left(\mathbf{B},B_{{i_0}}\right)\right)$.

Thus, it suffices to prove that the operators $\D_{\mathbf{A},A_{i_0}}$ and $R_{{i_0}}$ commute, which is essentially the same verification as the end of the proof of Theorem~\ref{freekernel} at Section~\ref{freekernelproof}.
\end{proof}

\subsubsection{Example} Let $t\geq 0$. Let $U_t$ be a free unitary Brownian motion at time $t$. Let us compute\label{utex} $\D_{\tau(U_t)} X^2=-t X^2,$ 
$\D_{U_t/\tau(U_t)} X^2=-t X \tr X,$ and $\left(\D_{U_t/\tau(U_t)}\right)^2 X^2=0.$ Thus,
$$
e^{\D_{U_t}}X^2=e^{\D_{U_t/\tau(U_t)} }e^{\D_{\tau(U_t)}}X^2\\
=e^{\D_{U_t/\tau(U_t)} }\left( e^{-t}X^2\right)\\
= e^{-2t}\left(X^2-t X \tr X+0\right).
$$
Using Theorem~\ref{freekernelmult}, we have for all $B\in \mathcal{A}$ free from $U_t$, 
$\tau\left(\left(U_tB\right)^2|B\right)=e^{-2t}\left(B^2-t B \tau (B)\right)$.
Let us derive a more general fact which will be useful in the proof of Theorem~\ref{fht}. For all $n\in \mathbb{N}$, set $W_n=\text{span} \{Q\tr R: Q\in \mathbb{C}_{n-1}[X], R\in \mathbb{C}\{X\} \}\subset \mathbb{C}\{X\}$. For all $n\in \mathbb{N}^*$, we have
$$
e^{\D_{U_t}}X^n=e^{\D_{U_t/\tau(U_t)} }e^{\D_{\tau(U_t)}}X^n\\
=e^{\D_{U_t/\tau(U_t)} }\left( e^{-nt/2}X^n\right)\\
=e^{-nt/2}X^n+P
$$
with $P\in W_{n-1}$. Indeed, $\D_{U_t/\tau(U_t)}$ maps the space $W_{n}$ into the space $W_{n-1}$, so the unique term which is not in $W_{n-1}$ in the sum $e^{-nt/2}\sum_{k=0}^\infty \frac{1}{k!}(\D_{U_t/\tau(U_t)})^kX^n$ is $e^{-nt/2}X^n$.

\section{Free Segal-Bargmann and free Hall transform}\label{fsb}

This section is devoted to provide a new construction of the free Segal-Bargmann transform and of the free Hall transform defined in~\cite{Biane1997} using free convolution or, more precisely, using the free transition operators of Theorems~\ref{freekernel} and~\ref{freekernelmult}. The free Segal-Bargmann transform is completely characterized in Theorem~\ref{fsbt} and the free Hall transform in Theorem~\ref{fht}.

For convenience, fix a (sufficiently large) $W^*$-probability space $\left(\mathcal{A},\tau\right)$ to use throughout this section. We start with a very succinct review of some results of~\cite{Nelson1974}. Let us suppose that $\left(\mathcal{A},\tau\right)$ is represented on a Hilbert space $\mathcal{H}$, in the sense that $\mathcal{A}$ is a subspace of bounded operators on $\mathcal{H}$. The operator norm of $\mathcal{A}$ is denoted by $\|\cdot \|_{\infty}$. A (not necessarily bounded or everywhere defined) operator $A$ on $\mathcal{H}$ is said to be affiliated with $\mathcal{A}$ if $AU = UA$ for every unitary operator $U$ in the commutant of $\mathcal{A}$.
The set of closed densely-defined operators affiliated with $\mathcal{A}$ is denoted by $\mathfrak{M}(\mathcal{A})$. For all $A\in \mathcal{A}$, let us denote by $\|A\|_2$ the norm $\|A\|_2^2=\tau \left(A^* A\right)$. If $A\in \mathfrak{M}(\mathcal{A})$, we can still define the norm $\|\cdot\|_2$ of $A$ (not necessarily finite), by extending the trace $\tau$ to general positive operators. The space $L^2(\mathcal{A},\tau)=\{A \in \mathfrak{M}(\mathcal{A}):\|A\|_2<\infty \}$ is a Hilbert space for the norm $\|\cdot\|_2$ in which $\mathcal{A}$ is dense. When we consider a Hilbert space completion for the norm $\|\cdot\|_2$, it will always be identified with a subset of $L^2(\mathcal{A},\tau)$, and so with a subset of $\mathfrak{M}(\mathcal{A})$.

\subsection{Semi-circular system and circular system}
 Let $H$ be a real Hilbert space, with inner product $\langle \cdot,\cdot\rangle_H$.
 
\subsubsection{Semi-circular system}A linear map $\mathbf{s}=(s(h))_{h\in H}$ from $H$ to $\mathcal{A}$ is called a semi-circular system if
\begin{enumerate}
\item for each $h\in H$, the element $s(h)$ is a semi-circular random variable of variance $\|h\|_H^2$,
\item for each orthogonal family $h_1,\ldots, h_n$ in $H$, the elements $s(h_1),\ldots,s( h_n)$ are free.
\end{enumerate}
Let $\mathbf{s}$ be a semi-circular system. We denote by $L^2(\mathbf{s},\tau)$ the Hilbert completion of the algebra generated by $s(H)$ for the norm $\|\cdot\|_2:A\mapsto \tau(A^* A)^{1/2}$. Let us compute the free cumulants of $\mathbf{s}$. Let $n\in \mathbb{N}$ and $h_1,\ldots, h_n\in H$. By Proposition~\ref{freeness}, by the definition of a semi-circular random variable in Section~\ref{sc}, and by the linearity of the free cumulants, we have $$\kappa(s(h_1),s(h_2))=\kappa(s(h_1)^*,s(h_2)^*)=\kappa(s(h_1)^*,s(h_2))=\kappa(s(h_1),s(h_2)^*)=\left\langle h_1,h_2\right\rangle_H,$$
and $\kappa(s(h_1)^{\varepsilon_1},\ldots,s(h_n)^{\varepsilon_n})=0$ for any $n\neq 2$ and $\varepsilon_1,\ldots, \varepsilon_n\in \{1,*\}$.\label{scs}
 
\subsubsection{Circular system}A linear map $\mathbf{c}=(c(h))_{h\in H}$ from $H$ to $\mathcal{A}$ is called a circular system if
\begin{enumerate}
\item the maps $\sqrt{2}\Re (c):h\mapsto \frac{1}{\sqrt{2}}(c(h)+c(h)^*)$ and $\sqrt{2}\Im (c):h\mapsto \frac{1}{\sqrt{2}i}(c(h)-c(h)^*)$ are semi-circular systems,
\item $\Re (c)(H)$ and $\Im (c)(H)$ are free.
\end{enumerate}
Let $\mathbf{c}$ be a circular system. We denote by $L^2_{\hol}(\mathbf{c},\tau)$ the Hilbert completion of the algebra generated by $c(H)$ for the norm $\|\cdot\|_2:A\mapsto \tau(A^* A)^{1/2}$.\label{cs} Let us compute the free cumulants of $\mathbf{c}$. Let $n\in \mathbb{N}$ and $h_1,\ldots, h_n\in H$. By the linearity of the free cumulants, we have \begin{eqnarray*}\kappa(c(h_1),c(h_2))&=&\kappa\left(\left(\Re (c)+i\Im (c)\right)(h_1),\left(\Re (c)+i\Im (c)\right)(h_2)\right)\\
&=&\frac{1}{2}\kappa\Big(\sqrt{2}\Re (c)(h_1),\sqrt{2}\Re (c)(h_2)\Big)-\frac{1}{2}\kappa\Big(\sqrt{2}\Im (c)(h_1),\sqrt{2}\Im (c)(h_2)\Big)\\
&=&0,
\end{eqnarray*}
and similarly, $\kappa(c(h_1)^*,c(h_2))=\kappa(c(h_1),c(h_2)^*)=\langle h_1,h_2\rangle_H$, $\kappa(c(h_1)^*,c(h_2)^*)=0$, and\\$\kappa(c(h_1)^{\varepsilon_1},\ldots,c(h_n)^{\varepsilon_n})=0$ for any $n\neq 2$ and $\varepsilon_1,\ldots, \varepsilon_n\in \{1,*\}$.

\subsection{Free Segal-Bargmann transform}

Let $\mathbf{s}$ be a semi-circular system, and $\mathbf{c}$ be a circular system from $H$ to $\mathcal{A}$. The free Segal-Bargmann transform is an isomorphism from $L^2\left(\mathbf{s},\tau\right)$ to $L^2_{\hol}\left(\mathbf{c},\tau\right)$ defined by Biane in~\cite{Biane1997}.

\subsubsection{Definition}
Let us define the Tchebycheff type II polynomials $\left(T_n\right)_{n\in \mathbb{N}}$ by their generating function: for all $|z|<1$ and $-2<x<2$, we have
$$\displaystyle\sum_{n=0}^\infty z^n T_n(x)=\frac{1}{1-xz+z^2}.$$
We remark that, for all $n\in \mathbb{N}$, the degree of $T_n$ is $n$ and the leading coefficient is the coefficient of $x^nz^n$ in the development and is therefore equal to $1$.

A quick way to define the free Segal-Bargmann transform is to define it on polynomials. Let $\left(h_j\right)_{j\in J}$ be an orthonormal basis of $H$. Let us define $\mathcal{G}$ by the unique linear operator on $\mathbb{C}\langle X_{h_j}:j\in J\rangle$ such that, for all $n\in \mathbb{N}$, $k_1,\ldots, k_n \in \mathbb{N}^*$ and $j(1)\neq j(2)\neq \cdots \neq j(n-1) \neq j(n)$ elements of $J$, one has $\mathcal{G}\Big(T_{k_1}(X_{h_{j(1)}})\cdots T_{k_n}(X_{h_{j(n)}})\Big)=X_{h_{j(1)}}\cdots X_{h_{j(n)}}$. The following theorem due to Biane corresponds to Definition~4 of~\cite{Biane1997}.
\begin{theo*}[Biane~\cite{Biane1997}]
The map $P \big(\mathbf{s} \big) \mapsto \mathcal{G}(P) \big(\mathbf{c} \big)  $ for all $P\in \mathbb{C}\langle X_{h_j}:j\in J\rangle$ is an isometric map which extends to a Hilbert isomorphism between $L^2\left(\mathbf{s},\tau\right)$ and $L^2_{\hol}\left(\mathbf{c},\tau\right)$ called the free Segal-Bargmann transform.
\end{theo*}

\subsubsection{Another construction}Theorem~\ref{fsbt} links the free Segal-Bargmann transform and the free convolution. The link with the free convolution is already present in~\cite{Biane1997} in the $1$-dimensional case.
\begin{theorem}
Let $\mathbf{s}$ be a semi-circular system, and $\mathbf{c}$ be a circular system from $H$ to $\mathcal{A}$.\label{fsbt}
For all $P\in \mathbb{C}\{X_h:h\in H\}$, $\mathcal{F}:P \left(\mathbf{s} \right) \mapsto (e^{\Delta_{\mathbf{s}}}P) \left(\mathbf{c} \right)  $ is an isometric map which extends to a Hilbert space isomorphism $\mathcal{F}$ between $L^2\left(\mathbf{s},\tau\right)$ and $L^2_{\hol}\left(\mathbf{c},\tau\right)$. Moreover, this isomorphism is the free Segal-Bargmann transform.

In particular, if $\mathbf{s}$ and $\mathbf{c}$ are free, for all $P\in \mathbb{C}\{X_h:h\in H\}$,
$$\mathcal{F}\Big(P\left(\mathbf{s}\right)\Big) =\tau\Big(P\left(\mathbf{s}+\mathbf{c}\right)\Big|\mathbf{c}\Big).$$
\end{theorem}
It should be remarked that, for the map $\mathcal{F}$ to be well-defined, it must be true that, for all $P, Q\in \mathbb{C}\{X_h:h\in H\}$, if $P \left(\mathbf{s} \right)=Q \left(\mathbf{s} \right)$, then $(e^{\Delta_{\mathbf{s}}}P) \left(\mathbf{c} \right)=(e^{\Delta_{\mathbf{s}}}Q) \left(\mathbf{c} \right)$. This fact is contained in the following proof.
\begin{proof}
We will prove in a first step that, for all $P\in \mathbb{C}\{X_h:h\in H\}$, we have
\begin{equation}\left\|P \left(\mathbf{s} \right)\right\|^2_{L^2\left(\mathbf{s},\tau\right)}=\left\| (e^{\Delta_{\mathbf{s}}}P) \left(\mathbf{c} \right)\right\|^2_{L^2_{\hol}\left(\mathbf{c},\tau\right)}. \label{auxi}\end{equation}
This proves that, for all $P\in \mathbb{C}\{X_h:h\in H\}$, if $P \left(\mathbf{s} \right)=Q \left(\mathbf{s} \right)$, then $$\left\| (e^{\Delta_{\mathbf{s}}}P) \left(\mathbf{c} \right)-(e^{\Delta_{\mathbf{s}}}Q) \left(\mathbf{c} \right)\right\|^2_{L^2_{\hol}\left(\mathbf{c},\tau\right)}=\left\| (e^{\Delta_{\mathbf{s}}}(P-Q)) \left(\mathbf{c} \right)\right\|^2_{L^2_{\hol}\left(\mathbf{c},\tau\right)}=\left\|(P-Q) \left(\mathbf{s} \right)\right\|^2_{L^2\left(\mathbf{s},\tau\right)}=0,$$and consequently $(e^{\Delta_{\mathbf{s}}}P) \left(\mathbf{c} \right)=(e^{\Delta_{\mathbf{s}}}Q) \left(\mathbf{c} \right)$. Thus, $\mathcal{F}$ is a well-defined isometric map, and so it extends to a Hilbert isomorphism $\mathcal{F}$ between $L^2\left(\mathbf{s},\tau\right)$ and $\mathcal{F}\left(L^2\left(\mathbf{s},\tau\right)\right)\subset L^2_{\hol}\left(\mathbf{c},\tau\right)$. The surjectivity is clear since, for all $P\in \mathbb{C}\left\langle X_h:h\in H\right\rangle$,
$$P \left(\mathbf{c}\right)=(e^{\Delta_{\mathbf{s}}}(e^{-\Delta_{\mathbf{s}}}P)) \left(\mathbf{c}\right) =\mathcal{F}\Big(( e^{-\Delta_{\mathbf{s}}}P) \left(\mathbf{c}\right)\Big).$$

The equality with the free Segal-Bargmann transform will be made in a second step, and the last part of Theorem~\ref{fsbt} follows from Theorem~\ref{freekernel}, which says in this case that, for all $P\in \mathbb{C}\{X_h:h\in H\}$, $(e^{\Delta_{\mathbf{s}}}P) \left(\mathbf{c} \right)=\tau(P(\mathbf{s}+\mathbf{c})|\mathbf{c}).$

\subsubsection*{Step $1$}In the aim to prove \eqref{auxi}, we work on the *-algebra $\mathbb{C}\{X_h:h\in H\cup (H\times \{*\})\}=\mathbb{C}\{X_h,X_h^*:h\in H\}$ (see Section~\ref{involution}).

Let $\mathbf{s}^{+}$ and $\mathbf{s}^{-}$ be two semi-circular system from $H$ to $\mathcal{A}$ such that $\mathbf{s}^{+}$, $\mathbf{s}^{-}$ and $\mathbf{c}$ are free.
Let us define five auxiliary maps from $H\cup (H\times \{*\})$ to $\mathcal{A}$.
The first four are the natural extensions of $\mathbf{s}$, $\mathbf{c}$, $\mathbf{s}^{+}$ and $\mathbf{s}^{-}$. For all $h\in H$, we set $\mathbf{s}(h)=s(h)=\mathbf{s}((h,*))=s(h)^*$, $\mathbf{c}(h)=c(h)$, $\mathbf{c}((h,*))=c(h)^*$, $\mathbf{s}^{+}(h)=s^{+}(h)=\mathbf{s}^{+}((h,*))=s^{+}(h)^*$ and $\mathbf{s}^{-}(h)=s^{-}(h)=\mathbf{s}^{-}((h,*))=s^{-}(h)^*$.
The last map $\tilde{\mathbf{s}}$ is defined for all $h\in H$ by $\tilde{\mathbf{s}}(h)=c(h)+s^{+}(h)$, and $\tilde{\mathbf{s}}((h,*))=c(h)^*+s^{-}(h)$.
We remark that $\mathbf{s}$, $\mathbf{s}^{+}$ and $\mathbf{s}^{-}$ have the same distribution, and that $\tilde{\mathbf{s}}$ is equal to $\mathbf{c}+\mathbf{s}^{+}$ on $H$ and is equal to $\mathbf{c}+\mathbf{s}^{-}$ on $H\times \{*\}$.

Let $P\in \mathbb{C}\{X_h:h\in H\}$. Theorem~\ref{freekernelmult} gives us that
$$
\left\|P \left(\mathbf{s} \right)\right\|^2_{L^2\left(\mathbf{s},\tau\right)}=\tau\Big( P(\mathbf{s})^*P (\mathbf{s})\Big)\\
=\tau\Big(P^*P(\mathbf{s})\Big)\\
=\left(e^{\Delta_{\mathbf{s}}}\left(P^*P\right)\right)(0)
 $$
and
$$\begin{array}{rcl}
\left\| e^{\Delta_{\mathbf{s}}}P \left(\mathbf{c} \right)\right\|^2_{L^2_{\hol}\left(\mathbf{c},\tau\right)}&=&\tau\left( \Big(e^{\Delta_{\mathbf{s}}}P \left(\mathbf{c} \right)\Big)^* \Big(e^{\Delta_{\mathbf{s}}}P \left(\mathbf{c} \right)\Big)\right)\\
&=&\tau\left(\tau \Big(P\left(\mathbf{c}+\mathbf{s}^{-}\right)|\mathbf{c},\mathbf{s}^{+}\Big)^* \tau \Big(P\left(\mathbf{c}+\mathbf{s}^{+}\right)|\mathbf{c},\mathbf{s}^{-}\Big)\right)\\
&=&\tau\Big( P\left(\mathbf{c}+\mathbf{s}^{-}\right)^* P\left(\mathbf{c}+\mathbf{s}^{+}\right)\Big)\\
&=&\tau\Big( P^*\left(\mathbf{c}+\mathbf{s}^{-}\right) P\left(\mathbf{c}+\mathbf{s}^{+}\right)\Big)\\
&=&\tau\Big( P^*\left(\tilde{\mathbf{s}}\right) P\left(\tilde{\mathbf{s}}\right)\Big)\\
&=&\tau\Big( P^*P\left(\tilde{\mathbf{s}}\right) \Big)\\
&=&\left(e^{\Delta_{\tilde{\mathbf{s}}}}\left(P^*P\right)\right)(0)
\end{array}$$
It remains to prove that $\Delta_{\mathbf{s}}=\Delta_{\tilde{\mathbf{s}}}$ on $\mathbb{C}\{X_h,X_h^*:h\in H\}$. By definition, it suffices to prove that the free cumulants of $(\mathbf{s}(h))_{H\cup (H\times \{*\})}$ and $(\tilde{\mathbf{s}}(h))_{H\cup (H\times \{*\})}$ coincide.

Recall that $0_{\mathcal{A}}$ is free from all element of $\mathcal{A}$, so all free cumulants involving $0_{\mathcal{A}}$ are equal to $0$.
The data of the free cumulants of a semi-circular system (see Section~\ref{scs}), of a circular system (see Section~\ref{cs}) and the linearity of the free cumulants allows us to deduce the free cumulants of $(\mathbf{s}(h))_{H\cup (H\times \{*\})}$ and $(\tilde{\mathbf{s}}(h))_{H\cup (H\times \{*\})}$.

We observe first that all free cumulants of order different from $2$ involved are equal to $0$. Thus, for all $n\in \mathbb{N}$ such that $n\neq 2$ and $h_1,\ldots, h_n\in H\cup (H\times \{*\})$, we have$$\kappa\left(\mathbf{s}(h_1),\ldots,\mathbf{s}(h_n)\right)=0=\kappa\left(\tilde{\mathbf{s}}(h_1),\ldots,\tilde{\mathbf{s}}(h_n)\right).$$
Only cumulants of order $2$ are non trivial.
Let $h_1, h_2 \in H$. We have
$$\begin{array}{ll}
\kappa\Big(\mathbf{s}(h_1),\mathbf{s}(h_2)\Big)=\kappa\left(\mathbf{s}(h_1),\mathbf{s}((h_2,*))\right)=\kappa\left(\mathbf{s}((h_1,*)),\mathbf{s}((h_2,*))\right)=\kappa(s(h_1),s(h_2))&=\left\langle h_1,h_2\right\rangle_H,\\
\kappa\Big(\tilde{\mathbf{s}}(h_1),\tilde{\mathbf{s}}(h_2)\Big)=\kappa\Big(c(h_1)+s^{+}(h_1),c(h_2)+s^{+}(h_2)\Big)=\kappa\Big(s^{+}(h_1),s^{+}(h_2)\Big) &=\left\langle h_1,h_2\right\rangle_H,\\
\kappa\Big(\tilde{\mathbf{s}}(h_1),\tilde{\mathbf{s}}((h_2,*))\Big)=\kappa\Big(c(h_1)+s^{+}(h_1),c(h_2)^*+s^{-}(h_2)\Big)=\kappa\Big(c(h_1),c(h_2)^*\Big) &= \left\langle h_1,h_2\right\rangle_H,
\end{array}$$
and
$$\kappa\Big(\tilde{\mathbf{s}}((h_1,*)),\tilde{\mathbf{s}}((h_2,*))\Big)=\kappa\Big(c(h_1)^*+s^{-}(h_1),c(h_2)^*+s^{-}(h_2)\Big)=\kappa\Big(s^{-}(h_1),s^{-}(h_2)\Big) = \left\langle h_1,h_2\right\rangle_H.$$
At the end, for all $h_1, h_2 \in H\cup (H\times \{*\})$, we have $\kappa\Big(\mathbf{s}(h_1),\mathbf{s}(h_2)\Big)=\kappa\Big(\tilde{\mathbf{s}}(h_1),\tilde{\mathbf{s}}(h_2)\Big)$.

\subsubsection*{Step $2$}We prove now that the isomorphism $\mathcal{F}$ is indeed the free Segal-Bargmann transform. That can be done using the factorization of Section~\ref{fact}. More precisely, with the aim of working on polynomials, we will use the identity $P\left(\mathbf{c}\right)=\left.P\right|_{\mathbf{c}}\left(\mathbf{s}\right)$ for all $P\in \mathbb{C}\{ X_h:h\in H\}$.

For all $P\in \mathbb{C}\langle X_{h_j}:j\in J\rangle$, we have $(e^{\Delta_{\mathbf{s}}}P)|_{\mathbf{c}}=\mathcal{G}(P)$. This fact is easily obtained by induction on the degree of $P$. Indeed, $P\mapsto (e^{\Delta_{\mathbf{s}}}P)|_{\mathbf{c}}$ and $P\mapsto \mathcal{G}(P)$ both respect the degree of $P$, preserve leading coefficients and are isometries between $(\mathbb{C}\langle X_{h_j}:j\in J\rangle,\|\cdot\|_{\mathbf{s}})$ and $(\mathbb{C}\langle X_{h_j}:j\in J\rangle,\|\cdot\|_{\mathbf{c}})$ for the norms $\|\cdot\|_{\mathbf{s}}:P\mapsto \|P(\mathbf{s})\|_2$ and $\|\cdot\|_{\mathbf{c}}:P\mapsto \|P(\mathbf{c})\|_2$. 
This proves that, for all $P\in \mathbb{C}\langle X_{h_j}:j\in J\rangle$, we have $\mathcal{F}(P(\mathbf{s})) =e^{\Delta_{\mathbf{s}}}P\left(\mathbf{c}\right)=(e^{\Delta_{\mathbf{s}}}P)|_{\mathbf{c}}\left(\mathbf{c}\right)=\mathcal{G}(P)\left(\mathbf{c}\right)$, establishing the equality between $\mathcal{F}$ and the free Segal-Bargmann transform on the algebra generated by $(s(h_j))_{j\in J}$, and because they are isomorphisms, this equality extends to $L^2\left(\mathbf{s},\tau\right)$.
\end{proof}

\subsection{Free stochastic calculus}\label{fsc}
\subsubsection{Semi-circular and circular Brownian motion}
A free standard Brownian motion, or semi-circular Brownian motion, is a family $\left(X_t\right)_{t\geq 0}$ of self-adjoint elements in $\mathcal{A}$, such that
\begin{enumerate}
  \item $X_0=0$;
  \item  For all $0\leq s<t $, the element $X_t-X_s$ is semi-circular of variance $t-s$;
  \item For all $0\leq t_1<\ldots<t_n$, the elements $X_{t_1},X_{t_2-t_1},\ldots,X_{t_n-t_{n-1}}$ form a free family.
\end{enumerate}
Let $\mathcal{B}$ be a von Neumann subalgebra of $\mathcal{A}$, free from $W^*(X_t, t\geq 0)$. We denote by $\mathcal{B}_t$ the von Neumann algebra generated by $W^*(X_s,0\leq s\leq t)\cup \mathcal{B}$.

A free circular standard Brownian motion is a family $(Z_t)_{t\geq 0}$ of non-commutative random variables in $\mathcal{A}$ such that $\frac{1}{\sqrt{2}}(Z_t+Z_t^*)_{t\geq 0}$ and $(\frac{1}{\sqrt{2}i}(Z_t-Z_t^*))_{t\geq 0}$ are two free standard Brownian motions which are free from each other.

\subsubsection{Free stochastic integration}Let us recall briefly some basic definitions of free stochastic calculus. For simplicity, let us consider a very restricted situation. For further information, and more developments, see~\cite{Biane1997},~\cite{Biane1998a} or~\cite{Kargin2011}. We will treat in parallel the semi-circular and the circular case.

Let $t\geq 0$. Let $s\mapsto A_s$ and $s\mapsto B_s$ be maps from $[0,t]$ to $\mathcal{A}$, continuous and uniformly bounded for the operator norm $\| \cdot\|_{\infty}$. We suppose that, for all $s\in[0,t]$, $A_s$ and $B_s\in \mathcal{B}_s$ (respectively $W^*(Z_u,u\leq s)$). Such mappings are called bounded adapted semi-circular (respectively circular) processes on $[0,t]$.
Then we can define the free stochastic integral $\int_0^t A_s \diff X_s B_s \in \mathcal{B}_t$ with respect to $X$ (respectively the free stochastic integrals $\int_0^t A_s \diff Z_s B_s \in W^*(Z_s,s\leq t)$ and $\int_0^t A_s \diff Z_s^* B_s \in W^*(Z_s,s\leq t)$ with respect to $Z$ and to $Z^*$). A bounded adapted process on $\mathbb{R}_+=[0,\infty)$ is a bounded adapted process on all $[0,t]\subset\mathbb{R}_+$.

The following lemma includes properties corresponding to Lemma 10 and Proposition 7 of~\cite{Biane1997}, and to Proposition 3.2.3 of~\cite{Biane1998a}.
\begin{lemma}
\label{fmart}Let $A$ and $B$ be two bounded semi-circular (or respectively circular) adapted processes. We have for all $t\geq 0$, $\tau \left(\left.\int_0^t A_s \diff X_s B_s \right|\mathcal{B} \right)=0$ and
$ \left(\int_0^t A_s \diff X_s B_s \right)^*=\int_0^t B_s^* \diff X_s A_s^*$
in the semi-circular case.

We have for all $t\geq 0$, $\tau \left(\int_0^t A_s \diff Z_s B_s \right)=0$, $\tau \left(\int_0^t A_s \diff Z_s^* B_s \right)=0$ and
$\left(\int_0^t A_s \diff Z_s B_s \right)^*=\int_0^t B_s^* \diff Z_s^* A_s^*$
in the circular case.
\end{lemma}

We have It\^o formulas for free stochastic integrals (Proposition 8 of~\cite{Biane1997}, or Theorem 4.1.2 of~\cite{Biane1998a}) which are summed up below, using formal rules. Let $A$, $B$, $C$ and $D$ be bounded adapted semi-circular processes. We have 
$$A_t \diff t \cdot C_t \diff t = A_t \diff t \cdot C_t \diff X_t D_t= A_t\diff X_t B_t \cdot C_t \diff t = 0,$$
and
$$A_t\diff X_t B_t\cdot C_t \diff X_t D_t= \tau\left(B_t C_t\right)A_t D_t \diff t.$$
Let $A$, $B$, $C$ and $D$ be bounded adapted circular processes. We have
$$A_t \diff t \cdot C_t \diff t = A_t \diff t \cdot C_t \diff Z_t D_t= A_t\diff Z_t B_t \cdot C_t \diff t = A_t \diff t \cdot C_t \diff Z_t^* D_t= A_t\diff Z_t^* B_t \cdot C_t \diff t=0,$$
$$A_t \diff Z_t B_t  \cdot C_t \diff Z_t D_t =A_t \diff Z_t^* B_t  \cdot C_t \diff Z_t^* D_t=0$$
and
$$A_t \diff Z_t^* B_t \cdot  C_t \diff Z_t D_t =A_t \diff Z_t B_t \cdot  C_t \diff Z_t^* D_t=\tau\left(B_t C_t\right)A_t D_t \diff t .$$
As an example, let us write the It\^o formula in terms of free stochastic integrals in the semi-circular case. For all $t\geq 0$\label{Ito}, we have
\begin{eqnarray*}
\int_0^t A_s \diff X_s B_s\int_0^t C_s \diff X_s D_s &=&\int_0^t A_s \diff X_s \left[B_s\int_0^s C_u \diff X_u D_u\right]+\int_0^t \left[\int_0^s A_u \diff X_u B_uC_s\right] \diff X_s D_s \\
 &  & \hspace{2cm}+\int_0^t  \tau\left(B_s C_s\right) A_sD_s \diff s.
\end{eqnarray*}

\subsection{Free unitary Brownian motion}

The (right) free unitary Brownian motion $(U_t)_{t\geq 0}$ is defined to be the unique bounded adapted semi-circular process which is the solution of the following free stochastic differential equation\label{uFBMsection}
\begin{equation}
\left\{
\begin{array}{c c c}
    U_0 &=& \Id, \\
    \diff U_t &=& i \diff X_t U_t-\frac{1}{2}U_t \diff t. \label{uFBM} \\
\end{array}
\right.\end{equation}
It can be constructed using Picard iteration (see for example~\cite{Kargin2011}). From Lemma~\ref{fmart}, we know that $(U_t^*)_{t\geq 0}$ is a bounded adapted semi-circular process defined by the free stochastic differential equation
\begin{equation}
\left\{
\begin{array}{c c c}
    U_0^* &=& \Id, \\
    \diff U_t^* &=& -i U_t^* \diff X_t -\frac{1}{2}U_t^* \diff t. \label{uAFBM} \\
\end{array}
\right.\end{equation}
Thanks to the free It\^o formula, one obtains $U_0U _0^*=\Id$ and for all $t\geq 0$,
$$ \diff \left(U_t U_t^*\right)
= i  U_t \diff X_tU_t^*- i  U_t \diff X_t U_t^*+U_tU_t^*\diff t-\frac{1}{2}U_tU_t^*\diff t-\frac{1}{2}U_tU_t^*\diff t=0.
$$
Thus, for all $t\geq 0$, $U_t^*=U_t^{-1}$.
The distribution of $(U_t)_{t\geq 0}$ was computed in~\cite{Biane1997a}, and has been recalled in Section~\ref{fu}.

\subsubsection{Extension of $\D_{U_t}$}\label{deltau}The purpose of this section is to establish Proposition~\ref{distU} below, which is an extension of Theorem~\ref{freekernelmult} for the unitary Brownian motion. We begin by extending $\D_{U_t}$ on $\mathbb{C}\{X,X^*,X^{-1},{X^*}^{-1}\}$.

Let $k\in \mathbb{N}$, and $X_1,\ldots,X_k\in \{X,X^*,X^{-1},{X^*}^{-1}\}$. For all $1\leq i\leq k$, set $l (i)=0$, $r(i)=1$ and $\epsilon(i)=1$ if $X_i=X$ or $X_i={X^*}^{-1}$, and $l(i)=1$, $r(i)=0$ and $\epsilon(i)=-1$ if $X_i=X^*$ or $X_i=X^{-1}$. Informally, the numbers $l(i)$ and $r(i)$ indicate if $X_i$ is a left multiplicative process, or a right one, and the sign $\epsilon(i)$ reflects the coefficient in the stochastic equation verified by $X_i$.

For all $1\leq i<j\leq k$, the term $\overbracket{X_i\cdots X_j}$ refers to the product $X^{r(i)}_iX_{i+1}\cdots X_{j-1}X_j^{l(j)}$, that is to say the product $X_i\cdots X_j$  possibly excluding $X_i$ and/or $X_j$. We set
$$\Delta_U X_1\cdots X_k=- k X_1\cdots X_k-2\displaystyle\sum_{1\leq i< j\leq k} \epsilon(i)\epsilon(j)X_1\cdots  \widehat{\overbracket{X_i\cdots X_j}}\cdots X_k \tr \left(\overbracket{X_i\cdots X_j}\right),$$
where the hat means that we have omitted the term $\overbracket{X_i\cdots X_j} $ in the product $X_1\cdots  X_k$.
We extend $\Delta_U$ to all $\mathbb{C}\{X,X^*,X^{-1},{X^*}^{-1}\}$ by linearity and by the relation of derivation
\begin{equation}
\forall P,Q \in \mathbb{C}\{X,X^*,X^{-1},{X^*}^{-1}\}, \Delta_U \left(P\tr Q\right)=\left(\Delta_U P\right)\tr Q+P\tr \left(\Delta_U Q\right).\label{structDeltaU}
\end{equation}
Let $t\geq 0$. We observe that, if we consider $\frac{t}{2} \Delta_U$ on $\mathbb{C}\{X\}$, the term $\overbracket{X_i\cdots X_j}$ is always the product $X_iX_{i+1}\cdots X_{j-1}$, and consequently $\frac{t}{2} \Delta_U$ coincides with $\D_{U_t}$ on $\mathbb{C}\{X\}$, as defined in Section~\ref{dut} using the free log-cumulants of $U_t$ which are stated in Section~\ref{logcumut}.\label{DeltaD}

\begin{lemma}
Let $B$ be an invertible variable of $\mathcal{B}$. For all $P\in\mathbb{C}\{X,X^*,X^{-1},{X^*}^{-1}\}$,\label{difftau}
$$\frac{\diff}{\diff t} \tau\left(P\left(U_tB\right)|\mathcal{B}\right)=\tau\left(\left.\frac{1}{2}\Delta_U P\left(U_t B\right)\right|\mathcal{B}\right) .$$
\end{lemma}
\begin{proof}
Let us define $(U_t^B=U_tB)_{t\geq 0}$. We deduce from \eqref{uFBM} that $(U_t^B)_{t\geq 0}$ is a bounded adapted semi-circular process defined by the free stochastic differential equation
\begin{equation}
\left\{
\begin{array}{r c l}
    U_0^B &=& B, \\
    \diff U_t^B &=& i \diff X_t U_t^B-\frac{1}{2}U_t^B \diff t. \label{uBFBM} \\
\end{array}
\right.\end{equation}
The process $(U_t^B)_{t\geq 0}$ has an inverse $ (U_t^B)^{-1}=B^{-1}U_t^{-1}$ at any time $t\geq 0$. From \eqref{uAFBM}, we know that $((U_t^B)^{-1})_{t\geq 0}$ is the bounded adapted semi-circular process defined by the free stochastic differential equation
\begin{equation}
\left\{
\begin{array}{r c l}
   \left(U_0^B\right)^{-1} &=& B^{-1}, \\
    \diff \left(U_t^B\right)^{-1}&=& -i \left(U_t^B\right)^{-1}\diff X_t -\frac{1}{2}\left(U_t^B\right)^{-1} \diff t. \label{uBIFBM} \\
\end{array}
\right.\end{equation}
Similarly, from Lemma~\ref{fmart} we know that $((U_t^B)^*)_{t\geq 0}$ is a bounded adapted semi-circular process defined by the free stochastic differential equation
\begin{equation}
\left\{
\begin{array}{r c l}
   \left(U_0^B\right)^* &=& B^*, \\
    \diff \left(U_t^B\right)^*&=& -i \left(U_t^B\right)^*\diff X_t -\frac{1}{2}\left(U_t^B\right)^* \diff t, \label{uABFBM} \\
\end{array}
\right.\end{equation}
and that $({(U_t^B)^*}^{-1})_{t\geq 0}$ is a bounded adapted semi-circular process defined by the free stochastic differential equation
\begin{equation}
\left\{
\begin{array}{r c l}
     {\left(U_0^B\right)^*}^{-1} &=& {B^*}^{-1}, \\
    \diff {\left(U_t^B\right)^*}^{-1} &=& i \diff X_t  {\left(U_t^B\right)^*}^{-1} -\frac{1}{2}{\left(U_t^B\right)^*}^{-1} \diff t. \label{uABIFBM} \\
\end{array}
\right.\end{equation}

Let $k\in \mathbb{N}$, and $U_1,\ldots,U_k\in \{{U_t^B},{U_t^B}^*,{U_t^B}^{-1},{{U_t^B}^*}^{-1}\}$.\\
For all $1\leq i\leq k$, set $l (i)=0$, $r(i)=1$ and $\epsilon(i)=1$ if $U_i={U_t^B}$ or $U_i={{U_t^B}^*}^{-1}$, and $l (i)=1$, $r(i)=0$ and $\epsilon(i)=-1$ if $U_i={U_t^B}^*$ or $U_i={U_t^B}^{-1}$. For all $1\leq i<j\leq k$, the term $\overbracket{U_i\cdots U_j}$ refers to the product $U^{r(i)}_iU_{i+1}\cdots U_{j-1}U_j^{l(j)}$, that is to say the product $U_i\cdots U_j$  possibly excluding $U_i$ and/or $U_j$. We claim
\begin{equation}
\diff \left(U_1\cdots U_k\right)=\displaystyle\sum_{i=1}^k  U_1\cdots \diff U_i \cdots U_k
-\hspace{-0.3cm}\displaystyle\sum_{1\leq i< j\leq k}\hspace{-0.3cm} \epsilon(i)\epsilon(j) U_1\cdots  \widehat{\overbracket{U_i\cdots U_j}}\cdots U_k\cdot \tau \left(\overbracket{U_i\cdots U_j}\right)\diff t.\label{evsto}
\end{equation}
Indeed, using Itô's formula (see Section~\ref{Ito}), the equation \eqref{evsto} follows by induction on $k\in \mathbb{N}$. 
Let us show how it works. We suppose that the formula is true for $k\in \mathbb{N}$. Let $U_{k+1}={U_t^B}^{-1}$. We use the It\^{o} formula and evolution equations \eqref{uBFBM}, \eqref{uBIFBM}, \eqref{uABFBM} and \eqref{uABIFBM} to infer
$$\begin{array}{rcl}& &\hspace{-1.5cm}\diff \left(U_1\cdots U_k U_{k+1}\right)\\
&=&\displaystyle\sum_{i=1}^k  U_1\cdots \diff U_i \cdots U_{k+1}
-\displaystyle\sum_{1\leq i< j\leq k} \epsilon(i)\epsilon(j) U_1\cdots  \widehat{\overbracket{U_i\cdots U_j}}\cdots U_{k+1}\cdot \tau \left(\overbracket{U_i\cdots U_j}\right)\diff t\\
&&\hspace{1cm}+U_1\cdots U_k \diff U_{k+1}+\displaystyle\sum_{1\leq i\leq k} \epsilon(i) U_1\cdots  \widehat{\overbracket{U_i\cdots U_{k+1}}}\cdot \tau \left(\overbracket{U_i\cdots U_{k+1}}\right)\diff t,\\
\end{array}$$
which is the formula at step $k+1$.
The other cases $U_{k+1}=U_t^B,{U_t^B}^*$ or ${{U_t^B}^*}^{-1}$ are treated similarly, justifying \eqref{evsto}.

Evaluating the conditional expectation $\tau\left(\left.\cdot\right|\mathcal{B}\right)$ on both sides of \eqref{evsto}, and using Lemma~\ref{fmart} (i.e. that the conditional expectation of a stochastic integral with respect to $(X_t)_{t\geq 0}$ vanishes), we have
$$\begin{array}{rcl}
\displaystyle\frac{\diff}{\diff t}\tau\left(U_1\cdots U_k|\mathcal{B}\right)&=&-\displaystyle\frac{k}{2}  \tau\left(\left.U_1\cdots U_k\right|\mathcal{B}\right)\\
&&\hspace{1cm}-2\displaystyle\sum_{1\leq i< j\leq k}\epsilon(i)\epsilon(j)\tau \left(\left.U_1\cdots  \widehat{\overbracket{U_i\cdots U_j}}\cdots U_{k+1}\right|\mathcal{B}\right) \tau \left(\overbracket{U_i\cdots U_j}\right).
\end{array}$$
Equivalently, for all $k\in \mathbb{N}$ and $X_1,\ldots,X_k\in \{X,X^*,X^{-1},{X^*}^{-1}\}$, we have
$$\frac{\diff}{\diff t}\tau\left(X_1\cdots X_k\left(U_t^B\right)|\mathcal{B}\right)=\tau\left(\frac{1}{2}\Delta_U X_1\cdots X_k\left(U_t^B\right)|\mathcal{B}\right) .$$
We extend this identity to all $P\in\mathbb{C}\{X,X^*,X^{-1},{X^*}^{-1}\}$ by linearity and by the following induction. If $P$ and $Q\in\mathbb{C}\{X,X^*,X^{-1},{X^*}^{-1}\}$ verify the lemma, we have
$$\begin{array}{rl}
\displaystyle\frac{\diff}{\diff t} \tau\left(\left(P\tr Q\right)\left(U_t^B\right)|\mathcal{B}\right)&=\displaystyle\frac{\diff}{\diff t}\left( \tau\left(P(U_t^B)|\mathcal{B}\right)\tau\left(Q\left(U_t^B\right)\right)\right) \\
&=\left(\displaystyle\frac{\diff}{\diff t} \tau\left( P(U_t^B)|\mathcal{B}\right)\right)\tau\left(Q\left(U_t^B\right)\right)+\tau\left(P(U_t^B)|\mathcal{B}\right)\left(\displaystyle\frac{\diff}{\diff t}\tau\left( Q\left(U_t^B\right)\right)\right)\\
&= \tau\left(\Delta_U P(U_t^B)|\mathcal{B}\right)\tau\left(Q\left(U_t^B\right)\right)+\tau\left(P(U_t^B)|\mathcal{B}\right)\tau\left(\left(\displaystyle\frac{\diff}{\diff t}\tau\left( Q\left(U_t^B\right)|\mathcal{B}\right)\right)\right)\\
&= \tau\left(\displaystyle\frac{1}{2}\Delta_U P(U_t^B)|\mathcal{B}\right)\tau\left(Q\left(U_t^B\right)\right)+\tau\left(P(U_t^B)|\mathcal{B}\right)\tau\left(\tau\left(\displaystyle\frac{1}{2}\Delta_U Q\left(U_t^B\right)|\mathcal{B}\right)\right)\\
&= \tau\left(\left(\left(\displaystyle\frac{1}{2}\Delta_U P\right)\tr Q+P\tr \left(\displaystyle\frac{1}{2}\Delta_U Q\right)\right)(U_t^B)|\mathcal{B}\right)\\
&=\tau\left(\displaystyle\frac{1}{2}\Delta_U\left(P\tr Q\right)\left(U_t^B\right)|\mathcal{B}\right).\qedhere
\end{array}$$
\end{proof}

Since, for all $m\in \mathbb{N}$, $\frac{t}{2} \Delta_U$ leaves $\mathbb{C}_m\{X,X^*,X^{-1},{X^*}^{-1}\}$ invariant, let us define the endomorphism $e^{\frac{t}{2}\Delta_U}$ on $\mathbb{C}\{X,X^*,X^{-1},{X^*}^{-1}\}=\bigcup_{m\in \mathbb{N}}\mathbb{C}_m\{X,X^*,X^{-1},{X^*}^{-1}\}$ by the convergent series $\sum_{n=0}^{\infty}\frac{1}{n!} \left(\frac{t}{2}\right)^n \Delta_U^n$ on each finite dimensional space $\mathbb{C}_m\{X,X^*,X^{-1},{X^*}^{-1}\}$.

\begin{proposition}
\label{distU}Let $t\geq 0$ and $P\in\mathbb{C}\{X,X^*,X^{-1},{X^*}^{-1}\}$. For all invertible $B\in \mathcal{B}$, we have
$$\tau\left(\left.P\left(U_tB\right)\right|\mathcal{B}\right)=\left(e^{\frac{t}{2}\Delta_U}P\right)(B).$$
\end{proposition}
\begin{proof}
Let $B$ be an invertible element of $\mathcal{B}$. Let $P\in \mathbb{C}\{X_i:i\in I\}$. There exists a finite index set $J\subset I$ and $d\in \mathbb{N}$ such that $P\in\mathbb{C}_d\{X_i:i\in J\}$. We remark that $\frac{t}{2}\Delta_U$ is a linear operator on $\mathbb{C}_d\{X_i:i\in J\}$ and that  $\left(P\mapsto \tau\left(\left.P\left(U_tB\right)\right|\mathcal{B}\right)\right)_{t\geq 0}$ are linear maps from $\mathbb{C}_d\{X_i:i\in J\}$ to $\mathcal{A}$. Consequently, we deduce Lemma~\ref{distU} from Lemma~\ref{eqdiff} and Lemma~\ref{difftau}.
\end{proof}

\subsection{Free circular multiplicative Brownian motion} Following~\cite{Biane1997a}, let us define the (right) free circular multiplicative Brownian motion. It is the unique bounded adapted circular process $(G_t)_{t\geq 0}$ defined by the free stochastic differential equation\label{FBMonGL}
\begin{equation}
\left\{
\begin{array}{c c c}
    G_0 &=& \Id, \\
    \diff G_t &=& \diff Z_t G_t. \label{FBM} \\
\end{array}
\right.\end{equation}
According to~\cite{Biane1997a}, the process $(G_t)_{t\geq 0}$ has an inverse at any time. It is the bounded adapted circular process $(G_t^{-1})_{t\geq 0}$ defined by the free stochastic differential equation
\begin{equation}
\left\{
\begin{array}{c c c}
    G_0^{-1} &=& \Id, \\
    \diff G_t^{-1} &=&-G_t^{-1} \diff Z_t .\label{iFBM} \\
\end{array}
\right.\end{equation}
From Lemma~\ref{fmart}, we know that $(G_t^*)_{t\geq 0}$ is a bounded adapted circular process defined by the free stochastic differential equation
\begin{equation}
\left\{
\begin{array}{r c l}
   G_t^* &=& \Id, \\
    \diff G_t^*&=& G_t^*\diff Z_t^*, \label{AFBM} \\
\end{array}
\right.\end{equation}
and that $({G_t^*}^{-1})_{t\geq 0}$ is a bounded adapted semi-circular process defined by the free stochastic differential equation
\begin{equation}
\left\{
\begin{array}{r c l}
     {G_t^*}^{-1} &=& \Id, \\
    \diff {G_t^*}^{-1} &=& -\diff Z_t^*  {G_t^*}^{-1}  . \label{AIFBM} \\
\end{array}
\right.\end{equation}

\subsubsection{Extension of $\D_{G_t}$}We establish in this section a version of Proposition~\ref{distU} for $\left(G_t\right)_{t\geq 0}$.

Let us define a derivation $\Delta_{GL}$ on $(\mathbb{C}\{X,X^*,X^{-1},{X^*}^{-1}\},\cdot_{\tr})$ analogously to the definition of $\Delta_{U}$ in Section~\ref{deltau}. Let $k\in \mathbb{N}$, and $X_1,\ldots,X_k\in \{X,X^*,X^{-1},{X^*}^{-1}\}$. For all $1\leq i\leq k$, set $l (i)=0$, $r(i)=1$ if $X_i=X$ or $X_i={X^*}^{-1}$, and $l(i)=1$, $r(i)=0$ if $X_i=X^*$ or $X_i=X^{-1}$. For all $1\leq i<j\leq k$, the term $\overbracket{X_i\cdots X_j}$ refers to the product $X^{r(i)}_iX_{i+1}\cdots X_{j-1}X_j^{l(j)}$, that is to say the product $X_i\cdots X_j$  possibly excluding $X_i$ and/or $X_j$. For all $1\leq i\leq k$, set $\epsilon(i)=1$ if $X_i=X$ or $X_i={X^*}$, and $\epsilon(i)=-1$ if $X_i={X^*}^{-1}$ or $X_i=X^{-1}$. For all $1\leq i<j\leq k$, set $\delta(i,j)=0$ if $X_i,X_j\in \{X,X^{-1}\}$ or if $X_i,X_j\in \{X^*,{X^*}^{-1}\}$, and $\delta(i,j)=1$ otherwise. We set
$$\Delta_{GL} X_1\cdots X_k=4\displaystyle\sum_{1\leq i< j\leq k} \delta(i,j)\epsilon(i)\epsilon(j)
\cdot X_1\cdots  \widehat{\overbracket{X_i\cdots X_j}}\cdots X_k \tr \left(\overbracket{X_i\cdots X_j}\right).
$$
We extend $\Delta_{GL}$ to all $\mathbb{C}\{X,X^*,X^{-1},{X^*}^{-1}\}$ by linearity and by the relation\label{deltag}
\begin{equation}
\forall P,Q \in \mathbb{C}\{X,X^*,X^{-1},{X^*}^{-1}\}, \Delta_{GL} \left(P\tr Q\right)=\left(\Delta_{GL} P\right)\tr Q+P\tr \left(\Delta_{GL} Q\right).\label{structugl}
\end{equation}

\begin{lemma}
For all $P\in\mathbb{C}\{X,X^*,X^{-1},{X^*}^{-1}\}$,\label{difftauG}
$$\frac{\diff}{\diff t} \tau\left(P\left(G_t\right)\right)=\tau\left(\frac{1}{4}\Delta_{GL} P\left(G_t \right)\right) .$$
\end{lemma}
\begin{proof}
Recall that $(G_t)_{t\geq 0}$, $(G_t^{-1})_{t\geq 0}$, $(G_t^*)_{t\geq 0}$ and $({G_t^*}^{-1})_{t\geq 0}$ are defined respectively by \eqref{FBM}, \eqref{iFBM}, \eqref{AFBM} and \eqref{AIFBM}.

Let $k\in \mathbb{N}$, and $G_1,\ldots,G_k\in \{{G_t},{G_t}^*,{G_t}^{-1},{{G_t}^*}^{-1}\}$. For all $1\leq i\leq k$, set $l (i)=0$, $r(i)=1$ if $G_i=G_t$ or $G_i={G_t^*}^{-1}$, and $l(i)=1$, $r(i)=0$ if $G_i=G_t^*$ or $G_i=G_t^{-1}$.\\
For all $1\leq i\leq k$, set $\epsilon(i)=1$ if $G_i=G_t$ or $X_i={G_t^*}$, and $\epsilon(i)=-1$ if $G_i={G_t^*}^{-1}$ or $G_i=G_t^{-1}$. For all $1\leq i<j\leq k$, set $\delta(i,j)=0$ if $G_i,G_j\in \{G_t,G_t^{-1}\}$ or if $G_i,G_j\in \{G_t^*,{G_t^*}^{-1}\}$, and $\delta(i,j)=1$ otherwise. We claim
\begin{multline}\label{evstog}
\diff \left(G_1\cdots G_k\right)\\=\displaystyle\sum_{i=1}^k  G_1\cdots \diff G_i \cdots G_k
+\displaystyle\sum_{1\leq i< j\leq k}\delta(i,j)\epsilon(i)\epsilon(j)
\cdot G_1\cdots  \widehat{\overbracket{G_i\cdots G_j}}\cdots G_k \tr \left(\overbracket{G_i\cdots G_j}\right)\diff t.
\end{multline}
Indeed, using Itô's formula~\ref{Ito}, the equation \eqref{evstog} follows by induction on $k\in \mathbb{N}$. One can see the proof of Lemma~\ref{difftau} to understand how it works, analogously to the case of $U_t$.

Evaluating the conditional expectation $\tau$ on both sides of this equation, and using the fact that the trace of a stochastic integral with respect to $(Z_t)_{t\geq 0}$ vanishes  (Lemma~\ref{fmart}), we have
$$\frac{\diff}{\diff t}\tau\left(G_1\cdots G_k\right) =\displaystyle\sum_{1\leq i< j\leq k}\delta(i,j)\epsilon(i)\epsilon(j)
\cdot \tau \left(G_1\cdots  \widehat{\overbracket{G_i\cdots G_j}}\cdots G_k \right)\tau \left(\overbracket{G_i\cdots G_j}\right).
$$
Equivalently, for all $k\in \mathbb{N}$, and $X_1,\ldots,X_k\in \{X,X^*,X^{-1},{X^*}^{-1}\}$, we have
$$\frac{\diff}{\diff t}\tau\left(X_1\cdots X_k\left(G_t\right)\right)=\tau\left(\frac{1}{4}\Delta_{GL} X_1\cdots X_k\left(G_t\right)\right) .$$
We extend this identity to all $P\in\mathbb{C}\{X,X^*,X^{-1},{X^*}^{-1}\}$ by linearity and by the same induction as the end of the proof of Lemma~\ref{difftau}.
\end{proof}

Let $t\geq 0$. Since, for all $m\in \mathbb{N}$, $\frac{t}{2} \Delta_{GL}$ leaves $\mathbb{C}_m\{X,X^*,X^{-1},{X^*}^{-1}\}$ invariant, the endomorphism $e^{\frac{t}{4}\Delta_{GL}}$ on $\mathbb{C}\{X,X^*,X^{-1},{X^*}^{-1}\}=\bigcup_{m\in \mathbb{N}}\mathbb{C}_m\{X,X^*,X^{-1},{X^*}^{-1}\}$ is defined to be the convergent series $\sum_{n=0}^{\infty}\frac{1}{n!} \left(\frac{t}{4}\right)^n (\Delta_{GL})^n$ on each finite dimensional space $\mathbb{C}_m\{X,X^*,X^{-1},{X^*}^{-1}\}$.

\begin{proposition}
Let $t\geq 0$ and $P\in \mathbb{C}\{X,X^*,X^{-1},{X^*}^{-1}\}$. We have\label{distG}
$$\tau\left(P\left(G_t\right)\right)=(e^{\frac{t}{4} \Delta_{GL}}P)(1).$$
\end{proposition}
This proposition allows us to compute the distribution of $G_t$. For example, let us remark that $\Delta_{GL}$ vanishes on $\mathbb{C}\{X\}$. This implies that $e^{\frac{t}{4} \Delta_{GL}}=\id$ on $\mathbb{C}\{X\}$, and that, for all $P\in \mathbb{C}\{X\}$, we have $\tau\left(P\left(G_t\right)\right)=(e^{\frac{t}{4} \Delta_{GL}}P)(1)=P(1)$.
\begin{proof}
The demonstration follows the proof of Proposition~\ref{distU}: since the maps $\left(P\mapsto \tau\left(P\left(G_t\right)\right)\right)_{t\geq 0}$ are linear, we deduce Proposition~\ref{distG} from Lemma~\ref{eqdiff} and Lemma~\ref{difftauG}.
\end{proof}

\subsection{Free Hall transform}\label{hallfree}

Let $\left(U_t\right)_{t\geq 0}$ be a free unitary Brownian motion, and let $\left(G_t\right)_{t\geq 0}$ be a free circular multiplicative Brownian motion. Let $t\geq 0$. We denote by $L^2(U_t,\tau)$ the Hilbert completion of the $*$-algebra generated by $U_t$ for the norm $\|\cdot\|_2:A\mapsto \tau\left(A^* A\right)^{1/2}$, and by $L^2_{\hol}(G_t,\tau)$ the Hilbert completion of the algebra generated by $G_t$ and $G_t^{-1}$ for the norm $\|\cdot\|_2:A\mapsto \tau\left(A^* A\right)^{1/2}$.

\subsubsection{Definition}
For all $t\geq 0$,  let us define the polynomials $\left(P_n^t\right)_{n\in \mathbb{N}}$ by the generating series
$$\displaystyle\sum_{n=0}^\infty z^n P_n^t(x)
=\frac{1}{1-ze^{\frac{t}{2}\left(\frac{1+z}{1-z}\right)}x}.$$
We remark that, for all $t\geq 0$ and $n\in \mathbb{N}$, the degree of $P_n^t$ is $n$ and the leading coefficient is the coefficient of $x^nz^n$ in the development and is therefore equal to $e^{nt/2}$.

A quick way to define the free Hall transform is to define it on Laurent polynomials. Let us denote by $\mathbb{C}[X,X^{-1}]$ the space of Laurent polynomials. For all $t\geq 0$, let us define $\mathcal{G}_t$ by the unique linear operator on $\mathbb{C}[X,X^{-1}]$ such that, for all $n\in \mathbb{N}$, one has $\mathcal{G}_t\left(P_n^t(X)\right)=X^n$ and $\mathcal{G}_t\left(P_n^t(X^{-1})\right)=X^{-n}$.
The following theorem, due to Biane, corresponds to Theorem~9 and Lemma~18 of~\cite{Biane1997}.
\begin{theo*}[Biane~\cite{Biane1997}]
Let $t> 0$. The map $\mathcal{F}_t:P \big(U_t \big) \mapsto \mathcal{G}_t(P) \big(G_t \big)  $ for all $P\in \mathbb{C}[X,X^{-1}]$ is an isometric map which extends to a Hilbert space isomorphism $\mathcal{F}_t$ between $L^2(U_t,\tau)$ and $L^2_{\hol}(G_t,\tau)$ called the free Hall transform.
\end{theo*}
We notice that the polynomials $\left(P_n^t\right)_{n\in \mathbb{N}}$ are defined in~\cite{Biane1997} by another generating series
$$\sum_{n=0}^\infty z^n P_n^t(x)=\frac{1}{1-\frac{z}{1+z}e^{\frac{t}{2}(1+2z)}x},$$but Philippe Biane kindly pointed out (personal communication, 2012) that it is necessary to replace this generating series by the series $$\frac{1}{1-ze^{\frac{t}{2}\left(\frac{1+z}{1-z}\right)}x}$$ for the proofs of Lemmas~18 and~19 of~\cite{Biane1997} to be correct.

\subsubsection{Another construction}Theorem~\ref{fht} has to be read in parallel with the classical Hall transform definition, which will be stated in Section~\ref{hallclas}.
\begin{theorem}
\label{fht}Let $\left(U_t\right)_{t\geq 0}$ be a free unitary Brownian motion, and let $\left(G_t\right)_{t\geq 0}$ be a free circular multiplicative Brownian motion. Let $t> 0$.

For all $P\in \mathbb{C}\{X,X^{-1}\}$, $$\mathcal{F}_t:P (U_t ) \mapsto (e^{\frac{t}{2}\Delta_U}P) (G_t )  $$
is an isometric map which extends to an Hilbert isomorphism $\mathcal{F}_t$ between $L^2(U_t,\tau)$ and $L^2_{\hol}(G_t,\tau)$. Moreover, this isomorphism is the free Hall transform.

In particular, if $\left(U_t\right)_{t\geq 0}$ and $\left(G_t\right)_{t\geq 0}$ are free, for all $P\in \mathbb{C}\{X,X^{-1}\}$,
$$\mathcal{F}_t\left(P(U_t)\right) =\tau\left(\left.P(U_t G_t)\right|G_t\right).$$
\end{theorem}
It should be remarked that, for the map $\mathcal{F}_t$ to be well-defined, it must be true that, for all $P, Q\in \mathbb{C}\{X,X^{-1}\}$, if $P \left(U_t \right)=Q \left(U_t \right)$, then $(e^{\frac{t}{2}\Delta_U}P) (G_t )=(e^{\frac{t}{2}\Delta_U}Q) (G_t )$. This fact is contained in the proof below.

Theorem~\ref{fht} allows us to compute explicitly the free Hall transform on the *-algebra generated by $U_t$. For example, we have from Section~\ref{utex} that
$
e^{\frac{t}{2}\Delta_U}X^2=e^{\D_{U_t}}X^2
= e^{-2t}\left(X^2-t X \tr X\right)
$ 
from which we deduce that
$
\mathcal{F}_t(U_t^2)
= e^{-2t}\left(G_t^2-t G_t \tau (G_t)\right).
$
We computed in the previous section that, for all $P\in \mathbb{C}[X]$, we have $\tau(P(G_t))=P(1)$. Thus, we have
$
\mathcal{F}_t(U_t^2)
= e^{-2t}G_t^2-t e^{-2t} G_t .
$
\begin{proof}
We will prove in a first step that, for all $P\in \mathbb{C}\{X,X^{-1}\}$, we have
\begin{equation}\left\|P(U_t)\right\|^2_{L^2(U_t,\tau)}=\left\| e^{\frac{t}{2}\Delta_U}P(G_t)\right\|^2_{L^2_{\hol}(G_t,\tau)} \label{auxil}\end{equation}
This proves that, for all $P, Q\in \mathbb{C}\{X,X^{-1}\}$, if $P \left(U_t \right)=Q \left(U_t \right)$, then \begin{multline*}\left\| (e^{\frac{t}{2}\Delta_U}P) (G_t )-(e^{\frac{t}{2}\Delta_U}Q) (G_t )\right\|^2_{L^2_{\hol}(G_t,\tau)}\\=\left\| (e^{\frac{t}{2}\Delta_U}(P-Q)) (G_t )\right\|^2_{L^2_{\hol}(G_t,\tau)}=\left\|(P-Q)(U_t)\right\|^2_{L^2(U_t,\tau)}=0,\end{multline*}
and consequently $(e^{\frac{t}{2}\Delta_U}P) (G_t )=(e^{\frac{t}{2}\Delta_U}Q) (G_t )$. Thus, $\mathcal{F}_t$ is a well-defined isometric map, and so it extends to a Hilbert space isomorphism $\mathcal{F}$ between $L^2(U_t,\tau)$ and $\mathcal{F}(L^2(U_t,\tau))\subset L^2_{\hol}(G_t,\tau)$. The surjectivity is clear since, for all $P\in \mathbb{C}\{X,X^{-1}\}$,
$$P (G_t )=\left(e^{\frac{t}{2}\Delta_U}\left(e^{-\frac{t}{2}\Delta_U}P\right)\right) (G_t ) =\mathcal{F}_t\left( \left(e^{-\frac{t}{2}\Delta_U}P\right) (U_t )\right).$$
The equality with the free Hall transform will be made in a second step, and the last part of Theorem~\ref{fht} follows Proposition~\ref{distU}, which says that, for all $P\in \mathbb{C}\{X,X^{-1}\}$,
$$(e^{\frac{t}{2}\Delta_U}P)(G_t)=\tau\left(\left.P(U_t G_t)\right|G_t\right).$$

\subsubsection*{Step $1$}With the aim of proving \eqref{auxil}, we work on the *-algebra $\mathbb{C}\{X,X^{-1},X^*,{X^*}^{-1}\}$. Let us define the subspace $\mathcal{C}$ of $\mathbb{C}\{X,X^{-1},X^*,{X^*}^{-1}\}$ which is generated by elements of
$$\left\{P_0 \tr(P_1)\cdots \tr(P_n): n\in \mathbb{N}, P_0,\ldots,P_n \text{ are monomials of }\mathbb{C}[X,X^{-1}]  \right\} .$$
Replacing the pairs $X X^{-1}=X^{-1}X$ by $1$, we observe that, for all $P\in\mathbb{C}\{X,X^{-1},X^*,{X^*}^{-1}\}$ there exists $\tilde{P}\in \mathcal{C}$ such that
$\left\|P(U_t)\right\|^2_{L^2(U_t,\tau)} =\|\tilde{P}(U_t)\|^2_{L^2(U_t,\tau)}.$
Moreover,
\begin{multline*}\left\| \left(e^{\frac{t}{2}\Delta_U}P\right)(G_t)\right\|^2_{L^2_{\hol}(G_t,\tau)} \\=\left\| \tau(\left.P(U_t G_t)\right|G_t)\right\|^2_{L^2_{\hol}(G_t,\tau)}  =\left\| \tau(\left.\tilde{P}(U_t G_t)\right|G_t)\right\|^2_{L^2_{\hol}(G_t,\tau)}=\left\| \left(e^{\frac{t}{2}\Delta_U}\tilde{P}\right)(G_t)\right\|^2_{L^2_{\hol}(G_t,\tau)}. \end{multline*}
Thus, it remains to prove that, for all $P\in \mathcal{C}$, we have $\left\|P(U_t)\right\|^2_{L^2(U_t,\tau)}=\| e^{\frac{t}{2}\Delta_U}P(G_t)\|^2_{L^2_{\hol}(G_t,\tau)}$.

Let $P\in \mathcal{C}$. Proposition~\ref{distU} gives us that
$$\begin{array}{rcl}
\left\|P(U_t)\right\|^2_{L^2(U_t,\tau)}&=&\tau\Big( P(U_t)P(U_t)^*\Big)\\
&=&\tau\Big((PP^*)(U_t)\Big)\\
&=&\left(e^{\frac{t}{2}\Delta_U}\left(PP^*\right)\right)(1)
\end{array}$$
and Proposition~\ref{distG} gives us that
$$\begin{array}{rcl}
\left\| e^{\frac{t}{2}\Delta_U}P(G_t)\right\|^2_{L^2_{\hol}(G_t,\tau)} &=&\tau\left( \Big(e^{\frac{t}{2}\Delta_U}P(G_t)\Big) \Big(e^{\frac{t}{2}\Delta_U}P(G_t)\Big)^*\right)\\
&=&\tau\left(\Big(e^{\frac{t}{2}\Delta_U}P\left(e^{\frac{t}{2}\Delta_U}P\right)^*\Big)(G_t)\right)\\
&=&\left(e^{\frac{t}{4}\Delta_{GL}}\Big(e^{\frac{t}{2}\Delta_U}P\left(e^{\frac{t}{2}\Delta_U}P\right)^*\Big)\right)(1).
\end{array}$$
It remains to prove that, for all $P\in \mathcal{C}$, $e^{\frac{t}{2}\Delta_U}\left(P^*P\right)=e^{\frac{t}{4}\Delta_{GL}}\Big(e^{\frac{t}{2}\Delta_U}P\left(e^{\frac{t}{2}\Delta_U}P\right)^*\Big).$

Let us define the subspace $\mathcal{E}$ of $\mathbb{C}\{X,X^{-1},X^*,{X^*}^{-1}\}$ which is generated by elements of
$$\left\{P_0Q_0^* \tr(P_1Q_1^*)\cdots \tr(P_nQ_n^*): n\in \mathbb{N}, P_0,Q_0,\ldots,P_n,Q_n \text{ are monomials of }\mathbb{C}[X,X^{-1}]  \right\} .$$
We remark that the space $\mathcal{E}$ is invariant by the product $\cdot_{\tr}$ and by the operators $\Delta_U$ and $\Delta_{GL}$. Moreover, the algebra $\left(\mathcal{E},\cdot_{\tr}\right)$ is generated as an algebra by the elements of $$\left\{PQ^*: P,Q \text{ monomials of }\mathbb{C}[X,X^{-1}]  \right\}.$$
Let us define the auxiliary derivations $\Delta_U^+$ and $\Delta_U^-$ on $\left(\mathcal{E},\cdot_{\tr}\right)$ in the following way. For all $P,Q$ monomials of $\mathbb{C}[X,X^{-1}]$, we set $\Delta_U^+(PQ^*)=(\Delta_UP)\cdot Q^*$ and $\Delta_U^-(PQ^*)=P\cdot (\Delta_UQ)^*$, and we extend $\Delta_U^+$ and $\Delta_U^-$ to all $\mathcal{E}$ by linearity and by the relations
$$\forall P,Q \in \mathcal{E}, \Delta_U^\pm \left(P\tr Q\right)=\left(\Delta_U^\pm P\right)\tr Q+P\tr \left(\Delta_U^\pm Q\right).$$
It should be remarked that the operators are well-defined because $X_i$ and $X_j^*$ are algebraically free: thus, if $P_1Q_1^*=P_2Q_2^*$, then there is a non-zero constant $\lambda$ so that $P_1=\lambda P_1$ and $Q_2=\bar{\lambda}Q_1$.

The operators $\Delta_U^+$ and $\Delta_U^-$ are such that, for all $P\in \mathcal{C}$, we have $$\left(e^{\frac{t}{2}\Delta_U}P\right)\left(e^{\frac{t}{2}\Delta_U}P\right)^*=e^{\frac{t}{2}\Delta_U^+}\left(P(e^{\frac{t}{2}\Delta_U}P^*)\right)=e^{\frac{t}{2}\Delta_U^+}\left(e^{\frac{t}{2}\Delta_U^-}(PP^*)\right).$$
Let us observe that, for all $P\in \mathcal{C}$, we have by linearity that $PP^*\in \mathcal{E}$. We can now reformulate what we wanted to prove. It remains to verify that, for all $Q\in \mathcal{E}$, $e^{\frac{t}{2}\Delta_U}Q=e^{\frac{t}{4}\Delta_{GL}}e^{\frac{t}{2}\Delta_U^+}e^{\frac{t}{2}\Delta_U^-}Q$. More precisely, the end of the demonstration is devoted to proving that $\Delta_U= \Delta_{GL}/2+\Delta_U^++\Delta_U^-$, and that $\Delta_{GL}$, $\Delta_U^+$ and $\Delta_U^-$ commute on $\mathcal{E}$ (they do not commute on all $\mathbb{C}\{X,X^{-1},X^*,{X^*}^{-1}\}$).

Because $\Delta_U$,$\Delta_{GL}$, $\Delta_U^+$ and $\Delta_U^-$ are derivations on $\mathcal{E}$, it suffices to verify the properties on elements of $\left\{PQ^*: P,Q \text{ monomials of }\mathbb{C}[X,X^{-1}]  \right\}.$ For example, let $m,n\in \mathbb{N}$, and set $P=X^n$ and $Q=X^{-m}$. We have
$$\begin{array}{rcrl}
\Delta_U(PQ^*)&=&(m+n)PQ^*&-2\displaystyle\sum_{1\leq i<j\leq n}X^{n-i+j}\tr(X^{i-j})Q^*\\
  & & &-2\displaystyle\sum_{1\leq i<j\leq m}P( {X^{-1}}^{*})^{m-i+j}\tr(({X^{-1}}^{*})^{i-j})\\
  & & &+2\displaystyle\sum_{\substack{1\leq i\leq n\\ 1\leq j \leq m}}X^{i-1}({X^{-1}}^{*})^{m-j}\tr\left(X^{n+1-i} ({X^{-1}}^{*})^{j}\right)\\
  &=&\Bigg( nP&-2\displaystyle\sum_{1\leq i<j\leq n}X^{n-i+j}\tr(X^{i-j})\Bigg)Q^*\\
  &  &+P\Bigg( mQ^*&-2\displaystyle\sum_{1\leq i<j\leq m}( {X^{-1}}^{*})^{m-i+j}\tr(({X^{-1}}^{*})^{i-j})\Bigg)\\
  & & +\displaystyle\frac{1}{2}\Delta_{GL}(PQ^*)&\\
&=&\displaystyle\frac{1}{2} \Delta_{GL}(PQ^*)&+\Delta_U^+(PQ^*)+\Delta_U^-(PQ^*).
\end{array}$$
The other three cases (depending on the signs of the exponents of $P$ and $Q$) are treated similarly, and we have $\Delta_U= \Delta_{GL}/2+\Delta_U^++\Delta_U^-$.

The commutativity is less obvious. For all $m,n\in \mathbb{N}$, and $P=X^n$ and $Q=X^{-m}$, we have
$$\begin{array}{rcl}
\Delta_U^+(PQ^*)&=&nPQ^*-2\displaystyle\sum_{1\leq i<j\leq n}X^{n-i+j}\tr(X^{i-j})Q^*\\
&=&nPQ^*-2\displaystyle\sum_{\substack{1\leq i,j \leq n\\ i+j=n}}i X^i\tr(X^{j})Q^*.\\
\end{array}$$
Let us fix $m,n\in \mathbb{N}$, and $P=X^n$ and $Q=X^{-m}$. We have
$$\begin{array}{rcl}
\Delta_U^+\Delta_{GL}(PQ^*)&=&\Delta_U^+\left(-4 \displaystyle\sum_{\substack{1\leq i\leq n\\ 1\leq l \leq m}}X^{i-1}({X^{-1}}^{*})^{m-l}\tr\left(X^{n+1-i} ({X^{-1}}^{*})^{l}\right)\right)\\
&=&\vspace{-0.3cm}-4n \displaystyle\sum_{\substack{1\leq i\leq n\\ 1\leq l \leq m}}X^{i-1}({X^{-1}}^{*})^{m-l}\tr\left(X^{n+1-i} ({X^{-1}}^{*})^{l}\right)\\
& &\vspace{-0.3cm}\hspace{1cm}+8\overbrace{\displaystyle\sum_{\substack{1\leq i\leq n\\ 1\leq l \leq m}}\displaystyle\sum_{\substack{1\leq j,k \leq i-1\\ j+k=i-1}}j X^j\tr(X^{k})({X^{-1}}^{*})^{m-l}\tr\left(X^{n+1-i} ({X^{-1}}^{*})^{l}\right)}^{=\displaystyle\Sigma_1}\\

& &\hspace{1cm}+8\overbrace{\displaystyle\sum_{\substack{1\leq i\leq n\\ 1\leq l \leq m}}\displaystyle\sum_{\substack{1\leq j,k \leq n+1-i\\ j+k=n+1-i}}j X^{i-1}({X^{-1}}^{*})^{m-l}\tr\left(X^{j} ({X^{-1}}^{*})^{l}\right)\tr\left(X^k\right)}^{=\displaystyle\Sigma_2},
\end{array}$$
while
$$\begin{array}{rcl}
\Delta_{GL}\Delta_U^+(PQ^*)&=
&\Delta_{GL}\left(nPQ^*-2\displaystyle\sum_{\substack{1\leq i,j \leq n\\ i+j=n}}i X^i\tr(X^{j})Q^*\right)\\
&=&\vspace{-0.3cm}-4n \displaystyle\sum_{\substack{1\leq i\leq n\\ 1\leq l \leq m}}X^{i-1}({X^{-1}}^{*})^{m-l}\tr\left(X^{n+1-i} ({X^{-1}}^{*})^{l}\right)\\
& &\hspace{1cm}+8\overbrace{\displaystyle\sum_{\substack{1\leq i,j \leq n\\ i+j=n}}\displaystyle\sum_{\substack{1\leq k\leq i\\ 1\leq l \leq m}}i \tr(X^{j})X^{k-1}({X^{-1}}^{*})^{m-l}\tr\left(X^{i+1-k} ({X^{-1}}^{*})^{l}\right)}^{=\displaystyle\Sigma_3}.
\end{array}$$
We compute
$$\begin{array}{rclr}
\Sigma_3&= &\displaystyle\sum_{\substack{1\leq i,j \leq n\\ i+j=n}}\displaystyle\sum_{\substack{1\leq k\leq i\\ 1\leq l \leq m}}i \tr(X^{j})X^{k-1}({X^{-1}}^{*})^{m-l}\tr\left(X^{i+1-k} ({X^{-1}}^{*})^{l}\right)&\\
& =&\displaystyle\sum_{\substack{1\leq i,j,k \leq n\\ i+j+k=n+1}}\displaystyle\sum_{1\leq l \leq m}(i+k-1) \tr(X^{j})\tr\left(X^{i} ({X^{-1}}^{*})^{l}\right)X^{k-1}({X^{-1}}^{*})^{m-l}&[i\leftarrow (i-k+1)]\\
& =&\displaystyle\sum_{\substack{1\leq i,j,k \leq n\\ i+j+k=n+1}}\displaystyle\sum_{1\leq l \leq m}i \tr(X^{j})\tr\left(X^{i} ({X^{-1}}^{*})^{l}\right)X^{k-1}({X^{-1}}^{*})^{m-l}&\\
&  &\hspace{1cm}+\displaystyle\sum_{\substack{1\leq i,j,k \leq n\\ i+j+k=n+1}}\displaystyle\sum_{1\leq l \leq m}(k-1) \tr(X^{j})\tr\left(X^{i} ({X^{-1}}^{*})^{l}\right)X^{k-1}({X^{-1}}^{*})^{m-l}.&
 \end{array}$$
 Finally, we conclude that $\Delta_{GL}\Delta_U^+(PQ^*)=\Delta_U^+\Delta_{GL}(PQ^*)$, using the following computation:
 $$\begin{array}{rclr}
\Sigma_3&  =&\displaystyle\sum_{\substack{1\leq i\leq n\\ 1\leq l \leq m}}\displaystyle\sum_{\substack{1\leq j,k \leq n+1-i\\ j+k=n+1-i}}j X^{i-1}({X^{-1}}^{*})^{m-l}\tr\left(X^{j} ({X^{-1}}^{*})^{l}\right)\tr\left(X^k\right)&\left[\begin{array}{l}i\leftarrow k\\j\leftarrow i\\ k\leftarrow i\end{array}\right]\\
&  &\hspace{1cm}+\displaystyle\sum_{\substack{1\leq i,k\leq n\\ 1\leq l \leq m}}\displaystyle\sum_{\substack{0\leq j \leq n-1\\ j+k=i-1}}j X^j\tr(X^{k})({X^{-1}}^{*})^{m-l}\tr\left(X^{n+1-i} ({X^{-1}}^{*})^{l}\right)&\left[\begin{array}{l} i\leftarrow n+1-i \\j\leftarrow k-1\\ k\leftarrow j\end{array}\right]\\
&=&\Sigma_2+\Sigma_1
 \end{array}$$
 where we remark that the case $j=0$ does not contribute in the last sum.
 
The other cases of $P, Q$ are similar, as well as the cases of the other couples of operators among the operators $\Delta_{GL}$, $\Delta_U^+$ and $\Delta_U^-$.

\subsubsection*{Step $2$}We prove now that the isomorphism $\mathcal{F}_t$ is indeed the free Hall transform. This can be done using the factorization of Section~\ref{fact}. More precisely, with the aim of working on polynomials, we will use the identity $P\left(G_t\right)=\left.P\right|_{G_t}\left(G_t\right)$ for all $P\in \mathbb{C}\{ X\}$.

Recall that $\frac{t}{2} \Delta_U$ coincides with $\D_{U_t}$ on $\mathbb{C}\{X\}$ (see Section~\ref{deltau}). For all $P\in \mathbb{C}[X]$, we have $(e^{\D_{U_t}}P)|_{G_t}=\mathcal{G}_t(P)$. Indeed, $P\mapsto (e^{\D_{U_t}}P)|_{G_t}$ and $P\mapsto \mathcal{G}_t(P)$ both respect the degree of $P$, and are isometries between $(\mathbb{C}[X],\|\cdot\|_{U_t})$ and $(\mathbb{C}[X],\|\cdot\|_{G_t})$ for the norms $\|\cdot\|_{U_t}:P\mapsto \|P(U_t)\|_2$ and $\|\cdot\|_{G_t}:P\mapsto \|P(G_t)\|_2$. The equality is then obtained by induction on the degree of $P$, provided that the leading coefficients are multiplied by the same factor. This is, in fact, the case: the coefficient of the leading term is multiplied by $e^{-nt/2}$ if the degree of the polynomial is $n$ (see the computation at the end of Section~\ref{utex} for $P\mapsto (e^{\D_{U_t}}P)|_{G_t}$ and the definition of $\mathcal{G}_t$ in Section~\ref{hallclas} for $P\mapsto \mathcal{G}_t(P)$).

This proves that, for all $P\in \mathbb{C}[X]$, we have $\mathcal{F}_t\left(P\left(U_t\right)\right) =e^{\frac{t}{2} \Delta_U}P(G_t)=e^{\D_{U_t}}P(G_t)=(e^{\D_{U_t}}P)|_{G_t}(G_t)=\mathcal{G}_t\Big(P\Big)(G_t)$, establishing the equality between $\mathcal{F}_t$ and the free Hall transform on the algebra generated by $U_t$. The same reasoning can be made on $\mathbb{C}\left[X^{-1}\right]$, and finally, $\mathcal{F}_t$ and the free Hall transform coincide on the algebra generated by $U_t$ and $U_t^{-1}$. Because they are isomorphisms, this equality extends to $L^2(U_t,\tau)$.
\end{proof}

\section{Random matrices}

In this last part of the paper, we will give two applications of the previous results. The first one, Theorem~\ref{Glim}, is the characterization of the distribution of a Brownian motion on $GL_N(\mathbb{C})$ at each fixed time in large-$N$ limit. The second one, Theorem~\ref{Flim}, establishes that the free Hall transform is the limit of the classical Hall transform for Laurent polynomial calculus.

Let us fix first some notation. Let $N$ be an integer and let $M_N(\mathbb{C})$ be the space of matrices of dimension $N$, with its canonical basis $(E_{a,b})_{1\leq a,b\leq N}$. If $M\in M_N(\mathbb{C})$, we denote
by $M^*$ the adjoint of $M$. Let us denote by $\Tr: M_N(\mathbb{C})\rightarrow \mathbb{R}$ the usual trace, and by $\tr: M_N(\mathbb{C})\rightarrow \mathbb{R}$ the normalized trace $\frac{1}{N}\Tr$, which makes $\left(M_N(\mathbb{C}),\tr\right)$ a non-commutative probability space.\label{granma}

\subsection{Brownian motion on $U(N)$}

Let $U(N)$ be the real Lie group of unitary matrices of dimension $N$:$$U(N)=\{U\in M_N(\mathbb{C}): U^*U=I_N\}.$$
Its Lie algebra $\frak{u}(N)$ consists of skew-hermitian matrices:
$$\frak{u}(N)=\{M\in M_N(\mathbb{C}): M^*+M=0\}.$$
We consider the following inner product on $\frak{u}(N)$:
$$(X,Y) \mapsto \left\langle X,Y\right\rangle_{\frak{u}(N)}=N \Tr (X^*Y).$$
This real-valued inner product is invariant under the adjoint action. Thus, it determines a bi-invariant Riemannian metric on $U(N)$ and a bi-invariant Laplace operator $\Delta_{U(N)}$.
Let $(X_1,\ldots,X_{N^2})$ be an orthonormal basis of $\frak{u}(N)$. Let us denote by $(\tilde{X}_1,\ldots,\tilde{X}_{N^2})$ the right-invariant vector fields on $U(N)$ which agree with $(X_1,\ldots,X_{N^2})$ at the identity. The Laplace operator $\Delta_{U(N)}$ is the second-order differential operator $\sum_{i=1}^{N^2}\tilde{X}_i^2$.

The (right) Brownian motion on $U(N)$ is a Markov process starting at the identity, and with generator $\frac{1}{2}\Delta_{U(N)}$.

\subsubsection{Computation of $\Delta_{U(N)} $}In this section, we will compute $\Delta_{U(N)} $ with the help of the space $\mathbb{C}\{X,X^{-1}\}$. It should be remarked that the operator $\Delta_{U(N)} $ has been calculated in other frameworks (see the recent~\cite{Driver2013}, and also~\cite{LEVY2008},~\cite{Rains1997} and~\cite{Sengupta2008}) using similar computations. We have at our disposal an operator $\Delta_U$ on $\mathbb{C}\{X,X^{-1}\}$, defined in Section~\ref{delta}. We will see in Lemma~\ref{delta} that the operator $\Delta_{U(N)} $ differs from $\Delta_U$ by an auxiliary operator $\frac{1}{N^2}\tilde{\Delta}_U$ when acting on the image of the $\mathbb{C}\{X,X^{-1}\}$-calculus.

Let us define this operator $\tilde{\Delta}_U$ on $\mathbb{C}\{X,X^{-1}\}$.
Let $k,l\in \mathbb{N}$, and $X_1,\ldots,X_{k+l}\in \{X,X^{-1}\}$. 
For all $1\leq i\leq k+l$, set $l (i)=0$, $r(i)=1$ and $\epsilon(i)=1$ if $X_i=X$, and $l(i)=1$, $r(i)=0$ and $\epsilon(i)=-1$ if $X_i=X^{-1}$. For all $1\leq i\leq k<  j\leq k+l$, the term $\overbracket{X_i\cdots X_j}$ refers to the product $X^{r(i)}_iX_{i+1}\cdots X_{j-1}X_j^{l(i)}$, that is to say the product $X_i\cdots X_j$  possibly excluding $X_i$ and/or $X_j$ according to $X_i$ and $X_j$, and the term $\overbracket{X_j\cdots X_i}$ refers to the product $X^{r(j)}_jX_{j+1}\cdots X_{k+l} X_1\cdots  X_{i-1}X_i^{l(i)}$, that is to say the product $X_j\cdots X_{k+l} X_1\cdots X_i$  possibly excluding $X_j$ and/or $X_i$. We set
\begin{multline*}\left\langle\nabla_U  X_1 \cdots X_k,\nabla_U X_{k+1}\cdots X_{k+l} \right\rangle\\
=- \displaystyle\sum_{1\leq i\leq k<  j\leq k+l} \epsilon(i)\epsilon(j)X_1 \cdots \widehat{\overbracket{X_i\cdots X_j}}\cdots X_{k+l}X_{k+1} \cdots \widehat{\overbracket{X_j\cdots X_i}}\cdots X_k,
\end{multline*}
where the hats mean that we have omitted the term $\overbracket{X_i\cdots X_j}$ in the product $X_1 \cdots X_{k+l}$, and the term $\overbracket{X_j\cdots X_i}$ in the product $X_{k+1} \cdots X_{k+l}X_{1} \cdots X_k$. Notice that $\left\langle\nabla_U  \cdot,\nabla_U \cdot \right\rangle$ is just a symbol for a bilinear map on $\mathbb{C}\left\langle X,X^{-1}\right\rangle$, but this notation is justified by Lemma~\ref{delta} (see the remark that follows Lemma~\ref{delta}).

Let $n\in \mathbb{N}$, and $P_0,\ldots,P_n$ be monomials of $\mathbb{C}\left\langle X,X^{-1}\right\rangle$. We set
\begin{equation} \label{structut}\begin{array}{rcl}
& &\hspace{-1.5cm} \tilde{\Delta}_U P_0 \tr P_1\cdots \tr P_n \\
&=& 2 \displaystyle\sum_{m=1}^n\left\langle\nabla_U  P_0 , \nabla_UP_m \right\rangle \tr P_1\cdots \widehat{\tr P_m}\cdots \tr P_n\\
& & \hspace{1cm}+\displaystyle\sum_{1\leq m,m'\leq n}P_0 \tr P_1\cdots \widehat{\tr P_m}\cdots \widehat{\tr P_{m'}}\cdots \tr P_n  \tr  \left\langle \nabla_U P_m , \nabla_UP_{m'} \right\rangle,
\end{array}\end{equation}
and we extend $\tilde{\Delta}_U$ to all $\mathbb{C}\{X,X^{-1}\}$ by linearity.

Let us denote by $U$ the identity function of $U(N)$.
\begin{lemma}
\label{delta}For all $P\in \mathbb{C}\{X,X^{-1}\}$, we have $$\displaystyle\Delta_{U(N)} \left(P(U)\right)=\left((\Delta_U+\frac{1}{N^2}\tilde{\Delta}_U)P\right)(U).$$
\end{lemma}
In particular, for $P,Q$ two monomials of $\mathbb{C}\left\langle X,X^{-1}\right\rangle$, we have $$\Delta_{U(N)} ((PQ)(U))=\left((\Delta_UP)\cdot_{\tr} Q+\frac{2}{N^2}\left\langle \nabla_U P , \nabla_UQ \right\rangle+P\cdot_{\tr} (\Delta_UQ)\right)(U).$$
\begin{proof}

Recall that $(X_1,\ldots,X_{N^2})$ is an orthonormal basis of $\frak{u}(N)$, that $(\tilde{X}_1,\ldots,\tilde{X}_{N^2})$ are the right-invariant vector fields on $U(N)$ which agree with $(X_1,\ldots,X_{N^2})$ at the identity, and that the Laplace operator $\Delta_{U(N)}$ is the second-order differential operator $\sum_{a=1}^{N^2}\tilde{X}_a^2$.

It suffices to prove that, for all $P$ and $Q$ monomials of $\mathbb{C}\left\langle X,X^{-1}\right\rangle$, we have
\begin{equation}\Delta_{U(N)} P(U) =\left(\Delta_{U} P\right)(U)\label{auxdelta} \end{equation}
and
\begin{equation}\displaystyle\sum_{a=1}^{N^2} \tilde{X_a}\left(P(U)\right)\tr \tilde{X_a}\left(Q(U)\right)=\frac{1}{N^2} \left\langle  \nabla_U P, \nabla_U Q \right\rangle (U).\label{auxnabla} \end{equation}
Indeed, let $P_0,\ldots,P_n$ be monomials of $\mathbb{C}\left\langle X,X^{-1}\right\rangle$. It follows from \eqref{auxdelta} and \eqref{auxnabla} that
$$\begin{array}{rcl}
& &\hspace{-1.5cm}\Delta_{U(N)} \left(P_0 \tr P_1\cdots \tr P_n\right)(U)\\
&=&\left(\Delta_{U(N)} P_0(U) \right) \tr P_1(U)\cdots \tr P_n(U)\\
& &\hspace{1cm}+\displaystyle\sum_{m=1}^n P_0(U) \tr P_1(U) \cdots \tr \left(\Delta_{U(N)} P_m(U) \right)\cdots \tr P_n(U)\\
&&\hspace{1cm}+ 2\displaystyle\sum_{m=1}^n\displaystyle\sum_{a=1}^{N^2}  \left( \tilde{X}_aP_0(U)  \right) \tr P_1(U)\cdots \tr \left(\tilde{X}_aP_m(U) \right)\cdots \tr P_n(U)\\
& &\hspace{1cm}+ \displaystyle\sum_{1\leq m,m'\leq n}\displaystyle\sum_{a=1}^{N^2} P_0 (U)\tr P_1(U)\cdots \tr \left(\tilde{X}_a P_m(U) \right)\cdots \tr \left(\tilde{X}_a P_m'(U) \right)\cdots \tr P_n(U).
\end{array}$$
Thanks to the structure equations \eqref{structDeltaU} and \eqref{structut} of $\Delta_U$ and $\tilde{\Delta}_U$, we deduce that
$$\begin{array}{rcl}
& &\hspace{-1.5cm}\Delta_{U(N)} \left(P_0 \tr P_1\cdots \tr P_n\right)(U)\\
&=&\left(\Delta_U P_0\right)(U)  \tr P_1(U)\cdots \tr P_n(U)\\
& &\hspace{1cm}+\displaystyle\sum_{m=1}^n P_0(U) \tr P_1(U) \cdots \tr \left(\Delta_U P_m\right)(U) \cdots \tr P_n(U)\\
&&\hspace{1cm}+\displaystyle\frac{2}{N^2} \displaystyle\sum_{m=1}^n\displaystyle\sum_{a=1}^{N^2}  \left\langle \nabla_U P_0, \nabla_U P_m \right\rangle (U) \tr P_1(U)\cdots \widehat{\tr P_m(U)} \cdots \tr P_n(U)\\
& &\hspace{1cm}+ \displaystyle\frac{1}{N^2} \displaystyle\sum_{1\leq m,m'\leq n}\displaystyle\sum_{a=1}^{N^2} P_0 (U)\tr P_1(U)\cdots \widehat{\tr P_m(U)}\cdots \widehat{\tr P_{m'}(U)}\cdots \tr P_n(U)\\
& &\hspace{5cm}\cdot \tr  \left\langle \nabla_U P_m ,\nabla_U P_{m'}\right\rangle(U) \\
&=&\Delta_U \left(P_0 \tr P_1\cdots \tr P_n \right)(U) +\frac{1}{N^2}\tilde{\Delta}_U  \left(P_0 \tr P_1\cdots \tr P_n \right)(U),
\end{array}$$
and Lemma~\ref{delta} follows by linearity. Thus, we need only prove \eqref{auxdelta} and \eqref{auxnabla}.

The end of the proof is devoted to proving \eqref{auxdelta} and \eqref{auxnabla}. For $1\leq a\leq N^2$, we have $\tilde{X}_a(U)=X_a U$ and $\tilde{X}_a(U^{-1})=\frac{\diff}{\diff t}|_{t=0} (e^{t X_a} U)^{-1}=\frac{\diff}{\diff t}|_{t=0} (U^{-1} e^{-t X_a}) =-U^{-1}X_a.$
Therefore, the expression of $\sum_{a=1}^{N^2}X_a\otimes X_a$ will be useful.
Let us fix the following orthonormal basis of $\mathfrak{u}(N)$: \begin{multline*}
\left\{X_1,\ldots,X_{N^2}\right\}\\
=\left\{\frac{i}{\sqrt{N}}E_{a,a},1\leq a\leq N\right\}\cup\left\{\frac{E_{a,b}-E_{b,a}}{\sqrt{2N}},1\leq a<b\leq N\right\}\cup\left\{i\frac{E_{a,b}+E_{b,a}}{\sqrt{2N}},1\leq a<b\leq N\right\}.
\end{multline*}
We have
$$\begin{array}{rcl}
N \displaystyle\sum_{a=1}^{N^2}X_a\otimes X_a &=&-\displaystyle\sum_{a=1}^N E_{a,a}\otimes E_{a,a}+\displaystyle\frac{1}{2}\displaystyle\sum_{1\leq a<b\leq N} \left(E_{a,b}-E_{b,a}\right)\otimes\left(E_{a,b}-E_{b,a}\right)\\
& &\hspace{1cm}-\displaystyle\frac{1}{2}\displaystyle\sum_{1\leq a<b\leq N}\left(E_{a,b}+E_{b,a}\right)\otimes \left(E_{a,b}+E_{b,a}\right)\\
&=&-\displaystyle\sum_{a=1}^N E_{a,a}\otimes E_{a,a}-\displaystyle\sum_{1\leq a<b\leq N} E_{b,a}\otimes E_{a,b}C-\displaystyle\sum_{1\leq a<b\leq N} E_{a,b}\otimes E_{b,a}\\
&=&-\displaystyle\sum_{1\leq a,b\leq N}E_{a,b} \otimes E_{b,a}.
\end{array}$$
From which we deduce that
\begin{equation*}\displaystyle\sum_{a=1}^{N^2}\left(\tilde{X}_a(U)\otimes \tilde{X}_a(U) \right)=-\displaystyle\frac{1}{N}\displaystyle\sum_{1\leq a,b\leq N}E_{a,b}U \otimes E_{b,a}U,\end{equation*}
\begin{equation*}\displaystyle\sum_{a=1}^{N^2}\left(\tilde{X}_a(U)\otimes \tilde{X}_a(U^{-1}) \right)=\displaystyle\frac{1}{N}\displaystyle\sum_{1\leq a,b\leq N}E_{a,b}U \otimes UE_{b,a},\end{equation*}
and \begin{equation*}\displaystyle\sum_{a=1}^{N^2}\left(\tilde{X}_a(U^{-1})\otimes \tilde{X}_a(U^{-1}) \right)=-\displaystyle\frac{1}{N}\displaystyle\sum_{1\leq a,b\leq N}U^{-1}E_{a,b} \otimes U^{-1}E_{b,a}.\end{equation*}
We can now compute \eqref{auxdelta} and \eqref{auxnabla}. Let $k,l\in \mathbb{N}$, and $X_1,\ldots,X_{k+l}\in \{X,X^{-1}\}$. 
For all $1\leq i\leq k+l$, set $\epsilon(i)=1$ if $X_i=X$, and $\epsilon(i)=-1$ if $X_i=X^{-1}$.
Thus, we have
$$\begin{array}{rcl}
\Delta_{U(N)} \left(P(U)\right)
&=&\displaystyle\sum_{a=1}^{N^2}\tilde{X}_a^2\left( U^{\epsilon (1)}\cdots U^{\epsilon (k)}\right)\\
&=&\displaystyle\sum_{a=1}^{N^2} \displaystyle\sum_{1\leq i \leq k} U^{\epsilon (1)}\cdots U^{\epsilon (i-1)} \tilde{X}_a^2\left( U^{\epsilon (i)} \right) U^{\epsilon (i+1)}\cdots U^{\epsilon (k)}\\
& &\hspace{1cm}+\displaystyle\sum_{a=1}^{N^2}2\displaystyle\sum_{1\leq i< j\leq k}  U^{\epsilon (1)}\cdots  \tilde{X}_a\left( U^{\epsilon (i)} \right)\cdots  \tilde{X}_a\left( U^{\epsilon (j)} \right)\cdots U^{\epsilon (k)} \\
&=&- k \displaystyle\frac{N}{N}U^{\epsilon (1)}\cdots U^{\epsilon (k)}\\
& &\hspace{1cm}-2\displaystyle\sum_{1\leq i< j\leq k} \epsilon(i)\epsilon(j)\left(X_1\cdots  \widehat{\overbracket{X_i\cdots X_j}}\cdots X_k \tr \left(\overbracket{X_i\cdots X_j}\right)\right)(U)\\
&=&\left(\Delta_U P\right)
(U)\end{array}$$
and
$$\begin{array}{rcl}
& &\hspace{-1.2cm}\displaystyle\sum_{a=1}^{N^2}\tilde{X_a}\left(P(U)\right)\tr \tilde{X_a}\left(Q(U)\right)\\
&=&\displaystyle\sum_{a=1}^{N^2} \tilde{X_a}\left(U^{\epsilon (1)}\cdots U^{\epsilon (k)}\right)\tr \tilde{X_a}\left(U^{\epsilon (k+1)}\cdots U^{\epsilon (k+l)}\right)\\
&=&\displaystyle\frac{1}{N}   \displaystyle\sum_{1\leq i\leq k<  j\leq k+l}\displaystyle\sum_{a=1}^{N^2}\left(U^{\epsilon (1)}\cdots \tilde{X_a}\left(U^{\epsilon (i)}\right) \cdots U^{\epsilon (k)}\right)\Tr \left( U ^{\epsilon (k+1)}\cdots  \tilde{X_a}\left(U^{\epsilon (j)}\right) \cdots U^{\epsilon (k+l)}\right)\\
&=&-\displaystyle\frac{1}{N^2} \displaystyle\sum_{1\leq i\leq k<  j\leq k+l} \epsilon(i)\epsilon(j)\left(X_1 \cdots \widehat{\overbracket{X_i\cdots X_j}}\cdots X_{k+l}X _{k+1}\cdots \widehat{\overbracket{X_j\cdots X_i}}\cdots X_k\right)(U)\\
&=&\displaystyle\frac{1}{N^2} \left\langle \nabla_U P, \nabla_UQ \right\rangle (U)
\end{array}$$
where we have used that,  for all $A,B,C,D \in M_N(\mathbb{C})$, we have
\begin{equation} \displaystyle\sum_{1\leq a,b\leq N}A E_{a,b}B E_{b,a}C= \Tr(B)AC\text{ and }\displaystyle\sum_{1\leq a,b\leq N}A E_{a,b}B \Tr(C E_{b,a}D)= ADCB.\label{magic} \qedhere \end{equation}
\end{proof}
\begin{corollary}
Let $\left(U_t^{(N)}\right)_{t\geq 0}$ be a Brownian motion on $U(N)$.

Let $t\geq 0$ and $P\in \mathbb{C}\{X,X^{-1}\}$. We have
$$\mathbb{E}\left[P\left(U_t^{(N)}\right)\right]=\left(e^{\frac{t}{2} (\Delta_U+\frac{1}{N^2}\tilde{\Delta}_U)}P\right)(1).$$
\end{corollary}
\begin{proof}For all $t\geq 0$ and $P\in \mathbb{C}\{X,X^{-1}\}$, thanks to Lemma~\ref{delta}, we have
$$\frac{\diff}{\diff t}\mathbb{E}\left[P\left(U_t^{(N)}\right)\right]= \mathbb{E}\left[\left(\frac{1}{2}\Delta_{U(N)}\Big(P(U)\Big)\right)\left(U_t^{(N)}\right)\right]=\mathbb{E}\left[\left(\frac{1}{2}(\Delta_U+\frac{1}{N^2}\tilde{\Delta}_U)P\right)\left(U_t^{(N)}\right)\right].$$
We conclude with Lemma~\ref{eqdiff}, since the maps $(t\mapsto \mathbb{E}\left[P\left(U_t^{(N)}\right)\right])_{t\geq 0}$ are linear.
\end{proof}

\subsection{Brownian motion on $GL_N(\mathbb{C})$}

Let $GL_N(\mathbb{C})$ be the Lie group of invertible $N\times N$ matrices. We regard $GL_N(\mathbb{C})$ as a real Lie group, with its Lie algebra $\frak{gl}_N(\mathbb{C})=M_N(\mathbb{C})$. We endow $\frak{gl}_N(\mathbb{C})$ with the following real-valued inner product:
$$(X,Y) \mapsto \left\langle X,Y\right\rangle_{\frak{gl}_N(\mathbb{C})}=N \Re \Tr (X^*Y).$$
Unfortunately, this inner product is not invariant under the adjoint action: the left-invariant Riemannian metric which it determines does not coincide with the right-invariant Riemannian metric.
Let us choose the right-invariant Riemannian metric. Let $\Delta_{GL_N(\mathbb{C})}$ be the Laplace operator corresponding to the right-invariant metric.
Let $(Z_1,\ldots,Z_{2 N^2})$ be an orthonormal basis of $M_N(\mathbb{C})$. Let us denote by $(\tilde{Z}_1,\ldots,\tilde{Z}_{2 N^2})$ the right-invariant vector fields on $GL_N(\mathbb{C})$ which agree with $(Z_1,\ldots,Z_{2 N^2})$ at the identity. The Laplace operator $\Delta_{GL_N(\mathbb{C})}$ is the second-order differential operator $\sum_{i=1}^{2N^2}\tilde{Z}_i^2$.

The (right) Brownian motion on $GL_N(\mathbb{C})$ is a Markov process starting at the identity, and with generator $\frac{1}{4}\Delta_{GL_N(\mathbb{C})}$. As mentioned in the introduction, we have taken a definition of the Brownian motion on $GL_N(\mathbb{C})$ which differs by a factor of $2$ from the usual definition for reasons of convenience.

\subsubsection{Computation of $\Delta_{GL_N(\mathbb{C})} $}Similarly to $\Delta_{U(N)} $, we will see in Lemma~\ref{deltagl} that the operator $\Delta_{GL_N(\mathbb{C})}$ differs from $\Delta_{GL}$ defined in Section~\ref{deltag} by an auxiliary operator $\frac{1}{N^2}\tilde{\Delta}_{GL}$ when acting on the image of the $\mathbb{C}\{X,X^*,X^{-1},{X^*}^{-1}\}$-calculus.

Let us define this operator $\tilde{\Delta}_{GL}$ on $\mathbb{C}\{X,X^*,X^{-1},{X^*}^{-1}\}$.
Let $k,l\in \mathbb{N}$, and $X_1,\ldots,X_{k+l}\in \{X,X^*,X^{-1},{X^*}^{-1}\}$.
For all $1\leq i\leq k+l$, set $l (i)=0$, $r(i)=1$ if $X_i=X$ or $X_i={X^*}^{-1}$, and $l(i)=1$, $r(i)=0$ if $X_i=X^*$ or $X_i=X^{-1}$.
For all $1\leq i\leq k<  j\leq k+l$, the term $\overbracket{X_i\cdots X_j}$ refers to the product $X^{r(i)}_iX_{i+1}\cdots X_{j-1}X_j^{l(i)}$, that is to say the product $X_i\cdots X_j$  possibly excluding $X_i$ and/or $X_j$, and the term $\overbracket{X_j\cdots X_i}$ refers to the product $X^{r(j)}_jX_{j+1}\cdots X_{k+l} X_1\cdots  X_{i-1}X_i^{l(i)}$, that is to say the product $X_j\cdots X_{k+l} X_1\cdots X_i$  possibly excluding $X_j$ and/or $X_i$.
For all $1\leq i\leq k+l$, set $\epsilon(i)=1$ if $X_i=X$ or $X_i={X^*}$, and $\epsilon(i)=-1$ if $X_i={X^*}^{-1}$ or $X_i=X^{-1}$. For all $1\leq i\leq k< j\leq k+l$, set $\delta(i,j)=0$ if $X_i,X_j\in \{X,X^{-1}\}$ or if $X_i,X_j\in \{X^*,{X^*}^{-1}\}$, and $\delta(i,j)=1$ otherwise. We set
\begin{eqnarray*}& &\hspace{-1.5cm}\left\langle \nabla_{GL} X_1\cdots X_k,\nabla_{GL} X_{k+1}\cdots X_{k+l} \right\rangle\\
&=&2\displaystyle\sum_{1\leq i\leq k< j\leq k+l}\delta(i,j)\epsilon(i)\epsilon(j)X_1 \cdots \widehat{\overbracket{X_i\cdots X_j}}\cdots X_{k+l}X _{k+1}\cdots \widehat{\overbracket{X_j\cdots X_i}}\cdots X_k.
\end{eqnarray*}
Here again, it is important to note that $\langle \nabla_{GL} \cdot,\nabla_{GL} \cdot \rangle$ is just a notation for a bilinear map on $\mathbb{C}\langle X,X^*,X^{-1},{X^*}^{-1}\rangle$ which may be justified thanks to Lemma~\ref{deltagl} (see the remark that follows Lemma~\ref{deltagl}).

Let $n\in \mathbb{N}$, and $P_0,\ldots,P_n$ be monomials of $\mathbb{C}\langle X,X^*,X^{-1},{X^*}^{-1}\rangle$. We set
\begin{equation}\label{structutgl}\begin{array}{rcl}
& &\hspace{-1.5cm} \tilde{\Delta}_{GL} (P_0 \tr P_1\cdots \tr P_n ) \\
&=& 2 \displaystyle\sum_{m=1}^n\left\langle \nabla_{GL} P_0 , P_m \right\rangle\nabla_{GL} \tr P_1\cdots \widehat{\tr P_m}\cdots \tr P_n\\
& &\hspace{1cm}+ \displaystyle\sum_{1\leq m,m'\leq n}P_0 \tr P_1\cdots \widehat{\tr P_m}\cdots \widehat{\tr P_{m'}}\cdots \tr P_n  \tr  \left\langle \nabla_{GL} P_m , \nabla_{GL}P_{m'} \right\rangle.
\end{array}\end{equation}
We extend $\tilde{\Delta}_{GL}$ to all $\mathbb{C}\{X,X^*,X^{-1},{X^*}^{-1}\}$ by linearity.

Let us denote by $G$ the identity function of $GL_N(\mathbb{C})$. 
\begin{lemma}
\label{deltagl}For all $P\in \mathbb{C}\{X,X^*,X^{-1},{X^*}^{-1}\}$, we have
$$\Delta_{GL_N(\mathbb{C})}  \left(P(G)\right)=\left((\Delta_{GL}+\frac{1}{N^2}\tilde{\Delta}_{GL})P\right)(G).$$
\end{lemma}
In particular, for $P,Q$ two monomials of $\mathbb{C}\left\langle X,X^*,X^{-1},{X^*}^{-1}\right\rangle$, we have $$\Delta_{GL_N(\mathbb{C})} ((PQ)(U))=\left((\Delta_{GL}P)\cdot_{\tr} Q+\frac{2}{N^2}\left\langle \nabla_{GL} P , \nabla_{GL}Q \right\rangle+P\cdot_{\tr} (\Delta_{GL}Q)\right)(U).$$
\begin{proof}The proof follows the demonstration of Lemma~\ref{delta}.

Recall that $(Z_1,\ldots,Z_{2N^2})$ is an orthonormal basis of $M_N(\mathbb{C})$, that $(\tilde{Z}_1,\ldots,\tilde{Z}_{2N^2})$ are the right-invariant vector fields on $U(N)$ which agree with $(Z_1,\ldots,Z_{2N^2})$ at the identity, and that the Laplace operator $\Delta_{GL_N(\mathbb{C})}$ is the second-order differential operator $\sum_{a=1}^{2N^2}\tilde{Z}_a^2$.

It suffices to prove that, for all $P$ and $Q$ monomials of $\mathbb{C}\langle X,X^*,X^{-1},{X^*}^{-1}\rangle$, we have
\begin{equation}\Delta_{GL_N(\mathbb{C})} P(G) =\left(\Delta_{GL} P\right)(G)\label{auxdeltagl}\end{equation}
and
\begin{equation}\displaystyle\sum_{a=1}^{2N^2} \tilde{Z_a}\left(P(G)\right)\tr \tilde{Z_a}\left(Q(G)\right)=\frac{1}{N^2} \left\langle\nabla_{GL}  P, \nabla_{GL}Q \right\rangle (G).\label{auxnablagl} \end{equation}
Indeed, let $P_0,\ldots,P_n$ be monomials of $\mathbb{C}\langle X,X^*,X^{-1},{X^*}^{-1}\rangle$. It follows from \eqref{auxdeltagl} and \eqref{auxnablagl} that
$$\begin{array}{rcl}
& &\hspace{-1.5cm}\Delta_{GL_N(\mathbb{C})} \left(P_0 \tr P_1\cdots \tr P_n\right)(G)\\
&=&\left(\Delta_{GL_N(\mathbb{C})} P_0(G) \right) \tr P_1(G)\cdots \tr P_n(G)\\
& &\hspace{1cm}+\displaystyle\sum_{m=1}^n P_0(G) \tr P_1(G) \cdots \tr \left(\Delta_{GL_N(\mathbb{C})} P_m(G) \right)\cdots \tr P_n(G)\\
&&\hspace{1cm}+ 2\displaystyle\sum_{m=1}^n\displaystyle\sum_{a=1}^{2N^2}  \left( \tilde{Z}_aP_0(G)  \right) \tr P_1(G)\cdots \tr \left(\tilde{Z}_aP_m(G) \right)\cdots \tr P_n(G)\\
& &\hspace{1cm}+ \displaystyle\sum_{1\leq m,m'\leq n}\displaystyle\sum_{a=1}^{2N^2} P_0 (G)\tr P_1(G)\cdots \tr \left(\tilde{Z}_a P_m(G) \right)\cdots \tr \left(\tilde{Z}_a P_m'(G) \right)\cdots \tr P_n(G).
\end{array}$$
Thanks to the structure equations \eqref{structugl} and \eqref{structutgl} of $\Delta_{GL}$ and $\tilde{\Delta}_{GL}$, we deduce that
$$\begin{array}{rcl}
& &\hspace{-1.5cm}\Delta_{GL_N(\mathbb{C})} \left(P_0 \tr P_1\cdots \tr P_n\right)(G)\\
&=&\left(\Delta_{GL} P_0\right)(G)  \tr P_1(G)\cdots \tr P_n(G)\\
& &\hspace{1cm}+\displaystyle\sum_{m=1}^n P_0(G) \tr P_1(G) \cdots \tr \left(\Delta_{GL} P_m\right)(G) \cdots \tr P_n(G)\\
&&\hspace{1cm}+\displaystyle\frac{2}{N^2} \displaystyle\sum_{m=1}^n\displaystyle\sum_{a=1}^{2N^2}  \left\langle \nabla_{GL} P_0,\nabla_{GL} P_m \right\rangle (G) \tr P_1(G)\cdots \widehat{\tr P_m(G)} \cdots \tr P_n(G)\\
& &\hspace{1cm}+ \displaystyle\frac{1}{N^2} \displaystyle\sum_{1\leq m,m'\leq n}\displaystyle\sum_{a=1}^{2N^2} P_0 (G)\tr P_1(G)\cdots \widehat{\tr P_m(G)}\cdots \widehat{\tr P_{m'}(G)}\cdots \tr P_n(G) \\
& &\hspace{5cm}\cdot\tr  \left\langle\nabla_{GL}  P_m ,\nabla_{GL} P_{m'}\right\rangle(G) \\
&=&\Delta_{GL} \left(P_0 \tr P_1\cdots \tr P_n \right)(G) +\displaystyle\frac{1}{N^2}\tilde{\Delta}_{GL}  \left(P_0 \tr P_1\cdots \tr P_n \right)(G),
\end{array}$$
and Lemma~\ref{deltagl} follows by linearity. Thus, we need only prove \eqref{auxdeltagl} and \eqref{auxnablagl}.
For $1\leq a\leq 2N^2$, we have $\tilde{Z}_a\left(G\right)=Z_a G$,
$$
\tilde{Z}_a\left(G^{-1}\right)=\left.\displaystyle\frac{\diff}{\diff t}\right|_{t=0} \left(e^{t Z_a} G\right)^{-1}=\left.\displaystyle\frac{\diff}{\diff t}\right|_{t=0} \left(G^{-1} e^{-t Z_a}\right) =-G^{-1}Z_a,$$
and similarly $\tilde{Z}_a\left(G^*\right)=G^*Z_a^*,$
and
$\tilde{Z}_a\left({G^*}^{-1}\right)=-Z_a^*G^{-1}$.
Therefore, the expressions of $\sum_{a=1}^{2N^2}Z_a\otimes Z_a$, $\sum_{a=1}^{2N^2}Z_a\otimes Z_a^*$ and $\sum_{a=1}^{2N^2}Z_a^*\otimes Z_a^*$ will be useful.
Let us fix the basis $$\left\{Z_1,\ldots,Z_{2N^2}\right\}=\{\frac{1}{\sqrt{N}}E_{a,b},\frac{i}{\sqrt{N}}E{a,b}:1\leq a,b\leq N\}.$$
We have
%
\begin{equation*}\displaystyle\sum_{a=1}^{N^2}Z_a\otimes Z_a =\displaystyle\sum_{1\leq a,b\leq N} \frac{1}{N}E_{a,b}\otimes E_{a,b}+\displaystyle\sum_{1\leq a,b\leq N} \frac{-1}{N}E_{a,b}\otimes E_{a,b}=0,\end{equation*}
\begin{equation*}\displaystyle\sum_{a=1}^{N^2}Z_a\otimes Z_a^* =\displaystyle\sum_{1\leq a,b\leq N} \frac{1}{N}E_{a,b}\otimes E_{b,a}+\displaystyle\sum_{1\leq a,b\leq N} \frac{1}{N}E_{a,b}\otimes E_{b,a}=\displaystyle\frac{2}{N}\displaystyle\sum_{1\leq a,b\leq N}E_{a,b}\otimes E_{b,a},\end{equation*}
and\begin{equation*}\displaystyle\sum_{a=1}^{N^2}Z_a^*\otimes Z_a^* =\displaystyle\sum_{1\leq a,b\leq N} \frac{1}{N}E_{b,a}\otimes E_{b,a}+\displaystyle\sum_{1\leq a,b\leq N} \frac{-1}{N}E_{b,a}\otimes E_{b,a}=0.\end{equation*}
From this, we deduce that $\sum_{a=1}^{2N^2}\tilde{Z}^2_a(G^{\varepsilon})=0$ for all $\varepsilon \in \{1,-1,*,*-1\}$. Moreover,
%
\begin{equation*}\displaystyle\sum_{a=1}^{2N^2}\left(\tilde{Z}_a(G)\otimes \tilde{Z}_a(G^{-1}) \right)=\displaystyle\sum_{a=1}^{N^2}\left(\tilde{Z}_a(G^*)\otimes \tilde{Z}_a({G^*}^{-1}) \right)=0,\end{equation*}
\begin{equation*}\displaystyle\sum_{a=1}^{2N^2}\left(\tilde{Z}_a(G)\otimes \tilde{Z}_a(G^*) \right)=\displaystyle\frac{2}{N}\displaystyle\sum_{1\leq a,b\leq N}E_{a,b}G \otimes G^*E_{b,a},\end{equation*}
\begin{equation*}\displaystyle\sum_{a=1}^{2N^2}\left(\tilde{Z}_a(G)\otimes \tilde{Z}_a({G^*}^{-1}) \right)=-\displaystyle\frac{2}{N}\displaystyle\sum_{1\leq a,b\leq N}E_{a,b}G \otimes E_{b,a}{G^*}^{-1},\end{equation*}
\begin{equation*}\displaystyle\sum_{a=1}^{2N^2}\left(\tilde{Z}_a(G^{-1})\otimes \tilde{Z}_a(G^*) \right)=-\displaystyle\frac{2}{N}\displaystyle\sum_{1\leq a,b\leq N}G^{-1}E_{a,b} \otimes G^*E_{b,a},\end{equation*}
and
\begin{equation*}\displaystyle\sum_{a=1}^{2N^2}\left(\tilde{Z}_a(G^{-1})\otimes \tilde{Z}_a({G^*}^{-1}) \right)=\displaystyle\frac{2}{N}\displaystyle\sum_{1\leq a,b\leq N}G^{-1}E_{a,b} \otimes E_{b,a}{G^*}^{-1}.\end{equation*}
Let $k,l\in \mathbb{N}$, and $X_1,\ldots,X_{k+l}\in \{X,X^*,X^{-1},{X^*}^{-1}\}$. For all $1\leq i\leq k+l$, set $\epsilon(i)=1$ if $X_i=X$ or $X_i={X^*}$, and $\epsilon(i)=-1$ if $X_i={X^*}^{-1}$ or $X_i=X^{-1}$.
For all $1\leq i\leq k< j\leq k+l$, set $\delta(i,j)=0$ if $X_i,X_j\in \{X,X^{-1}\}$ or if $X_i,X_j\in \{X^*,{X^*}^{-1}\}$, and $\delta(i,j)=1$ otherwise. For all $1\leq i\leq k+l$, set $G_i=X_i(G)$. We can now compute
$$\begin{array}{rcl}
\Delta_{GL_N(\mathbb{C})} \left(P(G)\right)
&=&\displaystyle\sum_{a=1}^{2N^2}\tilde{Z}_a^2\left( G_1\cdots G_k\right)\\
&=&\displaystyle\sum_{a=1}^{2N^2} 2 \displaystyle\sum_{1\leq i< j\leq k} G_1\cdots \tilde{Z}_a\left(G_i\right) \cdots \tilde{Z}_a\left(G_j\right) \cdots G_k  \\
&=&- 4\displaystyle\sum_{1\leq i< j\leq k} \delta(i,j)\epsilon(i)\epsilon(j)
\cdot \left( X_1\cdots  \widehat{\overbracket{X_i\cdots X_j}}\cdots X_k \tr \left(\overbracket{X_i\cdots X_j}\right)\right)(G)\\
&=&\left(\Delta_{GL} P\right)
(G)\end{array}$$
and
$$\begin{array}{rcl}
& &\hspace{-1cm}\displaystyle\sum_{a=1}^{2N^2}\tilde{Z_a}\left(P(G)\right)\tr \tilde{Z_a}\left(Q(G)\right)\\
&=&\displaystyle\sum_{a=1}^{2N^2} \tilde{Z_a}\left(G_1\cdots G_k\right)\tr \tilde{Z_a}\left(G_{k+1}\cdots G_{k+l}\right)\\
&=&\displaystyle\frac{1}{N}   \displaystyle\sum_{1\leq i\leq k<  j\leq k+l}\displaystyle\sum_{a=1}^{N^2}\left(G_1\cdots \tilde{Z_a}\left(G_i\right) \cdots G_k\right)\Tr \left( G_{k+1}\cdots  \tilde{Z_a}\left(G_j\right) \cdots G_{k+l}\right)\\
&=&\displaystyle\frac{2}{N^2}\displaystyle\sum_{1\leq i\leq k< j\leq k+l}\delta(i,j)\epsilon(i)\epsilon(j)\left(X_1 \cdots \widehat{\overbracket{X_i\cdots X_j}}\cdots X_{k+l}X _{k+1}\cdots \widehat{\overbracket{X_j\cdots X_i}}\cdots X_k\right)(G)\\
&=&\displaystyle\frac{1}{N^2} \left\langle \nabla_{GL} P, \nabla_{GL}Q \right\rangle (G),
\end{array}$$
where we have used~\eqref{magic}.
\end{proof}
\begin{corollary}
Let $\left(G_t^{(N)}\right)_{t\geq 0}$ be a Brownian motion on $GL_N(\mathbb{C})$.

Let $t\geq 0$ and $P\in \mathbb{C}\{X,X^*,X^{-1},{X^*}^{-1}\}$. We have\label{GLN}
$$\mathbb{E}\left[P\left(G_t^{(N)}\right)\right]=\left(e^{\frac{t}{4} (\Delta_{GL}+\frac{1}{N^2}\tilde{\Delta}_{GL})}P\right)(1).$$
\end{corollary}
\begin{proof}For all $t\geq 0$ and $P\in\mathbb{C}\{X,X^*,X^{-1},{X^*}^{-1}\}$, thanks to Lemma~\ref{deltagl}, we have
$$\frac{\diff}{\diff t}\mathbb{E}\left[P\left(G_t^{(N)}\right)\right]= \mathbb{E}\left[\left(\frac{1}{4}\Delta_{GL_N(\mathbb{C})}\Big(P(G)\Big)\right)\left(G_t^{(N)}\right)\right]=\mathbb{E}\left[\left(\frac{1}{4}(\Delta_{GL}+\frac{1}{N^2}\tilde{\Delta}_{GL})P\right)\left(G_t^{(N)}\right)\right].$$
We conclude with Lemma~\ref{eqdiff}, since the maps $(t\mapsto \mathbb{E}\left[P\left(G_t^{(N)}\right)\right])_{t\geq 0}$ are linear.
\end{proof}

\subsection{Hall transform}In this section, we introduce the Hall transform, and study it in the particular case of the unitary group $U(N)$. We refer to~\cite{Driver1995} or~\cite{Hall1994} for more details. The major difference is that Hall and Driver work with left-invariant metrics and Laplace operators. Because we chose to work with right Brownian motion, we will only consider right-invariant metrics and Laplace operators. The necessary modifications are minor.\label{hallclas}
\subsubsection{Definition}
Let $K$ be a compact connected Lie group, with Lie algebra $\frak{k}$. Let $\langle\cdot,\cdot\rangle_{\frak{k}}$ denote a fixed $\Ad(K)$-invariant inner product on $\frak{k}$. It determines a bi-invariant Riemannian metric on $K$ and a bi-invariant Laplace operator $\Delta_{K}$. For all $t\geq 0$, let $\diff \rho_t$ be the heat kernel measure on $K$ at time $t$; that is to say, the probability measure on $K$ which corresponds to the law at time $t$ of a Markov process starting at the identity and with generator $\frac{1}{2}\Delta_{K}$. For all $t\geq 0$, the operator $e^{\frac{t}{2}\Delta_{K}}$ over $L^2(K,\diff \rho_t)$ is given for all $f\in L^2(K,\diff \rho_t)$ by $e^{\frac{t}{2}\Delta_{K}}f:x\mapsto \int_K f(yx)\diff \rho_t(y)$.

Let $G$ be the complexification of $K$, with Lie algebra $\frak{g}=\frak{k}\oplus i\frak{k}$. Let us endow $\frak{k}$ with the (real-valued) inner product given by $\langle K_1+iK_2,K_3+iK_4\rangle_{\frak{g}}=\langle K_1,K_3\rangle_{\frak{k}}+\langle K_2,K_4\rangle_{\frak{k}}$ for all $K_1,K_2,K_3,K_4\in \frak{k}$: it determines a right-invariant Riemannian metric on $G$ and a right-invariant Laplace operator $\Delta_{G}$.  For all $t\geq 0$, let $\diff \mu_t$ be the "half heat kernel" measure on $G$ at time $t$; that is to say, the probability measure on $G$ which corresponds to the law at time $t$ of a Markov process starting at the identity and with generator $\frac{1}{4}\Delta_{G}$.


For all $t> 0$, we denote by $L^2_{\hol}(G,\diff \mu_t)$ the Hilbert space of holomorphic function in $L^2(G,\diff \mu_t)$. The fact that $L^2_{\hol}(G,\diff \mu_t)$ is a Hilbert space is not trivial. It is a part of Hall's theorem which may be stated as follows (see Theorem~1 of~\cite{Hall1994} or Theorem~1.16 of~\cite{Driver1995}).
\begin{theo*}[Hall~\cite{Hall1994}]Let $t> 0$. For all $f\in L^2(K,\diff \rho_t)$, the function $e^{\frac{t}{2}\Delta_{K}}f$ has an analytic continuation to a holomorphic function on $G$, also denoted by $e^{\frac{t}{2}\Delta_{K}}f$. Moreover, $e^{\frac{t}{2}\Delta_{K}}f\in L^2_{\hol}(G,\diff \mu_t)$ and the linear map $B_t:f\mapsto e^{\frac{t}{2}\Delta_{K}}f$
is an isomorphism of Hilbert spaces between $L^2(K,\diff \rho_t)$ and $L^2_{\hol}(G,\diff \mu_t)$.
\end{theo*}

\subsubsection{The Hall transform on $U(N)$}In this section, we take for $K$ the Lie group $U(N)$, whose complexification is $GL_N(\mathbb{C})$. Let $(U_t^{(N)})_{t\geq 0}$ be a Brownian motion on $U(N)$ (with generator $\frac{1}{2}\Delta_{U(N)}$), and $(G_t^{(N)})_{t\geq 0}$ be a Brownian motion on $GL_N(\mathbb{C})$ (with generator $\frac{1}{4}\Delta_{GL_N(\mathbb{C})}$).\label{rvfr} With the notation of the previous section, for all $t\geq 0$, $\diff \rho_t$ is the law of $U_t^{(N)}$ and $\diff \mu_t$ is the law of $G_t^{(N)}$. Hall's theorem gives us an isomorphism of Hilbert spaces between $L^2(U(N),\diff \rho_t)$ and $L^2_{\hol}(GL_N(\mathbb{C}),\diff \mu_t)$. We remind of the probabilistic formulation of Hall's theorem in the introduction: it expresses the Hall transform $B_t$ as an isomorphism of Hilbert spaces between the spaces of random variables $L^2(U^{(N)}_t)$ and $L^2_{\hol}(G^{(N)}_t)$ as defined in the introduction. Let us explicit the identification made between the spaces of functions $L^2(U(N),\diff \rho_t)$ and $L^2_{\hol}(GL_N(\mathbb{C}),\diff \mu_t)$ and the spaces of random variables $L^2(U^{(N)}_t)$ and $L^2_{\hol}(G^{(N)}_t)$.

Let $t> 0$. For all complex Borel function $f$, we have $\mathbb{E}[|f(U^{(N)}_t)|^2]=\int_{U(N)}|f|^2\diff \rho_t$ in $[0,+\infty]$. Thus, if $f$ and $g$ are two complex Borel function such that $f(U_t^N)=g(U_t^N)\ a.s.$, then ${\int_{U(N)}|f-g|^2\diff \rho_t=}0$ and therefore $f=g\ a.s$. Moreover, the map $f\mapsto f(U_t^{(N)})$ is an isomorphism between $L^2(U(N),\diff \rho_t)$ and $L^2(U^{(N)}_t)$. Similarly, $F\mapsto F(G_t^N)$ is an isomorphism between $L^2_{\hol}(GL_N(\mathbb{C}),\diff \mu_t)$ and $L^2_{\hol}(G^{(N)}_t)$.
Identifying $L^2(U(N),\diff \rho_t)$ with $L^2(U^{(N)}_t)$ and $L^2_{\hol}(GL_N(\mathbb{C}),\diff \mu_t)$ with $L^2_{\hol}(G^{(N)}_t)$, we deduce the version of Hall's theorem as formulated in the introduction: the Hall transform $B_t$ is an isomorphism of Hilbert spaces between $L^2(U^{(N)}_t)$ and $L^2_{\hol}(G^{(N)}_t)$, and for all Borel function $f$ on $U(N)$ such that $f(U^{(N)}_t)\in L^2(U^{(N)}_t)$, we have $$B_t\Big(f(U^{(N)}_t)\Big)=\left(e^{\frac{t}{2}\Delta_{U(N)}}f\right)\left(G^{(N)}_t\right).$$

\subsubsection{Boosting}
For all $N\in \mathbb{N}^*$, we endow $M_N(\mathbb{C})$ with the inner product $\langle X,Y\rangle_{M_N(\mathbb{C})}=\frac{1}{N} \Tr (X^*Y)$. For all $t> 0$ and $N\in \mathbb{N}^*$, we denote by $B^{(N)}_t$ the Hilbert space isomorphism $B_t\otimes \Id_{M_N(\mathbb{C})}$ from $L^2(U^{(N)}_t)\otimes M_N(\mathbb{C})$ into $L^2_{\hol}(G^{(N)}_t)\otimes M_N(\mathbb{C})$.
\begin{proposition}
Let $t> 0$. For all $P\in \mathbb{C}\{X,X^{-1}\}$, we have\label{btpun}
$$B^{(N)}_t\left(P\left(U^{(N)}_t\right)\right)=\left(e^{\frac{t}{2} (\Delta_U+\frac{1}{N^2}\tilde{\Delta}_U)}P\right)\left(G^{(N)}_t\right) .$$
\end{proposition}
\begin{proof}
Recall that we denote by $U$ and $G$ the identity functions of respectively $U(N)$ and $GL_N(\mathbb{C})$.
We know that, for all $P\in \mathbb{C}\{X,X^{-1}\}$, we have $\frac{\diff}{\diff t} \left(e^{\frac{t}{2}\Delta_{U(N)}}P(U)\right)=e^{\frac{t}{2}\Delta_{U(N)}}\left(\frac{1}{2}\Delta_{U(N)}P(U)\right)$. Thus, thanks to Lemma~\ref{delta}, we have
$$\frac{\diff}{\diff t} \left(e^{\frac{t}{2}\Delta_{U(N)}}(P(U))\right)= e^{\frac{t}{2}\Delta_{U(N)}}\left(\left(\frac{1}{2}(\Delta_U+\frac{1}{N^2}\tilde{\Delta}_U)P\right)(U)\right).$$
Since the maps $(P\mapsto e^{\frac{t}{2}\Delta_{U(N)}}(P(U)) )_{t\geq 0}$ are linear, Lemma~\ref{eqdiff} tells us that, for all $P\in \mathbb{C}\{X,X^{-1}\}$ and $t\geq 0$, we have $e^{\frac{t}{2}\Delta_{U(N)}}(P(U))=(e^{\frac{t}{2} (\Delta_U+\frac{1}{N^2}\tilde{\Delta}_U)}P)(U).$

Let $P\in \mathbb{C}\{X,X^{-1}\}$. Set $f=P(U)$. We have $e^{\frac{t}{2}\Delta_{U(N)}}f=(e^{\frac{t}{2} (\Delta_U+\frac{1}{N^2}\tilde{\Delta}_U)}P)\left(U\right)$ which has an analytic continuation to a holomorphic function on $GL_N(\mathbb{C})$, denoted by $F$, and by definition $$B^{(N)}_t\left(P\left(U^{(N)}_t\right)\right)=F\left(G^{(N)}_t\right).$$
It remains to prove that $F=(e^{\frac{t}{2} (\Delta_U+\frac{1}{N^2}\tilde{\Delta}_U)}P)(G)$, or equivalently to prove that $(e^{\frac{t}{2} (\Delta_U+\frac{1}{N^2}\tilde{\Delta}_U)}P)(G)$ is holomorphic, which is evident since $e^{\frac{t}{2} (\Delta_U+\frac{1}{N^2}\tilde{\Delta}_U)}P\in \mathbb{C}\{X,X^{-1}\}$, and coordinates of $G$ and $G^{-1}$ are holomorphic.
\end{proof}

\subsection{Large-$N$ limit}

\subsubsection{Brownian motion on $GL_N(\mathbb{C})$}The free circular multiplicative Brownian motion $(G_t)_{t\geq 0}$ is defined with a free stochastic differential equation driven by a free Brownian motion (see Section~\ref{FBMonGL}). It can be proved that the Brownian motion on $GL_N(\mathbb{C})$ verifies the same stochastic differential equation driven by a Brownian motion on $M_N(\mathbb{C})$ (see~\cite{Biane1997}). The following theorem is not surprising and reflects a well-known phenomenon: the free stochastic calculus is intuitively the limit of the stochastic calculus on $M_N(\mathbb{C})$ when $N\rightarrow \infty$. See also~\cite{Kemp2013} and~\cite{Kemp2013a} for similar results with a larger class of heat kernel measures on $GL_N(\mathbb{C})$.
\begin{theorem}\label{Glim}
For all $N\in \mathbb{N}^*$, let $(G_t^{(N)})_{t\geq 0}$ be a Brownian motion on $GL_N(\mathbb{C})$ and let $(G_t)_{t\geq 0}$ be a free circular multiplicative Brownian motion. For all $n\in \mathbb{N}$, and all $P_0,\ldots,P_n\in \mathbb{C}\langle X,X^*,X^{-1},{X^*}^{-1}\rangle$, for each $t\geq 0$, as $N \rightarrow \infty$, we have
$$\mathbb{E}\left[\tr\left(P_0\left(G_t^{(N)}\right)\right)\cdots \tr\left(P_n\left(G_t^{(N)}\right)\right) \right]=\tau\left(P_0\left(G_t\right)\right)\cdots\tau\left(P_n\left(G_t\right)\right) +O(1/N^2).$$
\end{theorem}
\begin{proof}
It is equivalent to prove that, for all $P\in \mathbb{C}\{X,X^*,X^{-1},{X^*}^{-1}\}$, we have
$$\mathbb{E}\left[\tr\left(P\left(G_t^{(N)}\right)\right)\right]=\tau\left(P\left(G_t\right)\right)+O(1/N^2).$$
Indeed, Theorem~\ref{Glim} is the particular case when $P=P_0\tr P_1 \cdots \tr P_n$.

Let $P\in \mathbb{C}\{X,X^*,X^{-1},{X^*}^{-1}\}$ and $t\geq 0$. From Corollary~\ref{GLN}, we have
$$\mathbb{E}\left[\tr\left(P\left(G_t^{(N)}\right)\right)\right]=\tr\left(\left(e^{\frac{t}{4} (\Delta_{GL}+\frac{1}{N^2}\tilde{\Delta}_{GL})}P\right)(1)\right)=\left(e^{\frac{t}{4} (\Delta_{GL}+\frac{1}{N^2}\tilde{\Delta}_{GL})}P\right)(1).$$
From Proposition~\ref{distG}, we have $\tau\left(P\left(G_t\right)\right)=\left(e^{\frac{t}{4} \Delta_{GL}}P\right)(1).$

Let $d$ be the degree of $P$. Because the exponential map is differentiable, we have
$$ e^{\frac{t}{4} (\Delta_{GL}+\frac{1}{N^2}\tilde{\Delta}_{GL})}=e^{\frac{t}{4} \Delta_{GL}}+O(1/N^2)$$
in the finite dimensional space $\End(\mathbb{C}_d\{X,X^*,X^{-1},{X^*}^{-1}\})$. Since $A\mapsto (A(P))(1)$ is a linear map from $\End(\mathbb{C}_d\{X,X^*,X^{-1},{X^*}^{-1}\})$ to $\mathbb{C}$, it therefore is continuous and we have
$$ \left(e^{\frac{t}{4} (\Delta_{GL}+\frac{1}{N^2}\tilde{\Delta}_{GL})}P\right)(1)=\left(e^{\frac{t}{4} \Delta_{GL}}P\right)(1)+O(1/N^2), $$
which completes the proof.
\end{proof}

\subsubsection{Hall transform}
We recall that, for all $t> 0$, $\mathcal{G}_t$ is the free Hall transform at the level of Laurent polynomials $\mathbb{C}[X,X^{-1}]$, and that $e^{\frac{t}{2} \Delta_U}$ is the free Hall transform at the level of the space $\mathbb{C}\{X,X^{-1}\}$ (see Section~\ref{hallfree}). The non-commutative random variables $U^{(N)}_t$ and $G^{(N)}_t$ give approximations to the free unitary Brownian motion $U_t$ and the circular multiplicative Brownian motion $G_t$ at time $t$. The following result states that the Hall transform also goes to the limit as $N$ tends to $\infty$, on both the $\mathbb{C}[X,X^{-1}]$ and the $\mathbb{C}\{X,X^{-1}\}$-calculus.\label{Laurent} See also~\cite{Driver2013} for another proof of the first item.

\begin{theorem}\label{Flim}Let $t> 0$. For all $N\in \mathbb{N}^*$, let $(U_t^{(N)})_{t\geq 0}$ be a Brownian motion on $U(N)$, and $(G_t^{(N)})_{t\geq 0}$ be a Brownian motion on $GL_N(\mathbb{C})$. We have
\begin{enumerate}
\item for all Laurent polynomial $P\in \mathbb{C}[X,X^{-1}]$, as $N \rightarrow \infty$,
$$ \left\|B^{(N)}_t\left(P\left(U^{(N)}_t\right)\right)-\mathcal{G}_t(P)\left(G^{(N)}_t\right)\right\|^2_{L^2_{\rm{hol}}(G^{(N)}_t)\otimes M_N(\mathbb{C})}=O(1/N^2),$$

\item for all $P\in \mathbb{C}\{X,X^{-1}\}$, as $N \rightarrow \infty$, $$ \left\|B^{(N)}_t\left(P\left(U^{(N)}_t\right)\right)-\left(e^{\frac{t}{2} \Delta_U}P\right)\left(G^{(N)}_t\right)\right\|^2_{L^2_{\rm{hol}}(G^{(N)}_t)\otimes M_N(\mathbb{C})}=O(1/N^2).$$
\end{enumerate}

\end{theorem}
\begin{proof}
Let $t> 0$ and $P\in \mathbb{C}[X,X^{-1}]$. Let $N\in \mathbb{N}^*$. From Corollary~\ref{GLN} and Proposition~\ref{btpun}, we have
\begin{eqnarray*}
& &\hspace{-1cm}\left\|B^{(N)}_t\left(P\left(U^{(N)}_t\right)\right)-\mathcal{G}_t(P)\left(G^{(N)}_t\right)\right\|^2_{L^2_{\rm{hol}}(G^{(N)}_t)\otimes M_N(\mathbb{C})}\\
&=&\left\|e^{\frac{t}{2} (\Delta_U+\frac{1}{N^2}\tilde{\Delta}_U)}P\left(G^{(N)}_t\right)-\mathcal{G}_t(P)\left(G^{(N)}_t\right)\right\|^2_{L^2_{\rm{hol}}(G^{(N)}_t)\otimes M_N(\mathbb{C})}\\
&=&e^{\frac{t}{4} (\Delta_{GL}+\frac{1}{N^2}\tilde{\Delta}_{GL})}\left(\left(e^{\frac{t}{2} (\Delta_U+\frac{1}{N^2}\tilde{\Delta}_U)}P-\mathcal{G}_t(P)\right)^*\left(e^{\frac{t}{2} (\Delta_U+\frac{1}{N^2}\tilde{\Delta}_U)}P-\mathcal{G}_t(P)\right)\right)(1).
\end{eqnarray*}
Let $d$ be the degree of $P$. We will work in finite-dimensional spaces and consequently all the norms are equivalent and the linearity of a map implies its boundedness. The differentiability of the exponential map leads to
$$ e^{\frac{t}{4} (\Delta_{GL}+\frac{1}{N^2}\tilde{\Delta}_{GL})}=e^{\frac{t}{4} \Delta_{GL}}+O(1/N^2)\ \text{ and }e^{\frac{t}{2} (\Delta_{U}+\frac{1}{N^2}\tilde{\Delta}_{U})}=e^{\frac{t}{2} \Delta_{U}}+O(1/N^2)$$
in the finite dimensional space $\End(\mathbb{C}_d\{X,X^*,X^{-1},{X^*}^{-1}\})$. Since $A\mapsto A(P)$ is linear, we deduce that
$$\left(e^{\frac{t}{2} (\Delta_U+\frac{1}{N^2}\tilde{\Delta}_U)}P-\mathcal{G}_t(P)\right)=\left(e^{\frac{t}{2} \Delta_U}P-\mathcal{G}_t(P)\right)+O(1/N^2).$$
The linearity of $(A, P, Q)\mapsto (A(P^*Q))(1)$ from $\End(\mathbb{C}_d\{X,X^*,X^{-1},{X^*}^{-1}\})\times \mathbb{C}_d\{X,X^*,X^{-1},{X^*}^{-1}\}^2$ to $ \mathbb{C}$ implies that
\begin{eqnarray*}& &\hspace{-1cm}e^{\frac{t}{4} (\Delta_{GL}+\frac{1}{N^2}\tilde{\Delta}_{GL})}\left(\left(e^{\frac{t}{2} (\Delta_U+\frac{1}{N^2}\tilde{\Delta}_U)}P-\mathcal{G}_t(P)\right)^*\left(e^{\frac{t}{2} (\Delta_U+\frac{1}{N^2}\tilde{\Delta}_U)}P-\mathcal{G}_t(P)\right)\right)(1)\\
&=&e^{\frac{t}{4} \Delta_{GL}}\left(\left(e^{\frac{t}{2} \Delta_U}P-\mathcal{G}_t(P)\right)^*\left(e^{\frac{t}{2} \Delta_U}P-\mathcal{G}_t(P)\right)\right)(1)+O(1/N^2).
\end{eqnarray*}
From Proposition~\ref{distG}, we have
\begin{eqnarray*}& &\hspace{-1cm}e^{\frac{t}{4} \Delta_{GL}}\left(\left(e^{\frac{t}{2} \Delta_U}P-\mathcal{G}_t(P)\right)^*\left(e^{\frac{t}{2} \Delta_U}P-\mathcal{G}_t(P)\right)\right)(1)\\
&=&\tau\left(\left((e^{\frac{t}{2} \Delta_U}P)(G_t)-\mathcal{G}_t(P)(G_t)\right)^*\left((e^{\frac{t}{2} \Delta_U}P)(G_t)-\mathcal{G}_t(P)(G_t)\right)\right),
\end{eqnarray*}
and by Theorem~\ref{fht}, we have $e^{\frac{t}{2} \Delta_U}P(G_t)-\mathcal{G}_t\Big(P\Big)(G_t)=0$.
Finally, we have
\begin{equation*}
 \left\|B^{(N)}_t\left(P\left(U^{(N)}_t\right)\right)-\mathcal{G}_t(P)\left(G^{(N)}_t\right)\right\|^2_{L^2_{\rm{hol}}(G^{(N)}_t)\otimes M_N(\mathbb{C})}=0+O(1/N^2).
\end{equation*}
The second assertion follows using the same argument substituting $e^{\frac{t}{2} \Delta_U}P$ for $\mathcal{G}_t(P)$.
\end{proof}

\addtocontents{toc}{\protect\setcounter{tocdepth}{1}}
\begin{appendix}

\section*{Appendix. The construction of $\mathbb{C}\{X_i:i\in I\}$}

In this appendix, the algebra $\mathbb{C}\{X_i:i\in I\}$ is constructed as the direct limit of the inductive system formed by the family of finite tensor product $\mathbb{C}\left\langle X_i:i\in I\right\rangle\otimes  \mathbb{C}\left\langle X_i:i\in I\right\rangle^{\odot n}$ where $\odot$ is defined below. Then we verify that, when endowed with the appropriate product and center-valued expectation, $\mathbb{C}\{X_i:i\in I\}$ satisfies the universal property~\ref{up}.

\subsection*{Construction by direct limit}\label{consdeux}

Let $n\in \mathbb{N}$. Let $E$ be the subspace of $\mathbb{C}\left\langle X_i:i\in I\right\rangle\otimes \mathbb{C}\left\langle X_i:i\in I\right\rangle^{\otimes n}$ generated by the families
$$\{ P_0\otimes P_1\otimes \cdots \otimes P_n-P_0\otimes P_{\sigma(1)}\otimes \cdots \otimes P_{\sigma(n)}: P_0,\ldots,P_n\in \mathbb{C}\left\langle X_i:i\in I\right\rangle,\sigma\in \mathfrak{S}_n\}. $$
Let us denote by $\mathbb{C}\left\langle X_i:i\in I\right\rangle\otimes  \mathbb{C}\left\langle X_i:i\in I\right\rangle^{\odot n}$ the quotient space of $\mathbb{C}\left\langle X_i:i\in I\right\rangle\otimes \mathbb{C}\left\langle X_i:i\in I\right\rangle^{\otimes n}$ by its subspace $E$.
If $P_0,\ldots,P_n\in \mathbb{C}\left\langle X_i:i\in I\right\rangle$, we denote the equivalence class of $P_0\otimes P_1\otimes \cdots \otimes P_n$ in $\mathbb{C}\left\langle X_i:i\in I\right\rangle\otimes  \mathbb{C}\left\langle X_i:i\in I\right\rangle^{\odot n}$ by $P_0 \otimes  P_1 \odot \cdots \odot  P_n$.

For all $n\leq m \in \mathbb{N}$, we identify $\mathbb{C}\left\langle X_i:i\in I\right\rangle\otimes  \mathbb{C}\left\langle X_i:i\in I\right\rangle^{\odot n}$ as a subset of $\mathbb{C}\left\langle X_i:i\in I\right\rangle\otimes  \mathbb{C}\left\langle X_i:i\in I\right\rangle^{\odot m}$ thanks to the injective morphism
$$\begin{array}{lrcl}
i_{n,m}: & \mathbb{C}\left\langle X_i:i\in I\right\rangle\otimes  \mathbb{C}\left\langle X_i:i\in I\right\rangle^{\odot n}& \longrightarrow & \mathbb{C}\left\langle X_i:i\in I\right\rangle\otimes  \mathbb{C}\left\langle X_i:i\in I\right\rangle^{\odot m} \\
    & P_0\otimes P_1\odot \cdots \odot P_n & \longmapsto & P_0\otimes P_1\odot \cdots \odot P_n\odot 1 \odot \cdots \odot 1.
\end{array}$$
The inclusion maps $\left(i_{n,m}\right)_{n\leq m \in \mathbb{N}}$ make $\mathbb{C}\left\langle X_i:i\in I\right\rangle\otimes  \mathbb{C}\left\langle X_i:i\in I\right\rangle^{\odot n}$ an inductive system in the sense that, for all $l \leq n\leq m \in \mathbb{N}$, we have $i_{l,m}=i_{l,n}\circ i_{n,m}$.
The vector space $\mathbb{C}\{X_i:i\in I\}$ is the direct limit $\mathbb{C}\{X_i:i\in I\}=\varinjlim\mathbb{C}\left\langle X_i:i\in I\right\rangle\otimes  \mathbb{C}\left\langle X_i:i\in I\right\rangle^{\odot n}$. It is the minimal vector space which contains all the tensor products $\mathbb{C}\left\langle X_i:i\in I\right\rangle\otimes  \mathbb{C}\left\langle X_i:i\in I\right\rangle^{\odot n}$ as subspaces identified via the inclusion maps $\left(i_{n,m}\right)_{n\leq m \in \mathbb{N}}$. In particular, the space $\mathbb{C}\left\langle X_i:i\in I\right\rangle$ is a subspace of $\mathbb{C}\{X_i:i\in I\}$.

We have at our disposal a basis of $\mathbb{C}\{X_i:i\in I\}$. Indeed, for all $ n\in \mathbb{N}$, $$\{M_0\otimes  M_1 \odot \cdots \odot M_n:  M_0,\ldots M_n \text{ are monomials of }\mathbb{C}\left\langle X_i:i\in I\right\rangle\}$$ is a basis of $\mathbb{C}\left\langle X_i:i\in I\right\rangle\otimes  \mathbb{C}\left\langle X_i:i\in I\right\rangle^{\odot n}$. Thus, $$\{M_0\otimes  M_1 \odot \cdots \odot M_n: n\in \mathbb{N}, M_0,\ldots M_n \text{ are monomials of }\mathbb{C}\left\langle X_i:i\in I\right\rangle\}$$is a basis of $\mathbb{C}\{X_i:i\in I\}$, called the canonical basis.

Let us define now the product and the center-valued expectation of $\mathbb{C}\{X_i:i\in I\}$. For all $n\in \mathbb{N}$, $M_0,\ldots M_n$ monomials of $\mathbb{C}\langle X_i:i\in I\rangle$, and  $ n'\in \mathbb{N}$, $N_0,\ldots N_{n'}$ monomials of $\mathbb{C}\left\langle X_i:i\in I\right\rangle$, let us define the product $M\cdot N=MN$ of $M=M_0\otimes  M_1 \odot \cdots \odot M_n$ and $N=N_0\otimes  N_1 \odot \cdots \odot N_{n'}$ by
$$M \cdot N=MN=M_0 N_0 \otimes M_1 \odot \cdots M_n \odot N_1 \odot \cdots \odot  N_{n'}.$$
We extend this product to all $\mathbb{C}\{X_i:i\in I\}$ by linearity. Endowed with this product, $\left(\mathbb{C}\{X_i:i\in I\},\cdot \right)$ is a unital complex algebra.

Let $n\in \mathbb{N}$, and $M_0,\ldots,M_n\in  \mathbb{C}\left\langle X_i:i\in I\right\rangle$. For all $M=M_0\otimes  M_1 \odot \cdots \odot M_n\in \mathbb{C}\{X_i:i\in I\}$, we set
$$\tr M=1 \otimes M_0\odot  M_1 \odot \cdots \odot M_n,$$
and we extend the map $\tr$ to all $\mathbb{C}\{X_i:i\in I\}$ by linearity. Defined this way, $\tr$ is a center-valued expectation. Indeed, the center consists precisely of linear combinations of elements of the form $1 \otimes M_0\odot  M_1 \odot \cdots \odot M_n$.

We remark that, for all $n\in \mathbb{N}, P_0,\ldots P_n \in \mathbb{C}\left\langle X_i:i\in I\right\rangle$, we have $P_0\otimes  P_1 \odot \cdots \odot P_n=P_0 \tr P_1 \cdots \tr  P_n.$ Thus, the canonical basis of $\mathbb{C}\{X_i:i\in I\}$ can be rewritten
$$\{M_0 \tr M_1 \cdots \tr  M_n: n\in \mathbb{N}, M_0,\ldots M_n \text{ are monomials of }\mathbb{C}\left\langle X_i:i\in I\right\rangle\}.$$

\subsection*{Universal property}\label{constrois}
The triplet $\left(\mathbb{C}\{X_i:i\in I\},\tr,\left(X_i\right)_{i\in I}\right)$ is now constructed, and it remains to prove the universal property.

Let $\mathcal{A}$ be an algebra endowed with a center-valued expectation $\tau$, and with a family of $I$ elements $\mathbf{A}=\left(A_i\right)_{i\in I}$. Let $n\in \mathbb{N}$, and $M_0,\ldots M_n$ be monomials of $\mathbb{C}\left\langle X_i:i\in I\right\rangle$. Let us define $f\left(M_0 \tr M_1 \cdots \tr  M_n \right)$ by
$$f\left(M_0\tr M_1 \cdots \tr M_n\right)=\tau\left(M_1\left(\mathbf{A}\right)\right)\cdots\tau\left(M_n\left(\mathbf{A}\right)\right) \cdot M_0\left(\mathbf{A}\right),$$
and we extend $f$ to a map from $\mathbb{C}\{X_i:i\in I\}$ to $\mathcal{A}$ by linearity. The properties of polynomial calculus make $f$ an algebra morphism. Moreover, by definition, for all $i\in I$, we have $f(X_i)=A_i$. Finally, we verify that, for all $P\in \mathbb{C}\{ X_i:i\in I\}$, we have $\tau(f(P))=f(\tr (P))$. Let $n\in \mathbb{N}$, and $P_0,\ldots,P_n\in  \mathbb{C}\left\langle X_i:i\in I\right\rangle$. We have
\begin{eqnarray*}
\tau \left(f\left(P_0\tr P_1 \cdots \tr P_n\right)\right)&=&\tau \left(\tau\left(P_1\left(\mathbf{A}\right)\right)\cdots\tau\left(P_n\left(\mathbf{A}\right)\right) \cdot P_0\left(\mathbf{A}\right)\right)\\
&=&\tau\left(P_0\left(\mathbf{A}\right)\right)\tau\left(P_1\left(\mathbf{A}\right)\right)\cdots\tau\left(P_n\left(\mathbf{A}\right)\right)\\
&=&f\left(1\tr P_0 \cdots \tr P_n\right)\\
&=&f\left(\tr\left(P_0\tr P_1\cdots \tr P_n\right)\right)
\end{eqnarray*}
and by linearity, we deduce that, for all $P\in\mathbb{C}\{ X_i:i\in I\}$, we have $\tau \left(f(P)\right)=f\left(\tr\left(P\right)\right)$.

Thus, there exists an algebra homomorphism $f$ from $\mathbb{C}\left\langle X_i:i\in I\right\rangle$ to $\mathcal{A}$ such that
\begin{enumerate}
\item for all $i\in I$, we have $f(X_i)=A_i$;
\item for all $X\in \mathbb{C}\{ X_i:i\in I\}$, we have $\tau \left(f(X)\right)=f\left(\tr\left(X\right)\right)$.
\end{enumerate}
The uniqueness of such a morphism is clear since the action of $f$ on $\mathbb{C}\left\langle X_i:i\in I\right\rangle$ is uniquely determined by the polynomial calculus, and by consequence $f$ is uniquely determined on the basis $$\{M_0 \tr M_1 \cdots \tr  M_n: n\in \mathbb{N}, M_0,\ldots M_n \text{ are monomials of }\mathbb{C}\left\langle X_i:i\in I\right\rangle\}.$$Indeed, let $n\in \mathbb{N}$, and $M_0,\ldots M_n$ be monomials of $\mathbb{C}\left\langle X_i:i\in I\right\rangle$. We have
\begin{eqnarray*}
f\left(M_0\tr M_1\cdots\tr M_n\right)&=&f(M_0)f\left(\tr M_1\right)\cdots f\left(\tr M_n\right)\\
&=&f(M_0)\tau\left(f\left( M_1\right)\right)\cdots \tau\left(f\left( M_n\right)\right)\\
&=&\tau\left(M_1\left(\mathbf{A}\right)\right)\cdots\tau\left(M_n\left(\mathbf{A}\right)\right) \cdot M_0\left(\mathbf{A}\right).
\end{eqnarray*}

\end{appendix}

\subsection*{Acknowledgments}
The author is grateful to his PhD advisor Thierry L\'{e}vy for his helpful comments and corrections which led to improvements in this manuscript, and to Philippe Biane for clarifying to him the definition of the free Hall transform. Many thanks are due to Antoine Dahlqvist and Franck Gabriel for useful discussions related to this work.

\nocite{*}

\bibliographystyle{plain}
\bibliography{/Users/guillaumecebron/Documents/Bibtex/freeconv}

\end{document}